\RequirePackage{fix-cm}
\documentclass{article}     
\usepackage{graphicx}
\usepackage{amsmath}
\usepackage{amssymb}
\usepackage{eurosym}
\usepackage{float}
\usepackage{multirow}
\usepackage{amsfonts}
\usepackage{bm}
\usepackage[margin=1in]{geometry}
\newcommand{\ie}{\emph{i.e.,} }
\newcommand{\eg}{\emph{e.g.,} }

\graphicspath{{Figures/}{./}} 

\begin{document}

\title{ Bayesian optimization of variable-size design space problems}


\author{Julien Pelamatti, Lo\"ic Brevault, Mathieu Balesdent, El-Ghazali Talbi, Yannick Guerin
}


\maketitle

\begin{abstract}
Within the framework of complex system design, it is often necessary to solve mixed variable optimization problems, in which the objective and constraint functions can depend simultaneously on continuous and discrete variables. Additionally, complex system design problems occasionally present a variable-size design space. This results in an optimization problem for which the search space varies dynamically (with respect to both number and type of variables) along the optimization process as a function of the values of specific discrete decision variables. Similarly, the number and type of constraints can vary as well. In this paper, two alternative Bayesian Optimization-based approaches are proposed in order to solve this type of optimization problems. The first one consists in a budget allocation strategy allowing to focus the computational budget on the most promising design sub-spaces. The second approach, instead, is based on the definition of a kernel function allowing to compute the covariance between samples characterized by partially different sets of variables. The results obtained on analytical and engineering related test-cases show a faster and more consistent convergence of both proposed methods with respect to the standard approaches.
 
\end{abstract}

\section{Introduction}
\label{intro}
\subsection{Context}
The design of complex systems, such as launch vehicles, aircraft, automotive vehicles or electronic components, can usually be represented under the form of an optimization problem. In other words, for a given formulation of the problem in terms of objective and constraint functions, as well as design variables these functions depend on, the aim of the design process is to determine the values of the design variables which yield the best value of the objective function. Furthermore, this optimal solution must also comply with all the constraints the problem is subject to. In the most simple case, the architecture of the considered system is determined beforehand through an empirical process. As a consequence, the design variables characterizing the resulting problem are usually defined within a continuous search space. For illustrative purposes, the design of Re-usable Launch Vehicles (RLV) is used as an example throughout this paper. Within this design framework, typical continuous design variables are sizing parameters, combustion pressures and propellant masses. However, in most real world engineering design problems, the architecture of the considered system cannot be specified beforehand and must typically be defined during the early stages of the design process. In this case, it is therefore necessary to simultaneously optimize the architecture layout and the continuous design variables which characterize it. The most straightforward approach in order to perform this simultaneous optimization consists in iterating over every possible architecture definition and performing an optimization of the continuous design variables for each architecture \cite{prasadh2014systematic}. However, depending on the complexity of the considered system, it may be necessary to consider several thousands of different architectures \cite{Frank2016AVehicles}, thus resulting extremely time consuming, if not unfeasible.

Without loss of generality, the choices related to the architecture definition can be represented under the form of discrete design variables, sometimes referred to as categorical or qualitative variables, which characterize design alternatives and/or technological choices. Typical examples of technological choices which can be encountered within the context of RLV design are the type of material to use for a given sub-system, the type of propulsion to include in a given stage, the presence of a certain system component and the number of reinforcements to be included in a given structure. From an analytical perspective, part of these technological choices plays the role of standard discrete variables defined within a finite set of choices, whereas others play a more complex role, as they can also influence the definition of the objective and constraint functions, as well as the number and type of variables that characterize the problem. These particular choices are often related to the selection of sub-problem specific technologies. For instance, depending on whether a combustion, hybrid or electric engine is selected for the design of a car, the specific engine related design variables can vary considerably (\eg combustion chamber pressure, electric engine torque). Additionally, the associated constraint can also differ (\eg combustion temperature associated constraints).  Due to the inherent discrete and potentially non-numerical nature of these design variables, the concept of metrics is usually not definable within their domain, thus resulting in an unordered set of possible choices. In the literature, this kind of problems are referred to as mixed-variable optimization problems \cite{Lucidi2005AnProgramming} or \textbf{Variable-Size Design Space Problem} \cite{Nyew2015} (VSDSP). Most modern optimization algorithms are developed with the purpose of solving design problems essentially characterized  by continuous and integer variables and by consequence, the introduction of these categorical (qualitative) variables raises a number of additional challenges. Firstly, a large number of algorithms cannot be used due to the non-differentiability of the objective function and the possible absence of metrics in the discrete variables domains (\emph{e.g.,} gradient-based algorithms or Nelder-Mead simplex \cite{Nelder1964}). Further limitations arise from the need to initialize newly created variables during the optimization process in a dynamically varying design space (\eg Genetic algorithms \cite{Goldberg1989}) and the impossibility to relax the integrity constraint on some of the discrete variables (\emph{e.g.,} Branch and Bound \cite{Land2010}).  Additionally, for problems characterized by large numbers of technological choices and options, the combinatorial size of the discrete variables search space can be considerably large, thus rendering a complete exploration of the problem search space computationally particularly difficult. Finally, real-world complex system design problems are often characterized by computationally intensive objective and/or constraint functions (\eg Finite Element Models), which require any applicable algorithm to converge towards the optimum of the considered problem with a limited number of evaluations.

\subsection{Problem definition}
Without loss of generality, a generic variable-size design space optimization problem can be modeled as depending on three different types of design variables: continuous, discrete and dimensional.
\begin{itemize}
\item \textbf{Continuous variables}: $\mathbf{x} $ \\
Continuous variables refer to real numbers defined within a given interval.
\item \textbf{Discrete variables}: $\mathbf{z}$ \\
Discrete variables are non-relaxable variables defined within a finite set of choices.  Discrete variables are typically divided into 2 categories, quantitative and qualitative, depending on whether a relation of order between the possible values of a given variable can be defined.  For clarity and synthesis purposes, no distinction between quantitative and qualitative variables is made in this paper, and both types of variables are simply referred to as discrete variables.

\item   \textbf{Dimensional variables}:  $\mathbf{w} $ \\
Similarly to the discrete variables, dimensional variables are non-relaxable variables defined within a finite set of choices. The main distinction is that, depending on their values, the number and type of continuous and discrete variables the problem functions depend on can vary. Furthermore, they can also influence the number and type of constraints a given candidate solution is subject to. These particular variables often represent the choice between several possible sub-system architectures. Each architecture is usually characterized by a partially different set of design variables, and therefore, depending on the considered choice, different continuous and discrete design variables must be optimized. 
\end{itemize}
Although both discrete and dimensional variables may lack a conceptual numerical representation (\emph{i.e.,} they are not associated to measurable values), it is common practice to assign an integer value to every considered alternative of the given variable in order to provide a clear problem formulation as well as simplifying the implementation process \cite{FilomenoCoelho2014}. For instance, if 3 hypothetical choices for the type of material are considered (\emph{e.g.,} aluminum, steel, titanium), the associated discrete design variable can be defined as having 3 possible values: 0, 1 and 2. In this work, it is assumed that none of the discrete and dimensional variables are directly related to a measurable value and the integer values associated to the possible discrete choices are assigned arbitrarily. However, it is important to highlight the fact that said choice has no influence on working principle of the various algorithms presented in the following sections. Similarly to what is proposed by Lucidi \emph{et al.} \cite{Lucidi2005AnProgramming}, a generic variable-size design space problem can be formulated as follows:
\begin{eqnarray}
\label{PBDef}
\min &\;& f(\mathbf{x},\mathbf{z},\mathbf{w}) \; \ \qquad  \quad f: \mathbb{R}^{n_{x}(\mathbf{w})} \times \prod_{d = 1}^{n_{z}(\mathbf{w})} F_{z_d} \times F_w  \rightarrow F_f \subseteq \mathbb{R} \\ \nonumber
\text{w.r.t.} &\;& \mathbf{x} \in F_x(\mathbf{w}) \subseteq \mathbb{R}^{n_{x}(\mathbf{w})}  \\ \nonumber
&& \mathbf{z}  \in \prod_{d = 1}^{n_{z}(\mathbf{w})} F_{z_d} \\ \nonumber
&& \mathbf{w}  \in F_w  \\ \nonumber
\text{s.t.} &\;& \mathbf{g}  (\mathbf{x},\mathbf{z},\mathbf{w}) \leq  0 \\ \nonumber
&& g_i : F_{x_i}(\mathbf{w}) \times \prod_{d = 1}^{n_{z_i}(\mathbf{w})} F_{z_{d_i}} \times F_w \rightarrow F_{g_i} \subseteq \mathbb{R} \qquad \mbox{for} \qquad i = 1,...,n_g(\mathbf{w})
\end{eqnarray}
where $f(\cdot)$ and $\mathbf{g}(\cdot)$ are respectively the objective function and constraint vector, $\mathbf{x}$ is a $n_x(\mathbf{w})$-dimensional vector containing the continuous design variables, $\mathbf{z}$ is a $n_z(\mathbf{w})$-dimensional vector containing the discrete design variables, $\mathbf{w}$ is a $n_w$-dimensional vector containing the dimensional design variables and $n_g (\mathbf{w})$ is the number of constraints the problem is subject to. For the sake of clarity, the discrete variables domain is referred to as $F_z(\mathbf{w})$ in the remainder of the paper:
\begin{equation}
F_z(\mathbf{w}) = \prod_{d = 1}^{n_{z}(\mathbf{w})} F_{z_d}
\end{equation}
Similarly, the variable domains relative to the constraint $g_i$ are referred to as $ F_{z_i}(\mathbf{w})$ and $F_{x_i}(\mathbf{w})$. Finally, it is important to highlight that throughout this work, only single-objective optimization is considered. 
\\[12pt]
As can be noticed, the continuous and discrete search spaces of  Eq. \ref{PBDef} are not fixed throughout the optimization process, as they depend on the values of the dimensional variables. As a result, two different candidate solutions of the optimization problem defined above can be characterized by a partially different set of variables, depending on the values of their dimensional variables. In the same way, the feasibility domain can also vary according to the values of $\mathbf{w}$, which results in different candidate solutions possibly being subject to different constraints (in terms of both type and number of active constraints). Furthermore, it is important to notice that the constraint functions can also present a search space which can vary as a function of the dimensional variable values. By consequence, the optimization problem defined in Eq. \ref{PBDef} varies dynamically along the optimum search, both in terms of number and type of design variables as well as feasibility domain, depending on the values of the candidate solution dimensional variables. For the sake of clarity, the possibility of having the dimensional variables search space vary as a function of the dimensional variables themselves is not taken into consideration. Although this choice does not limit the applicability of the methods presented in this paper it could, however, require a more complex problem formulation. It is essential to highlight the fact that although the design space may vary depending on the dimensional variable values, the quantity represented by the objective function.
\\[12pt]
For the sake of conciseness, in the remainder of this work the following two terms are used:
\begin{itemize}
\item \textbf{levels}: the possible values of a discrete ($z \in \{ z_1,\dots,z_l\}$) or a dimensional ($w \in \{w_1,\dots,w_l\}$) variable, where $l$ indicates the total number of levels associated to a given variable

\item \textbf{categories}: possible combinations of levels. In other words, a category characterizes a candidate solution in the combinatorial discrete/dimensional search space.
\end{itemize}  
As an illustrative example, 3 hypothetical discrete variables associated to the RLV design can be considered:
\[
\mathbf{z} = \{z_1,z_2,z_3\}: \begin{cases}
z_1 \in \{ \mbox{aluminum}, \mbox{steel}, \mbox{titanium}\} \\
z_2  \in \{\mbox{liquid propulsion}, \mbox{solid propulsion}\} \\
z_3 \in \{1 \mbox{ engine}, 2 \mbox{ engines}, 4 \mbox{ engines},8 \mbox{ engines} \}
\end{cases}
\]
In this case, the choice \emph{aluminum} is one of the levels of the variable $z_1$, and the 3 discrete variables are respectively characterized by 3, 2 and 4 levels. A possible category of the problem above can, for instance, be characterized by \{steel, liquid propulsion, 8 engines\}. Let $l_d$ be the number of levels associated to the discrete variable $z_d$, the total number of categories $m$ associated to a given problem can be computed as:
\begin{equation}
m = \prod_{d = 1}^{n_z} l_d
\end{equation}
For instance, the number of categories associated to the illustrative example above is equal to $m = 3 *  2 * 4 = 24$.

Although VSDSP represent the majority of actual complex system design problems, very few optimization methods or algorithms allowing to solve such problems in their generic formulation exist. In most cases, the chosen approach in order to solve VSDSP consists in decomposing the global problem into several fixed-size sub-problems (\ie one per combination of dimensional variable values) and separately optimizing each one of them. Said approach is straightforward and easily implementable. However, it also tends to quickly become computationally overly demanding when confronted with VSDSP characterized by a large number of sub-problems, especially in the presence of computationally intensive objective and/or constraint functions. Alternatively, a few algorithms have been extended in order to deal with VSDSP. A comprehensive taxonomy of the existing approaches and algorithms allowing to deal with this kind of problems in both complete and simplified formulations can be found in \cite{Pelamatti2018}. In \cite{Abdelkhalik2013a} a so-called hidden gene adaptation of GA is proposed. In this  algorithm, each candidate solution is represented by a chromosome containing the entirety of genes that can characterize the problem at hand. However, not all the genes are taken into account when computing the value of the objective and constraint functions. The choice regarding which genes are considered and which genes are 'hidden' depends on the values of a limited number of so-called activation genes. This variant of GA has the advantage of being intuitive and easily implementable. However, cross-overs and mutations over hidden (\emph{i.e.,} with no influence on the problem functions) parts of the vector result in ineffective numerical operations, thus wasting computational effort. Furthermore, for problems characterized by a large number of discrete categories, the vector containing the entirety of possibly present variables might become considerably large and therefore inefficient memory-wise. In the same way, the authors also implemented the hidden gene approach within DE in \cite{Abdelkhalik2013}. A more complex, but theoretically more efficient adaptation of GA called  Structured-Chromosome Evolutionary Algorithm is proposed by H.M. Nyew \emph{et al.} in \cite{Nyew2015}. In this algorithm, the candidate solution is conceptually represented by a hierarchical multi-level chromosome structure rather than a linear one. Differently than in the standard formulation of GA, the genes of each chromosome are linked by both vicinity and hierarchy relationships. Compared to the hidden genes approach, this solution has the advantage of not performing computationally wasteful mutations and cross-overs. However, the encoding of the individual chromosomes is much more complex and often requires problem-specific knowledge in order to be efficiently implemented.

A second family of algorithms which have been extended in order to provide a solution to VSDSP are mesh-based optimization algorithms, as is discussed in \cite{Lucidi2005AnProgramming}, \cite{Abramson2009MeshOptimization}, \cite{Abramson2007FilterProblems} and \cite{Audet2000PatternProgramming}. Although each paper describes a different algorithm, they share the same approach to the optimization problem which consists in an alternation between an optional user-defined search phase and a poll phase, in which an exhaustive and schematic optimum search over a mesh is performed. When applying this family of algorithms to VSDSP, it is necessary to take into account the appearance and disappearance of design variables as a function of the dimensional variables characterizing the incumbent solution around which the mesh is centered. In order to solve this issue, the mesh algorithms mentioned above alternate between searches over a mesh defined in the dimensional design space and searches over a mesh defined in the continuous and discrete search space dependent on the dimensional variables values. The application of this algorithm requires the user to provide the definition of the neighborhood of a candidate solution in the dimensional design space in order to create the mesh. 

Both population and mesh-based algorithms for the solution of VSDSP tend to require a large number of function evaluations in order to converge to the considered problem optimum. Furthermore, they also rely on a penalization-based handling of constraints, which often results inefficient if the penalization weights are not properly tuned, especially in the presence of a large number of constraints. For these reasons, the possibility of extending Bayesian Optimization (BO) algorithms in order to deal with VSDSP is discussed in this paper. More specifically, two alternative BO-based approaches allowing to solve VSDSP with a lower number of function evaluations with respect to the existing alternatives are proposed.  The first one is based on the separate and independent optimization of every sub-problem characterized by a fixed-size design space coupled with a budget allocation strategy relying on the information provided by the surrogate models of the various sub-problems functions. The second method, instead, is based on the definition of a Gaussian Process kernel defined in the variable-size design space, allowing to compute the covariance between data samples characterized by partially different sets of variables. 

Following this introduction, in the second Section a theoretical overview of Bayesian Optimization is provided and the extension allowing to include discrete variables is discussed. In Section 3, the budget allocation strategy is presented. In Section 4, two alternative definitions of a Gaussian Process kernel allowing to compute the covariance between data samples belonging to partially different design spaces are proposed. Subsequently, in Section 5 the two aforementioned methods for the optimization of VSDSP are applied on both analytical and engineering related test-cases and the obtained results are presented and discussed. Finally, in Section 6 the relevant conclusions which can be drawn from the obtained results are provided and possible perspectives and improvements are discussed.

\section{Mixed-variable Bayesian optimization}
\label{MVBO}
\begin{figure}[h]
  \centering
  \includegraphics[width=.9\linewidth]{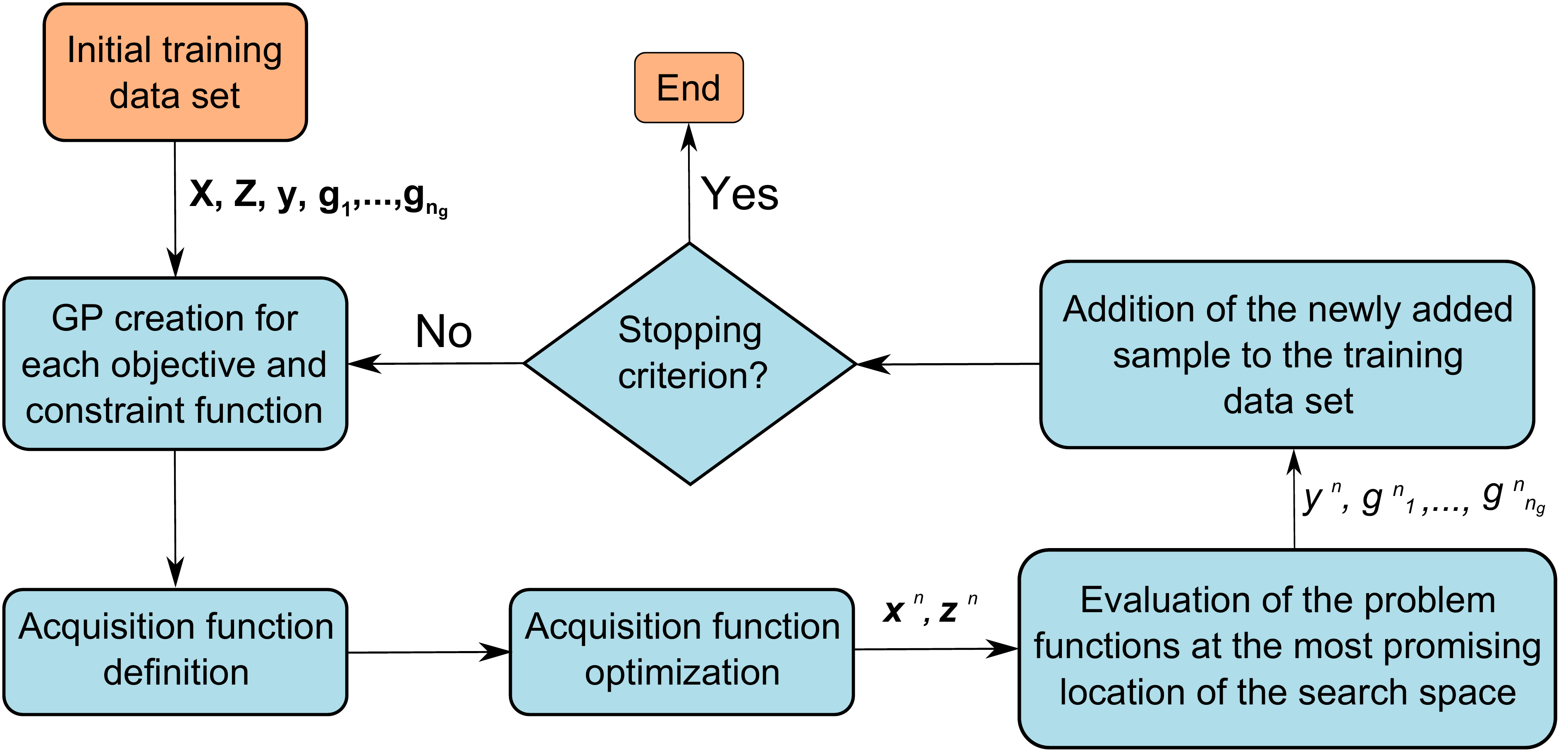}
   \caption{Schematic representation of the working principle of Bayesian Optimization. $\mathbf{X}$ and $\mathbf{Z}$ contain the continuous and discrete variable data sets, while $\mathbf{y}$ and $\textbf{g}_1,\dots,\textbf{g}_{n_g}$ contain the associated objective function and constraint responses.}
  \label{BO_framework}
\end{figure}

Bayesian Optimization \cite{Jones1998} is a Surrogate-Model Based Design Optimization (SMBDO) algorithm which relies on the use of Gaussian Process (GP) surrogate modeling \cite{Rasmussen2006}. Without loss of generality, BO can be decomposed into 2 main phases. The first phase consists in creating a separate and independent GP based surrogate model of the objective function as well as each constraint by relying on a finite training data set. During the second phase, the most promising data samples in terms of objective function value and feasibility, are identified, evaluated and added to the data set with the purpose of simultaneously refining  the surrogate model (\emph{i.e.,} improving the modeling accuracy) and exploring the areas of the design space which are more likely to contain the optimization problem optimum. This refinement process is sometimes referred to as \emph{infill}. The location at which each newly added data sample is computed is determined through an auxiliary optimization of a given acquisition function (or infill criterion). This routine is repeated until a given convergence criterion is reached. For clarity purposes, a generic mixed-variable BO algorithm is schematically represented in Figure \ref{BO_framework}.

 Relying on GP modeling provides 2 main advantages, the first one being that the surrogate models tend to be characterized by a negligible computational cost when compared to the actual functions of computationally intensive design problems. As a result, these models can be called a large number of times at a negligible cost. Additionally, GP also provide an estimate of the modeling error as a virtually free byproduct of the modeled function prediction, which results particularly useful as it allows to define acquisition functions that provide an efficient compromise between exploration and exploitation of the search space. In order to apply BO algorithms for the solution of VSDSP, a first necessary step is to extend the standard BO methods in order to deal with the presence of generic discrete design variables. In its original formulation, BO is defined for purely continuous problems, \ie problems characterized by functions depending solely on continuous variables \cite{Jones1998}. However, a few extension to the mixed continuous/discrete case exist in the literature \cite{pelamatti2019efficient}, \cite{pelamatti2019surrogate} and \cite{zhang2019bayesian}. The purpose of this paper is to extend the applicability of these mixed-variable BO algorithms in order to enable the optimization of VSDSP. In this section a brief overview of mixed-variable Gaussian Process surrogate modeling and its application within a mixed-variable BO framework is provided. In the first part of this section, the surrogate modeling of mixed-variable functions is discussed. Subsequently, the definition and optimization of an acquisition function in the mixed-variable space are presented.

\subsection{Mixed-variable Gaussian Process surrogate modeling}
The standard GP formulation allows to model purely continuous functions, however, a few variants enabling the modeling of functions which depend simultaneously on continuous and discrete variables exist in the literature. Reviews and comparisons from a theoretical and practical point of view of these existing techniques can be found in  \cite{Swiler2012}, \cite{Zhang2015} and \cite{pelamatti2020overview}. The most commonly used approach when dealing with this kind of functions consists in creating a separate and independent GP model for every category of the considered problem by relying solely on the training data relative to said category \cite{Swiler2012}. However, within the framework of computationally intensive design it is often unfeasible to provide the amount of data for each category of the problem necessary to model the considered function accurately enough for optimization purposes. This issue becomes particularly relevant when dealing with problems characterized by a large number of categories \cite{Swiler2012}. 

 In this section, the GP modeling of functions depending simultaneously on continuous and discrete variables is discussed. Please note that the possible dependency on dimensional variables $\mathbf{w}$ is not considered here, but is discussed in Section \ref{VarDimKernel}. A mixed-variable (or mixed continuous-discrete) GP  can be created by processing a training data set $\mathcal{D}$ of $n$ samples $\{\mathbf{x}^i,\mathbf{z}^i,y^i\}$ with $i\in\{1,...,n\}$, defined as follows:
\begin{equation*}
\mathcal{D} = \left\{ 
\mathbf{X} = \{   \mathbf{x}^1,..., \mathbf{x}^n\} \in F_x, \quad 
\mathbf{Z} =  \{    \mathbf{z}^1,..., \mathbf{z}^n\} \in F_z, \quad 
 \mathbf{y} =  \{ y^1,...,y^n\}  \in F_y \right\}
\end{equation*} 
where $\mathbf{X}$ and $\mathbf{Z}$ are the matrices containing the $n$ continuous and discrete vectors characterizing the training data set, while  $\mathbf{y}$ is the vector containing the associated responses (\emph{i.e.,} modeled function).  $F_x$,  $F_z$ and $F_y$ are the definition domains of the 3 previously mentioned matrices. 
 The core concept of Gaussian Process based surrogate modeling  is to predict the response value $y^*$ of a black-box function $f(\cdot)$ for a generic unmapped input through an inductive procedure by mapping the probability distribution of the possible regression functions. A generic Gaussian Process $Y(\mathbf{x},\mathbf{z})$ is characterized by its covariance function:
\begin{equation}
Cov\left( Y(\mathbf{x},\mathbf{z}), Y(\mathbf{x}',\mathbf{z}')\right)  = E[(Y(\mathbf{x},\mathbf{z}) - \mu(\mathbf{x},\mathbf{z}))(Y(\mathbf{x}',\mathbf{z}') - \mu(\mathbf{x}',\mathbf{z}'))]
\end{equation}
\noindent and its mean function $\mu$. Given that in most real-life engineering design cases, insufficient information regarding the global trend of the modeled functions is known,  it is common practice to consider the regression function $\mu$ as a being constant  with respect to the design space \cite{Simpson2001}:
\begin{equation}
\mu(\mathbf{x},\mathbf{z})  = \mu
\end{equation} 
The predicted value $f^*$ of a given function at an unmapped location $\{\mathbf{x}^*, \mathbf{z}^*\}$ is computed under the form of a Gaussian distribution conditioned on the data set \cite{Rasmussen2006}:
\begin{equation}
f^*|\mathbf{x}^*, \mathbf{z}^*,\mathbf{X},\mathbf{Z},\mathbf{Y} \sim \mathcal{N} \left(   \hat{y}(\mathbf{x}^*,\mathbf{z}^*) ,  \hat{s}^2 (\mathbf{x}^*,\mathbf{z}^*) \right) 
\end{equation}
with a mean value $\hat{y}(\mathbf{x}^*,\mathbf{z}^*)$:
\begin{eqnarray}
\hat{y}(\mathbf{x}^*,\mathbf{z}^*)   & = &  \mathbb{E}[f^*|\mathbf{x}^*, \mathbf{z}^*,\mathbf{X},\mathbf{Z},\mathbf{Y}] \\
 & = & \mu +\bm{\psi}^T(\mathbf{x}^*,\mathbf{z}^*) \mathbf{K}^{-1}(\mathbf{y}-\mathbf{1} \mu) \nonumber
\label{Prediction}
\end{eqnarray}
and associated variance $\hat{s}^2( \mathbf{x}^*, \mathbf{z}^*)$:
\begin{small}
\begin{eqnarray}
\hat{s}^2( \mathbf{x}^*, \mathbf{z}^*) & = & Var(f^*|\mathbf{x}^*, \mathbf{z}^*,\mathbf{X},\mathbf{Z},\mathbf{Y}) \\  \nonumber
& = & k(\{ \mathbf{x}^*, \mathbf{z}^*\},\{ \mathbf{x}^*, \mathbf{z}^*\})  - \bm{\psi}^T( \mathbf{x}^*, \mathbf{z}^*) \mathbf{K
}^{-1} \bm{\psi}( \mathbf{x}^*, \mathbf{z}^*) \nonumber   
\label{PredictionVar}
\end{eqnarray}
\end{small}
\noindent where $ \mathbf{K}$  is the $n\times n$ Gram covariance matrix containing the covariance values between the $n$ samples of the data set:
\begin{equation}
\mathbf{K}_{i,j} = Cov\left( Y(\mathbf{x},\mathbf{z}), Y(\mathbf{x}',\mathbf{z}' )\right) =k( \{ \mathbf{x}, \mathbf{z}\} , \{ \mathbf{x}', \mathbf{z}'\})
\end{equation}
$\mathbf{y}$ is a   $n\times 1$ vector containing the responses corresponding to the $n$ data samples:
\begin{equation}
y_i = f(\mathbf{x}^i,\mathbf{z}^i) \ \mbox{ for } i=1,\dots ,n
\end{equation}
  $\mathbf{1}$ is a $n\times 1$ vector of ones and finally $ \bm{\psi} $ is an $n\times 1$ vector containing the covariance values between each sample of the training data set and the point at which the function is predicted:
\begin{equation}
 \bm{\psi}_i( \mathbf{x}^*, \mathbf{z}^*) = Cov(Y(\mathbf{x}^*,\mathbf{z}^*),Y(\mathbf{x}^i,\mathbf{z}^i)) =   k(\{ \mathbf{x}^*, \mathbf{z}^*\},\{ \mathbf{x}^i, \mathbf{z}^i\}) \qquad \mbox{ for } i=1,\dots,n
 \label{CovVector}
\end{equation}
Within the GP framework, the covariance function is defined as an input space dependent parameterized function, as is discussed in the following section.

\subsection{Gaussian Process kernels}
\label{GPkernels}
The covariance function $k(\cdot)$ is the core component of a Gaussian Process based surrogate model \cite{Rasmussen2006}. Loosely speaking, the purpose of this function is to characterize the similarity between distinct data samples in the design space with respect to the modeled function.

\subsubsection{Kernel operators}
\label{Kerneloperators}
In complex design problems, a single kernel may not be sufficient in order to capture the different influences of the various design variables. The main reason for this limitation, is that the same single set of hyperparameters is used to characterize the covariance between every dimension of the compared samples. For this reason, it is common practice to combine kernels defined over various sub-spaces of the design space, thus resulting in a valid kernel defined over the entire search space which provides a more accurate modeling of the considered function. It can be shown that kernels can be combined while still resulting in a valid covariance function, as long as the chosen operator allows it \cite{steinwart2008support}. In this paper, the three following kernel operators are considered:

\begin{itemize}

\item \textbf{Sum}

Let $k_1(\cdot)$ be a kernel defined on the input space $F_{x_1}$ and $k_2(\cdot)$ be a kernel defined on $F_{x_2}$. It can be shown that $k(\cdot) = k_1(\cdot) + k_2(\cdot) $ is a valid kernel on the input space $F_x = F_{x_1} \times F_{x_2}$. 

\item \textbf{Product}

Let $k_1(\cdot)$ be a kernel defined on the input space $ F_{x_1} $ and $k_2(\cdot)$ be a kernel defined on $F_{x_2}$. It can be shown that $k(\cdot) = k_1(\cdot) \times k_2(\cdot) $ is a valid kernel on the input space $F_x = F_{x_1} \times F_{x_2}$. Furthermore, let $k(\cdot)$ be a kernel defined on the input space $F_x$ and $\alpha \in \mathbb{R}^+$, $\alpha k(\cdot) $ is also a valid kernel on $ F_x $.

\item \textbf{Mapping}

Let $k(\cdot)$ be a kernel on $F_x $, let $ \tilde{F}_x$ be a set and let $A : \tilde{F}_x \rightarrow F_x$ be a mapping function. Then $\tilde{k}(\cdot)$ defined as $\tilde{k}(\mathbf{x},\mathbf{x}') := k(A(\mathbf{x}),A(\mathbf{x}'))$ is a kernel on $ \tilde{F}_x$ 

\end{itemize}

The combination of kernels has 2 main advantages: first, it enables defining distinct kernels over different design variables, which allows to better capture and model the influence the various design variables on the modeled function. This is particularly important when dealing with problems in which different variables are related to considerably different trends of the modeled function, as often happens within the framework of system design. Furthermore, different kernels can also be combined in order to model more accurately different trends characterized by the same design variable or group of variables. 

A popular approach when dealing with multidimensional problems consists in defining a single separate kernel $k_d(\cdot)$ for each dimension $d$ the considered problem depends on, and computing the global kernel as the product between one-dimensional kernels \cite{Santner2003}. In the purely continuous case, this results in the following kernel definition:
\begin{equation}
\label{KernelProd}
k( \mathbf{x}, \mathbf{x}') = \prod_{d = 1}^{n_x}  k_d(x_d,x_d')
\end{equation}
where $n_x$ represents the total dimension of the input space $F_x$. By doing so, each variable of the modeled function can be associated to a specific kernel parameterization as well as a set of specific associated hyperparameters, thus providing a more flexible modeling of the considered function. 

A similar logic can also be relied on in order to deal with mixed-variable problems: rather than defining a single kernel on the mixed continuous/discrete design space, it is possible to define two distinct and independent kernels: $k_x(\cdot)$ with respect to the continuous variables $\mathbf{x}$ and $k_z(\cdot)$ with respect to the discrete variables $\mathbf{z}$. Subsequently, the mixed variable kernel $k(\{ \mathbf{x}, \mathbf{z}\},\{ \mathbf{x'}, \mathbf{z'} \})$ can be defined as:
\begin{equation}
k(\{ \mathbf{x}, \mathbf{z}\},\{ \mathbf{x}', \mathbf{z}' \}) = k_{x}(\mathbf{x},\mathbf{x}')  k_{z} (\mathbf{z},\mathbf{z}')
\end{equation}
Moreover, the kernel defined above can be further decomposed as a product of one-dimensional kernels, similarly to what is shown in Eq. \ref{KernelProd}. The resulting mixed-variable kernel can then be defined as:
\begin{equation}
\label{KernelProd2}
k(\{ \mathbf{x}, \mathbf{z}\},\{ \mathbf{x}', \mathbf{z}' \})  = \prod_{d = 1}^{n_x}  k_{x_d}(x_d,x_d') \prod_{d = 1}^{n_z}  k_{z_d}(z_d,z_d')
\end{equation}
A large number of  parameterizations for the continuous kernel $k_x$ exist in the literature, such as the  linear kernel \cite{Rasmussen2006}, the polynomial kernel \cite{scholkopf2001learning} and the Matern kernels class \cite{Minasny2005}. Due to the fact that this paper focuses on the discrete part of the kernels, the squared exponential parameterization (sometimes referred to as RBF) \cite{Rasmussen2006} is considered for all the continuous kernels:
\begin{equation}
\label{SquaredExpKern}
k_x(x,x') =  \sigma_x^2 \exp\left(-\theta||x-x'||^2\right) 
\end{equation}
where $\theta$ and $\sigma_x^2$ are respectively the lengthscale hyperparameter and the variance associated to the kernel.
\subsubsection{Discrete kernel construction}
Kernels characterizing the covariance between discrete variables are less commonly discussed in the literature when compared to the continuous ones, and only a handful of parameterizations have been proposed. The main difference with what is described in the previous paragraph lies within the fact that the discrete variables that are considered within the framework of system design usually present a finite number of levels as well as the fact that the numerical representation assigned to the considered variables levels is usually arbitrary, and does therefore not yield useful information. In order to avoid this issue, discrete kernels can be constructed by relying on the RKHS formalism. In order for a function $k_z(z,z')$ to represent a valid covariance, there are two main requirements \cite{aronszajn1950theory}. More specifically, the function must be symmetric:
\begin{equation}
k( z,z') = k( z',z) 
\end{equation}
and positive semi-definite over the input space, \emph{i.e.}:
\begin{eqnarray}
&& \sum_{i = 1}^{n} \sum_{j = 1}^{n} a_i a_j k( z,z' ) \geq 0 \\
&& \forall n \geq 1, \ \forall (a_1,\dots , a_n) \in \mathbb{R}^n \ \ \mbox{ and } \ \ \forall z,z' \in  F_z \nonumber
\end{eqnarray}
Alternatively, a positive semi-definite function can also be defined by ensuring that the resulting Gram matrix is positive semi-definite.

Thanks to the characteristics mentioned above, a valid covariance function can, by construction, be defined as a parameterizable Hilbert space kernel \cite{steinwart2008support}. Let $F_z $ be a non-empty set, a kernel function on $F_z$, \emph{i.e.,} $k: F_z \times F_z \rightarrow \mathbb{R}$, can be defined if there exists an $\mathbb{R}$-Hilbert space and a map $\phi : F_z \rightarrow \mathcal{H}$ such that $\forall \mathbf{z},\mathbf{z}' \in F_z$:
\begin{equation}
\label{KernelDef}
k(\mathbf{z},\mathbf{z}') := \langle \phi (\mathbf{z}) , \phi(\mathbf{z}') \rangle_\mathcal{H} 
\end{equation}
where $\langle \cdot, \cdot \rangle_\mathcal{H} $ is the inner product on the Hilbert space $\mathcal{H}$. By definition, an inner product defined in a Hilbert space must be \cite{steinwart2008support}, \cite{scholkopf2001learning}:

\begin{itemize}
\item bi-linear: $\langle \alpha_1 \phi (z) + \alpha_2 \phi (z'), \phi (z'') \rangle _\mathcal{H} = \alpha_1 \langle \phi (z),\phi (z'') \rangle _\mathcal{H} + \alpha_2 \langle \phi (z'), \phi (z'') \rangle _\mathcal{H} $
\item symmetric: $\langle \phi (z),\phi (z') \rangle _\mathcal{H} =  \langle \phi (z'),\phi (z) \rangle _\mathcal{H}$
\item positive: $\langle \phi (z),\phi (z) \rangle _\mathcal{H} \geq 0 $, \quad $\langle \phi (z),\phi (z) \rangle _\mathcal{H} = 0$ \  if and only if \ $ \phi (z) = 0$
\end{itemize}
for any $z,z',z'' \in F_z$. The Hilbert space $\mathcal{H}$ can be seen as a space onto which specific features of the considered design variables are mapped. If necessary, multiple mappings and associated kernels can be defined for the same design variable in case several distinct features of the considered variable must be taken into account. In order to better capture and model the dependence of the GP covariance function with respect to the various variables of the design space, both the mapping function $\phi(\cdot)$ and the Hilbert space inner product $\langle\cdot\rangle_\mathcal{H}$ can depend on a number parameters, known as hyperparameters \cite{Rasmussen2006}. The kernel hyperparameters can be tuned as a function of the available data in order to maximize the modeling accuracy of the surrogate model. 

The discrete kernel construction described above can be used in order to define various valid parameterizations. In this paper, two different parameterizations are considered: the compound symmetry and the latent variable kernels.

\subsubsection{Compound Symmetry}
\label{CS}
The first and most simple discrete kernel to be considered in this paper is the Compound Symmetry (CS), characterized by a single covariance value for any non-identical pair of inputs \cite{Roustant2018}:
\begin{equation}
k(z,z') =  \begin{cases} 
\label{CSkern}
\sigma^2_z \quad \mbox{if} \quad \quad z = z'  \\ 
\theta \cdot \sigma^2_z  \quad \mbox{if} \quad z \neq z'
\end{cases}
\end{equation}
where $\sigma^2_z$ and $0 < \theta < 1$ are respectively the variance and the hyperparameter associated to the CS kernel. By definition, in case the input data samples are identical, the kernel returns the associated variance value $\sigma_z^2$,  as can be seen in Eq. \ref{CSkern}. Alternatively, the variance computed between any pair of non-identical data samples is independent from the inputs, and is equal to a value ranging from $0$ to $\sigma^2_z$.  It should be pointed out that the CS kernel presented above is nearly identical to the one proposed by Hutter and Halstrup in \cite{Halstrup2016} and \cite{Hutter2009}, although the construction differs. The approach proposed in these works consists in defining a distance in the mixed-variable search space by relying on the concept of Gower distance \cite{Gower1971}.
\\[12pt]
The CS kernel presented in the paragraph above provides a very simple method allowing to model the effect of a given discrete design variable by relying on a single hyperparameter. However, the underlying assumption which is made when considering the CS kernel is that the covariance between any pair of non identical levels of a given discrete variable is the same, regardless of the considered levels. This assumption may often result overly simplistic, especially when dealing with discrete variables which present a large number of levels. In this case, the modeling error introduced by CS kernel can become problematic, and alternative kernels should be considered. It is also important to point out that the CS kernel formulation, as presented in Eq. \ref{CSkern}, can only return positive covariance values (by construction). This characteristic further limits the number of suitable applications for this particular covariance function. Finally, it is worth mentioning that Roustant \emph{et al.} have extended the CS kernel in order to model mixed-variable functions characterized by discrete variables with a large number of levels \cite{Roustant2018}. The underlying idea is to group levels with similar characteristics, thus allowing to compute the covariance between said groups rather than between the levels. However, given that in the scope of this work the considered problems present a limited number of levels, this approach is not further discussed.

\subsubsection{Latent variable kernel}
\label{LVSection}
\begin{figure}[h!]
\centering
 \includegraphics[width=0.55\linewidth]{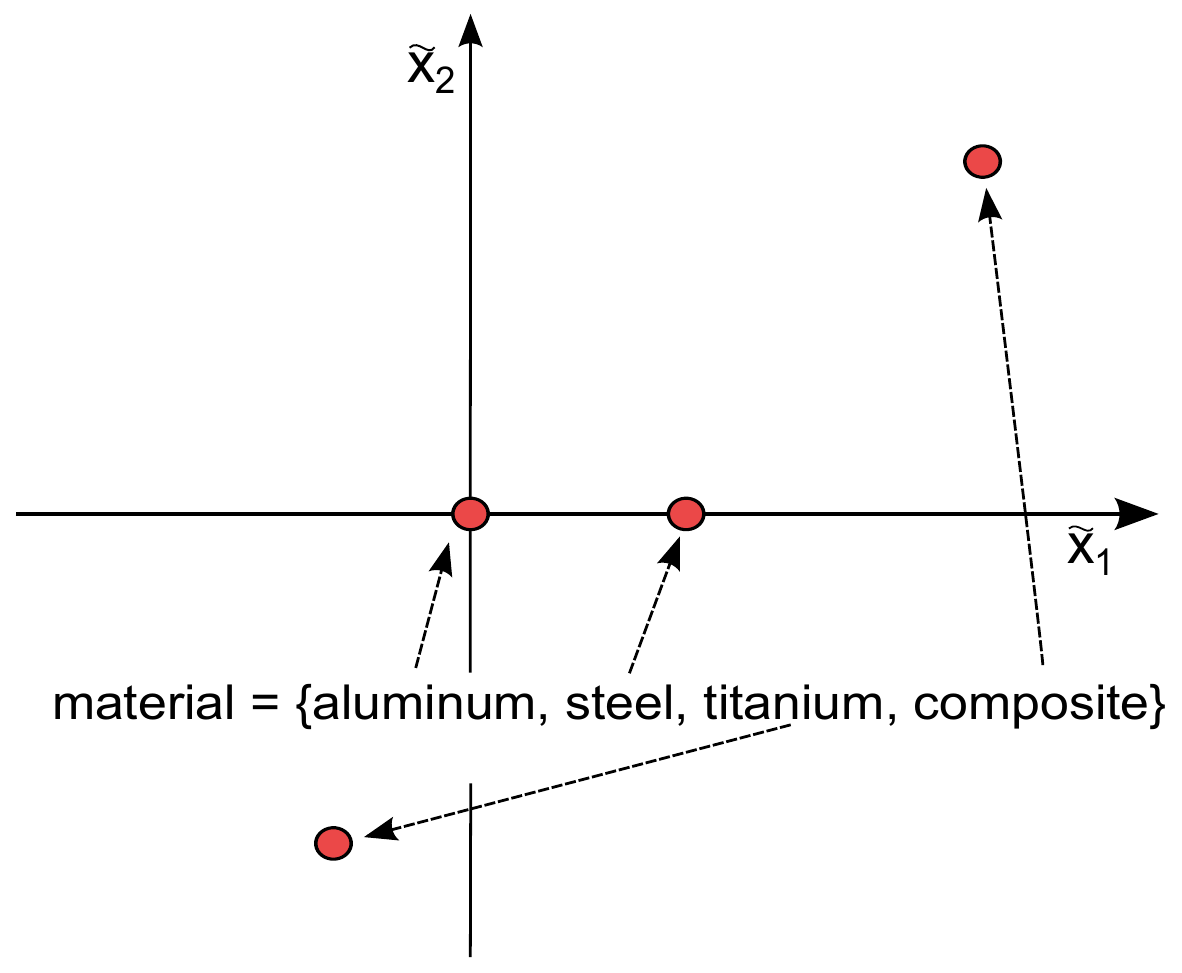}
\caption{Example of latent variable mapping for a generic discrete variable characterizing the 'material choice' characteristic.}
\label{LV}
\end{figure}

The Latent Variable (LV) kernel, first proposed by Zhang \emph{et al.} \cite{zhang2019latent}, is constructed by mapping the discrete variable levels onto a 2-dimensional Euclidean space regardless of the number of discrete levels. This mapping can be defined as follows:
\begin{eqnarray}
\label{LVMapping}
&& \phi(z) : F_z \rightarrow \mathbb{R}^2 \\  \nonumber
&& \phi(z = z_m) = [ \theta_{m,1}, \theta_{m,2} ]^T \quad \mbox{ for } m=1,\dots, l
\end{eqnarray}
where $\theta_{m,1}$ and $\theta_{m,2}$ are the hyperparameters representing coordinates in the $2$-dimensional latent variable space onto which the discrete variable level $m$ is mapped. By consequence, the set of hyperparameters characterizing this kernel is represented by the $l$ pairs of latent variable coordinates associated to a given discrete variable. For clarity purposes, an example of the mapping of a discrete variable characterizing the material choice property of a hypothetical system is provided in Figure \ref{LV}.
\\[12pt]
Let $ \phi : F_z \rightarrow \mathbb{R}^2$ be the mapping defined in Eq. \ref{LVMapping}, by applying the mapping rule discussed in Section \ref{Kerneloperators}, it can be shown that a kernel $\tilde{k}(\cdot)$ valid on $\mathbb{R}^2$ can be used in order to define a kernel $k(\cdot)$ valid on $F_z$ in the following fashion \cite{steinwart2008support}:
\begin{equation}
k(z,z') = \tilde{k}(\phi(z), \phi(z'))
\end{equation}
By consequence, any of the continuous kernels valid on $\mathbb{R}^2$ can be coupled with the latent variable mapping in order to define a valid kernel on the discrete search space $F_z$. Although in the original formulation of the method, the following squared exponential kernel is considered:
\begin{equation}
\label{LVKernel2}
k(z,z') = \sigma_z^2 \exp(  || \phi(z) - \phi(z')||^2_2 )
\end{equation}
alternative continuous kernels could technically be used without loss of generality. It can be noticed that no lengthscale parameter $\theta$ appears in the squared exponential kernel as defined in Eq. \ref{LVKernel2}. The reason behind this is that the distance between the latent variables already directly depends on the hyperparameter values (\emph{i.e.,} the latent variables coordinates), and by consequence a lengthscale parameter would result redundant. In order to remove the translation and rotation indeterminacies on the hyperparameter estimation, one of the latent variables coordinate pair is fixed to the origin of the 2-dimensional latent search space and a second pair on the $\theta_1$ axis. 

The latent variable kernel construction discussed above relies on mapping the discrete variable levels onto a $2$-dimensional Hilbert space. This choice is arbitrary as the mapping can technically be performed onto a higher dimensional space. However, Zhang \emph{et al.} \cite{zhang2019latent} state that the theoretical improvement of modeling performance does not compensate the increase in the number of hyperparameters to be tuned. Due to the fact the covariance values are computed as a function of the distance between latent variables in a continuous space, as shown in Eq. \ref{LVKernel2}, the returned values can not be negative, which partially limits the modeling capabilities of the LV discrete kernel. 

By mapping the levels of the considered discrete variable onto a $2$-dimensional latent space and by characterizing the covariance between said levels as the Euclidean distance between the latent variables, the latent variable kernel provides an intuitive visual representation of how the levels are correlated to each other, as is shown in Figure \ref{LV}. In said figure, the distance between the latent variables (\emph{i.e.,} red dots) associated to the various materials is inversely proportional to the covariance between the relative levels. For instance, with the considered hyperparameter values, the computed covariance between the aluminum and steel choices would result larger than the one between the aluminum and the composite choices.

As previously mentioned, additionally to the two kernels presented above, a few other discrete kernels exist, such as the hypersphere decomposition \cite{Zhou2011} and a kernel based on the used of a coregionalization matrix \cite{alvarez2012kernels}. More detailed analyses on the theoretical and practical differences between the various discrete kernels within the context of mixed-variable GP modeling can be found in \cite{Swiler2012}, \cite{Zhang2015} and \cite{pelamatti2020overview}. However, please note that although for synthesis purposes only the CS and the LV kernels are considered in this paper, the proposed approaches are valid regardless of the selected kernel parameterization.

\subsection{Mixed variable acquisition function}
Within the BO framework, the mean predicted value of the objective and constraint functions of the considered problem, as well as their variance can be used in order to define an acquisition function, which has the purposes of determining the most promising locations of the search space, in terms of both objective function value and feasibility. The location at which the true problem functions is to be evaluated at each iteration is determined by optimizing said function. Depending on the characteristics of the considered problem, different acquisition functions may be considered. The common approach consists in defining this function as being comprised of 2 terms. The first one is tasked with determining the location of the design space at which the best value of the objective function is most likely to be found. The second term, instead, has the purposes of driving the infill search towards regions of the search space where the largest number of constraints is most likely to be satisfied. 

\subsubsection{Objective function oriented infill criterion}
In this work, the first term of the acquisition function is defined as a mixed-variable adaptation of the Expected Improvement (EI) first defined in \cite{mockus1994application}. As the name suggests, the EI represents the expected value of the improvement  $I = \max \left(y_{min} - Y(\mathbf{x^*}, \mathbf{z}^*),0 \right)$ in terms of objective function value with respect to the data set. The original purely continuous formulation is the following:
\small
\begin{eqnarray}
\hspace{-0.7cm}
\mathbb{E}[I(\mathbf{x^*}, \mathbf{z}^*)] & = & \mathbb{E} \left[\max \left(y_{min} - Y(\mathbf{x^*}, \mathbf{z}^*),0 \right) \right] \\
 & = & (y_{min} - \hat{y}(\mathbf{x^*}, \mathbf{z}^*)\Phi \left(\frac{y_{min}-\hat{y}(\mathbf{x^*}, \mathbf{z}^*)}{\hat{s}(\mathbf{x^*}, \mathbf{z}^*)}\right) + \hat{s}(\mathbf{x^*}, \mathbf{z}^*)\phi \left(\frac{y_{min}-\hat{y}(\mathbf{x^*}, \mathbf{z}^*)}{\hat{s}(\mathbf{x^*}, \mathbf{z}^*)}\right) \nonumber
\label{EI}
\end{eqnarray}
\normalsize
\noindent where $y_{min}$ is the minimum feasible value present within the data set at the given BO iteration, $\hat{y}(\mathbf{x^*}, \mathbf{z}^*)$ and $\hat{s}(\mathbf{x^*}, \mathbf{z}^*)$ are the mean and standard deviation of the objective function prediction $Y(\mathbf{x^*})$, $\Phi(.)$ is the cumulative distribution function of a normal distribution and finally $\phi(.)$ is the probability density function of a normal distribution. 

The EI criterion is originally derived by considering purely continuous GP \cite{mockus1994application}. However, it can be easily shown that said derivation still holds in the mixed-variable case as long as the prediction provided by the GPs can be represented under the form of a normally distributed variable, \emph{i.e.,}  $\mathcal{N} \left(   \hat{y}(\mathbf{x}^*,\mathbf{z}^*) ,  \hat{s}^2 (\mathbf{x}^*,\mathbf{z}^*) \right)$. This property can be ensured by defining a valid mixed-variable kernel characterized by hyperparameters within their limit bounds. 

\subsubsection{Feasibility oriented infill criteria}
In order to take into account the presence of constraints and thus providing feasible final designs, the objective function oriented infill criterion discussed in the previous paragraph must be combined with an auxiliary criterion. In this paper, the Expected Violation (EV) criterion \cite{audet2000surrogate} is considered as it allows to fairly compare candidate solutions characterized by different constraints (in terms of both number and type). However, other solutions exist, as is for instance discussed in \cite{Durantin2016}, \cite{priem2019use} and  \cite{Sasena2002}. In a similar way to the EI, the EV represents the expected value of the violation of a given constraint, \emph{i.e.,} the difference between the predicted value and the maximum acceptable value, which is usually set to $0$:
\begin{equation}
V_i = \mbox{max} \left(G_i(\mathbf{x^*},\mathbf{z^*}) - 0,0  \right)
\end{equation}
The EV for a given constraint $g_i(\cdot)$ is defined as follows:
\small
\begin{eqnarray}
\hspace{-0.38cm}
\mathbb{E}[V_i(\mathbf{x^*},\mathbf{z^*})] & = & \mathbb{E} [ \mbox{max} \left(G_i(\mathbf{x^*},\mathbf{z^*}) - 0,0 \right)] \\
 & = & (\hat{g}_i(\mathbf{x^*},\mathbf{z^*})-0)\Phi \left(\frac{\hat{g}_i(\mathbf{x^*},\mathbf{z^*})-0}{\hat{s}_{g_i}(\mathbf{x^*},\mathbf{z^*})}\right) + \hat{s}_{g_i}(\mathbf{x^*},\mathbf{z^*})\phi \left(\frac{\hat{g}_i(\mathbf{x^*},\mathbf{z^*})}{\hat{s}_{g_i}(\mathbf{x^*},\mathbf{z^*})}\right) \nonumber
\label{EV}
\end{eqnarray}
\normalsize
In the same way as the previously discussed infill criterion, the EV derivation still holds in the mixed-variable case as long as the GP prediction of each constraint follows a normal distribution (\ie  $G_i(\mathbf{x^*},\mathbf{z^*}) \sim \mathcal{N} \left(   \hat{g}_i(\mathbf{x}^*,\mathbf{z}^*) ,  \hat{s}_{g_i}^2 (\mathbf{x}^*,\mathbf{z}^*) \right)$).

\subsubsection{Infill criterion optimization}

If EI and EV are considered, the  acquisition function is defined as a constrained problem. In other words, the location $\{\mathbf{x}^n,\mathbf{z}^n\}$ at which the actual objective and constraint functions of the considered problem are evaluated at a given BO iteration is determined through the optimization of the following constrained auxiliary problem:
\begin{eqnarray}
\{\mathbf{x}^n,\mathbf{z}^n\} & = & \mbox{argmax} \left(EI ( \mathbf{x},\mathbf{z})\right) \\
 \mbox{s.t. } && EV_i ( \mathbf{x},\mathbf{z}) \leq t_i \quad \mbox{ for } \quad i=1,\dots ,n_g \nonumber \\ \nonumber
\text{w.r.t.} &\;& \mathbf{x} \in F_x \subseteq  \mathbb{R}^{n_{x}}  \\ \nonumber
&& \mathbf{z}  \in F_z  \\ \nonumber
\end{eqnarray}
where $t_i$ is the maximum accepted violation for the constraint $g_i$. When dealing with mixed-variable problems, the acquisition function must be optimized within the mixed continuous/discrete search space and is therefore subject to the limitations related to the presence of discrete and unoredered variables. The results presented in this work are obtained by optimizing the infill criterion with  the help of a mixed continuous/discrete Genetic Algorithm (GA) similar to the one presented by Stelmack \cite{Stelmack1998} and coded with the help of the python based toolbox DEAP \cite{DEAP_JMLR2012}. The specific parameters of said mixed-variable GA, such as mutation and cross-over probabilities as well as number of generations and population size, vary between test-cases and are adapted to the size and characteristics of the considered problem.

As is schematically represented in Figure \ref{BO_framework}, once the values of $\{\mathbf{x}^*, \mathbf{z}^*\}$ that yield the optimal acquisition function value have been determined, the exact objective and constraint functions of the optimization problem can be computed at said location and the obtained data sample can be added to the GP training data set. Subsequently, the surrogate models must be trained anew in order to take into account the additional information provided by the newly computed data sample. This routine is repeated until a user-defined stopping criterion is reached. In order to better represent the conditions of real-life computationally intensive design problems, in this paper a pre-defined number of data samples to be infilled is selected and the BO process is repeated until said computational budget is exhausted. This also allows to provide a fair comparison between the considered mixed-variable BO variants in terms of convergence speed.

\section{Budget allocation strategy}
\label{BudgetAllocationStrat}
As mentioned in the introduction, due to the limitations of the existing algorithms allowing to solve VSDSP, the most commonly used approach consists in decomposing the global problem into $N_p$ fixed-size mixed continuous/discrete sub-problems which can be independently optimized by relying on more standard algorithms, such as the mixed-variable BO presented in the previous section. In practice, each sub-problem $q$ can be obtained by fixing the dimensional variables $\mathbf{w}$ to a given set of possible values  $\mathbf{w}_q$ and is formulated as follows:
\begin{eqnarray}
\label{SPBDef}
\min &\;& f(\mathbf{x},\mathbf{z},\mathbf{w}_q) \; \ \qquad  \quad f: \mathbb{R}^{n_{x}(\mathbf{w}_q)} \times \prod_{d = 1}^{n_{z}(\mathbf{w}_q)} F_{z_d} \times F_w  \rightarrow F_f \subseteq \mathbb{R} \\ \nonumber
\text{w.r.t.} &\;& \mathbf{x} \in F_x(\mathbf{w}_q) \subseteq \mathbb{R}^{n_{x}(\mathbf{w}_q)}  \\ \nonumber
&& \mathbf{z}  \in \prod_{d = 1}^{n_{z}(\mathbf{w}_q)} F_{z_d} \\ \nonumber
\text{s.t.} &\;& \mathbf{g}  (\mathbf{x},\mathbf{z},\mathbf{w}_q) \leq  0 \\ \nonumber
&& g_i : F_{x_i}(\mathbf{w}_q) \times \prod_{d = 1}^{n_{z_i}(\mathbf{w}_q)} F_{z_{d_i}} \times F_w \rightarrow F_{g_i} \subseteq \mathbb{R} \qquad \mbox{for} \qquad i = 1,...,n_g(\mathbf{w}_q)
\end{eqnarray}
The optimization problem presented in Eq. \ref{SPBDef}, obtained by fixing the values of $\mathbf{w}$ to $\mathbf{w}_q$ (\ie one of $\mathbf{w}$ categories), depends on a fixed number of continuous and discrete variables and is subject to a fixed number of constraints, and can therefore be solved by relying on the mixed-variable BO described in Section \ref{MVBO}. Alternatively, the sub-problems can be enumerated by regrouping all the dimensional variables as a single variable $w$ defined in the combinatorial space of $\mathbf{w}$. In this case, each sub-problem is simply defined by one of the possible levels of the scalar variable $w$. It is important to remember that although the various sub-problems defined as in Eq. \ref{SPBDef} can be characterized by a different set of design variables as well as different constraints, their objective function always models the same quantity with the same unit of measure, as otherwise the comparison between sub-problems would lose meaning.

In order to determine the global optimum of a VSDSP as defined in Eq. \ref{PBDef}, a separate and independent sub-problem optimization must be performed for every possible category of  $\mathbf{w}$. As a consequence, this approach quickly results computationally inefficient, if not unfeasible, when dealing with optimization problems which present a large dimensional variable combinatorial search space,  due to the fact that all the non-optimal sub-problems must be optimized until convergence regardless of the provisional results obtained along the optimization process. This issue becomes particularly problematic when dealing with computationally intensive problems, which are particularly common within the context of complex system design. In order to partially avoid this issue when dealing with VSDSP, the first approach for the optimization of variable-size design space problems which is proposed in this paper is based on coupling the separate and independent optimization of each sub-problem defined by Eq. \ref{SPBDef} with a computational budget allocation strategy. For the sake of conciseness, this approach is referred to as Strategy for the Optimization of Mixed Variable-Size design space Problems (SOMVSP). The underlying idea is to rely on the information provided by the surrogate models of each sub-problem objective and constraint functions along the optimization process in order to determine which sub-problems are more likely to contain the global optimum. At every SOMVSP iteration, the information provided by the various sub-problems surrogate models is exploited in order to allocate a different computational budget to each sub-problem (in terms of infilled data samples) proportionally to how promising it is. Moreover, in order to further optimize the usage of the overall computational budget, the sub-problems least likely to contain the global optimum can also discarded between SOMSVP iterations. In synthesis, the proposed strategy relies on two main axes: the allocation of a different optimization computational budget to each sub-problem as a function of the predicted optimal solution (in terms of both objective function value and feasibility), and the possibility of discarding a given sub-problem in case its predicted relative performance with respect to the other sub-problems crosses a given threshold. These two mechanisms are detailed in the following paragraphs.

\subsection{Discarding of non-optimal sub-problems}
The proposed computational effort allocation logic relies on both the dimension of the considered sub-problem and its predicted relative performance with respect to the objective function value as well as the solution feasibility. At each SOMVSP iteration, the relative performance of each sub-problem $q$ with respect to the others must be determined. In order to do so, three different predicted optima are computed for each sub-problem by considering different confidence interval scenarios: the \emph{Best Case} (BC), the \emph{Worst Case} (WC) and the \emph{Nominal Case} (NC). In practice, these scenarios correspond to the predicted feasible optimum value of the sub-problem $q$  by considering an optimistic, pessimistic and null value of the predicted variance, respectively, and they are defined as follows:

\begin{eqnarray}
&& BC_q = \min \hat{y}_q (\mathbf{x}, \mathbf{z},\mathbf{w}_q) - a \cdot \hat{s}_q (\mathbf{x}, \mathbf{z},\mathbf{w}_q) \\
&& \mbox{ s.t. } EV\left(g_c(\mathbf{x}, \mathbf{z},\mathbf{w}_q)\right) < t_{c} \ for \ c = 1,...,n_{g_i}(\mathbf{w}_q) \nonumber
\end{eqnarray}

\begin{eqnarray}
&&WC_q = \min \hat{y}_i (\mathbf{x}, \mathbf{z},\mathbf{w}_q) + a \cdot \hat{s}_q (\mathbf{x}, \mathbf{z},\mathbf{w}_q) \\ 
&& \mbox{ s.t. } EV\left(g_c(\mathbf{x}, \mathbf{z},\mathbf{w}_q)\right) < t_{c} \ for \ c = 1,...,n_{g_i}(\mathbf{w}_q) \nonumber
\end{eqnarray}

\begin{eqnarray}
&& NC_q = \min \hat{y}_q (\mathbf{x}, \mathbf{z},\mathbf{w}_q) \\ 
&& \mbox{ s.t. } EV\left(g_c(\mathbf{x}, \mathbf{z},\mathbf{w}_q)\right) < t_{c} \ for \ c = 1,...,n_{g_i}(\mathbf{w}_q) \nonumber
\end{eqnarray}
where $a \in \mathbb{R}^+$ is a tunable parameter representing how conservative are the definitions of the BC and WC scenarios. In practice, $a$ represents the confidence the user gives to the GP models modeling accuracy of the various sub-problems functions. In general, larger values of $a$ tend to result in more conservative and robust results, however, they also tend to reduce the convergence speed of the proposed SOMVSP. For the sake of simplicity, only values of $a=2$ and $a=3$ are considered in this paper, which approximately correspond to confidence intervals of 95\% and 99\%. 

\begin{figure}[h!]
\centering
 \includegraphics[width=.7\linewidth]{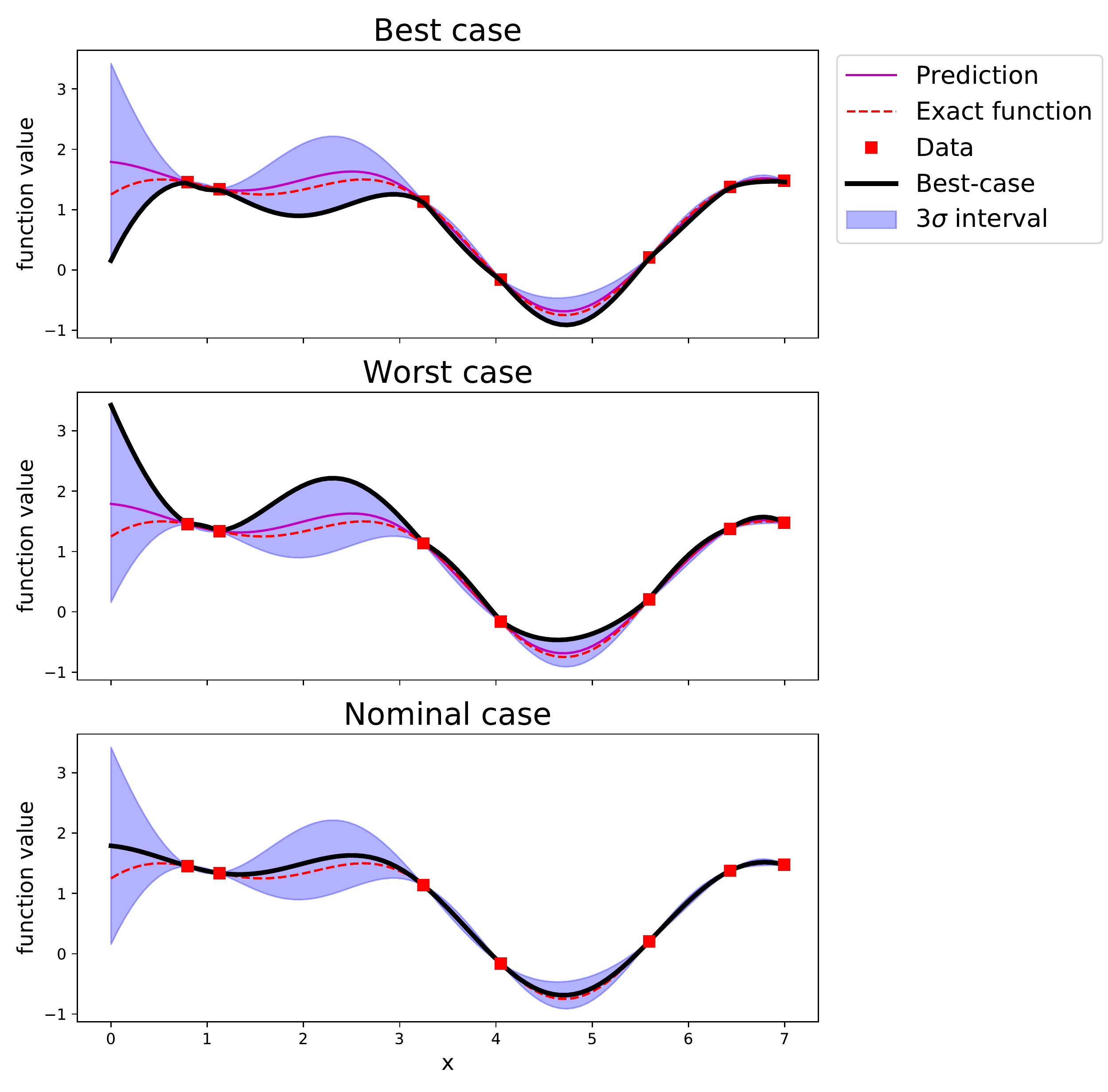}
\caption{Example of Best-case, Worst-case and Nominal-case functions.}
\label{CasesEx}
\end{figure}

\begin{figure}[h]
\centering
 \includegraphics[width=.8\linewidth]{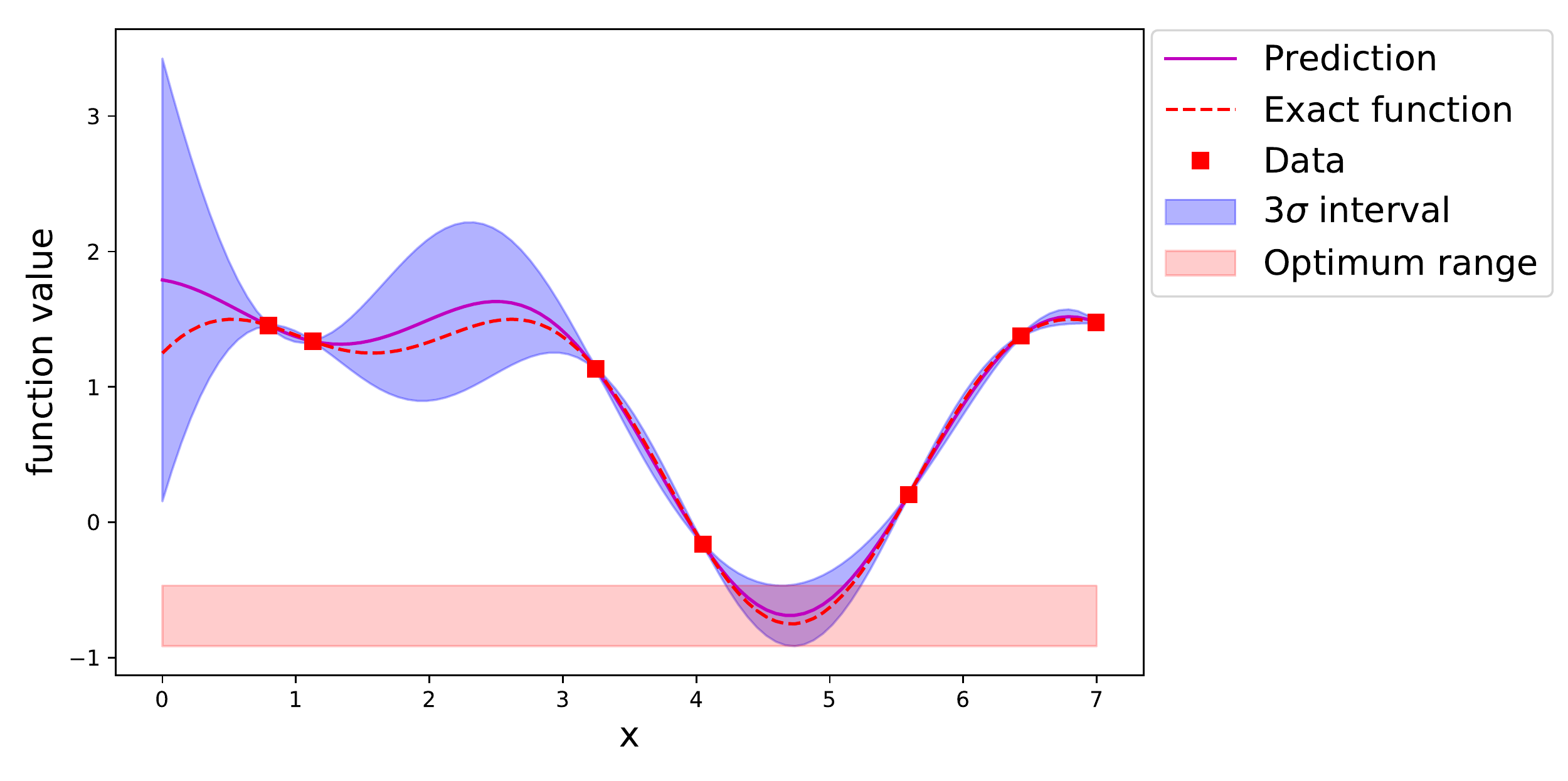}
\caption{Example of a function optimum range, defined as the interval between the BC and the WC}
\label{CasesEx_OptRange}
\end{figure}

For illustrative purposes, the following unconstrained one-dimensional continuous function is considered:
\begin{equation}
f(x) = \cos(x)^2 + \sin(x) + 0.25
\end{equation}
A GP model of said function is created by relying on 7 data samples, and the Best-case, Worst-case and Nominal-case functions associated to the given model are provided in Figure \ref{CasesEx}.
In practice,  by considering the minimum value of the Best-case and the Worst-case scenario, it is then possible to define an optimum range, in which the optimum of the modeled function is most likely to be, with a given confidence (proportional to the value of $a$). This is illustrated in Figure \ref{CasesEx_OptRange}. Please note that the considered example is unconstrained. In the presence of constraints, the presented scenario would only be defined in the regions of the design space for which the EV of each constraint is lower than the given limit.

\begin{figure}[h!]
\centering
 \includegraphics[width=1\linewidth]{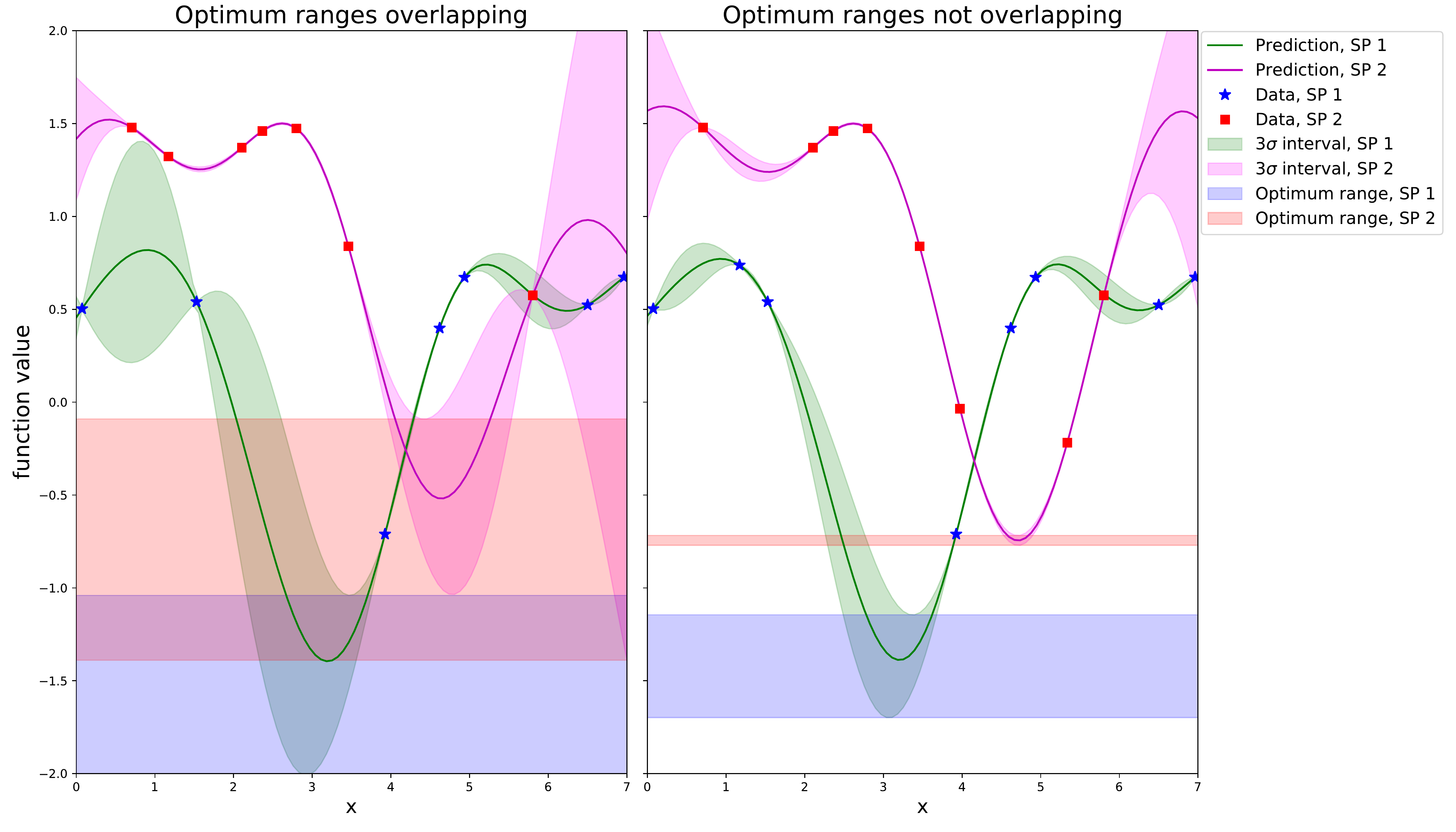}
\caption{Example of comparison between optimum ranges between two Sub-Problems (SP) over different iterations. In the left figure, the optimum ranges of the two SP overlap and it is therefore not possible to accurately predict whether one of the two contains the global optimum. In the right figure, instead, the two optimum ranges do not overlap and it is therefore possible to discard the SP2 as it does most likely not contain the global optimum.}
\label{SPComp}
\end{figure}

At the beginning of each SOMVSP iteration, the BC and the WC values  are used in order to compare the predicted performance of each sub-problem and use this information in order to discard the less promising ones. More specifically, if the WC predicted feasible optimum of a given sub-problem yields a lower value than the BC predicted feasible optimum of a different sub-problem, it is expected that the latter is not worth further exploring as it will most likely not contain the global optimum. In this case, the non-optimal sub-problem is discarded and is not further explored for the rest of the optimization process. This is equivalent to determining whether the optimum range associated to two given sub-problems overlap or not. In practice, a sub-problem $q$ is discarded if:
\begin{equation}
\label{SPComparison}
BC_q \geq WC_p \quad \mbox{ for } \quad p = 1,...,N_p \quad \mbox{ and } \quad p \neq q
\end{equation}
This condition must be evaluated for ever remaining sub-problem at the given iteration. For illustrative purposes, the comparison between 2 sub-problems is illustrated in Figure \ref{SPComp}. As can be noticed, on the left side of the figure the variance associated to the two models is considerable, and the optimum ranges overlap. Therefore, in this case none of the two sub-problems would be discarded. Instead, the right size of the figure shows the comparison between the GPs obtained with a few more data samples per sub-problem. In this case the models are more accurate and the associated variance is smaller. As a consequence, it can be seen that the predicted optimum of sub-problem 1 yields a lower value than the one of sub-problem 2 with a given confidence level. In this case, the sub-problem 2 will be discarded as it is most likely that it does not contain the global optimum.
Once the non-optimal sub-problems have been  discarded (if present), the computational budget for the given SOMVSP iteration is allocated among the remaining sub-problems.

\subsection{Computational budget allocation}
At every SOMVSP iteration, the discarding of non-optimal sub-problems is followed by the allocation of a different computational budget to each remaining sub-problem. More specifically, at every iteration a budget of $B_q$ data samples to be infilled is allocated to each remaining sub-problem $q$. In order to fairly compare sub-problems characterized by a potentially considerably different number of design variables, $B_q$ is computed by taking into account both the predicted performance of a given sub-problem as well as its dimension, and it is defined as follows: 
\begin{equation}
B_q = d_q\left(\frac{1+\Delta_q}{2}\right)
\end{equation}
where $d_q$ is the total dimension of the sub-problem $q$ (\emph{i.e.,} sum of the continuous and discrete dimensions), while $\Delta_q$ is a term representing the relative performance of the considered sub-problem with respect to the remaining others. It is computed as: 
\begin{equation}
\Delta_q = \frac{NC_{max}-NC_q}{NC_{max}-NC_{min}}
\end{equation}
where $NC_{max}$ and $NC_{min}$ are respectively the largest and lowest NC values among the remaining sub-problems. As a result, a budget equal to half its total dimension is assigned to the least promising sub-problem, whereas a computational budget equal to its total dimension is assigned to the most promising one. 

\subsection{Bayesian optimization of remaining sub-problems}
Following the computational budget allocation, each remaining sub-problem is independently optimized with the help of a mixed-variable BO similar to the one presented in Section \ref{MVBO}, and by infilling a number of data samples proportional to the allocated budget. The newly infilled data sample for each sub-problem is defined as follows:
\begin{eqnarray}
\label{SubProbInfill}
\{\mathbf{x}^n,\mathbf{z}^n\} & = & \mbox{argmax} \left(EI ( \mathbf{x},\mathbf{z})\right) \\
 \mbox{s.t. } && EV_i ( \mathbf{x},\mathbf{z}) \leq t_i \quad \mbox{ for } \quad i=1,\dots ,n_g \nonumber \\
 \mbox{w.r.t. } && \mathbf{x} \in F_x, \ \mathbf{z} \in F_z \nonumber
\end{eqnarray}
Please note that in this case, $\mathbf{x}$ and $\mathbf{z}$ refer to continuous and discrete variables the considered sub-problem depends on. By consequence, the acquisition function optimization is performed in a different search space for each considered sub-problem. Similarly, each infill process is only subject to the EV related to the constraints the considered sub-problem is subject to.

\subsection{Algorithm overview}
\begin{figure}[h!]
\centering
 \includegraphics[width=.9\linewidth]{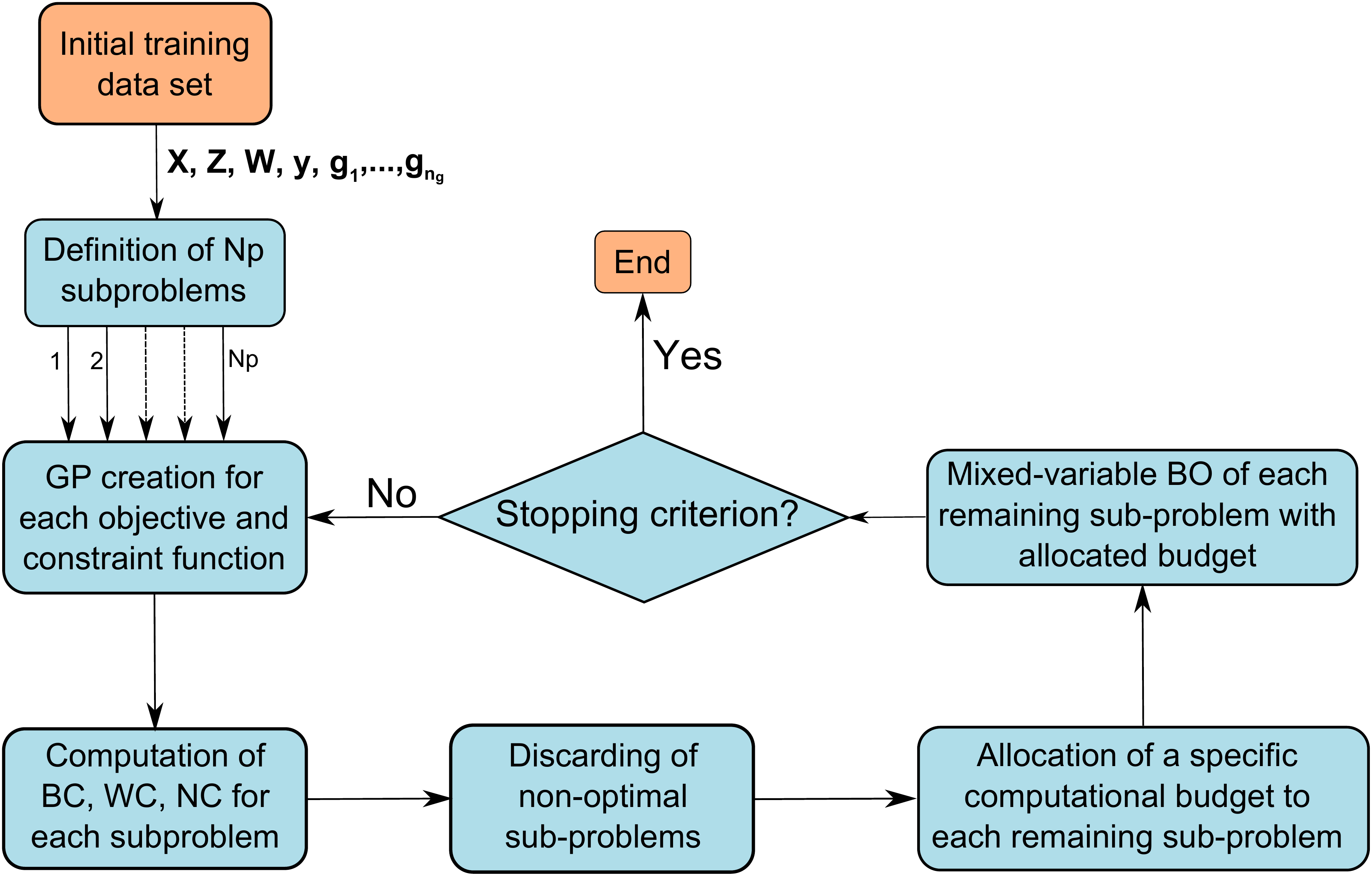}
\caption{Budget allocation strategy for the optimization of variable-size design space problems.}
\label{SOMVSP}
\end{figure}

In the previous paragraphs, the main steps comprising an iteration of the proposed SOMVSP are detailed. By repeating this process, the number of remaining sub-problems is expected to gradually decrease, thus allowing to focus the computational budget on the most promising ones. In general, the proposed optimization strategy continues until a stopping criterion is reached. In this work, a predefined total computational budget (\ie number of function evaluations) is allocated to the optimization process, which performs SOMVSP iterations until said computational budget is depleted. This choice allows to analyze the performance of the proposed algorithm in a context of computationally intensive design problems. For illustrative purposes, a visual representation of the SOMVSP algorithm is provided in Figure \ref{SOMVSP}.

The SOMVSP described in the previous paragraphs relies on the idea of using the information provided by the various sub-problems GP models in order to determine which ones are more likely to actually yield the global optimum, and allocating the computational budget accordingly. However, it might be important to note that the criteria which are used in order to determine how promising a given sub-problem is and how to allocate the computational budget among sub-problems are defined empirically. Therefore, these criteria might require to be adapted or redefined in order to better suit the characteristics of particular considered VSDSP, while still relying on the same working principle. The optimization results which can be obtained with the proposed SOMVSP on two different test-cases are shown in Section \ref{VSDSPRes}.

\section{Variable-size design space kernel}
\label{VarDimKernel}
The Budget allocation strategy for the optimization of variable-size design space problems presented in Section \ref{BudgetAllocationStrat} provides an efficient alternative than the separate and independent optimization of every fixed-sized sub-problem when dealing with VSDSP, as is shown in Section \ref{VSDSPRes}. However, for particularly complex and computationally costly problems characterized by a large dimensional variables combinatorial search space, the proposed method might still not be viable as it would still require a large amount of data in order to identify the most promising sub-problems. For this reason, an alternative approach for the optimization of VSDSP is proposed in this Section. The underlying idea is to adapt the mixed-variable BO discussed in Section \ref{MVBO} for the solution of variable-size design space problems by defining a single Gaussian Process for each variable-size design space function  characterizing a given problem (\ie objective and constraints). In order to do this, it is necessary to define a  kernel allowing to characterize the covariance between data samples defined in the variable-size design space:
\begin{equation}
k\left((\mathbf{x},\mathbf{z},\mathbf{w}),(\mathbf{x}',\mathbf{z}',\mathbf{w}')\right)
\end{equation}
where $\mathbf{x}$, $\mathbf{x}'$, $\mathbf{z}$ and $\mathbf{z}'$ may contain different sets of continuous and discrete variables, respectively. A few kernels allowing to compute the covariance between samples containing partially different sets of purely continuous variables within the context of hierarchical search spaces exist in the literature, such as the arc kernel \cite{hutter2013kernel}, the imputation kernel \cite{zaefferer2018first} and the wedge kernel \cite{horn2019surrogates}. These kernels are all based on the redefinition of the concept of distance within the hierarchical space, the construction of which varies depending of whether both, only one, or neither considered variables are active. However, all of the previously mentioned kernels are only defined within a purely continuous design space, and do not take into account the presence of discrete variables and the dependence of the continuous/discrete design space on the presence of dimensional variables. As a consequence, their applicability to actual design problems is limited. For these reasons, two alternative definitions for a variable-size design space kernel are proposed in this section, and the adaptations required to optimize an infill criterion in a variable-size design space are described. The first proposed approach relies on grouping the GP training data according to the sub-problem the samples belong to, in order to allow computing the covariance between different sub-problems. The second proposed approach, instead, relies on grouping the training data according to their dimensional variables values, thus providing a more complex but more accurate modeling of variable-size design space functions.
\subsection{Sub-problem-wise decomposition kernel}
As discussed in Section \ref{BudgetAllocationStrat}, a generic VSDSP can be decomposed into $N_p$ sub-problems, each one characterized by a fixed-sized mixed variable search space. Furthermore, each sub-problem can be seen as a level of a single (scalar) dimensional variable $w$ defined in the combinatorial space of $\mathbf{w}$, which characterizes the sub-problem the considered sample belongs to. A first possible way of defining a kernel in a variable-size design space consists in separately computing the within sub-problem covariance (\ie covariance between samples belonging to the same sub-problem) and the between sub-problems covariance (\ie covariance between data samples belonging to different sub-problems).  By definition, the within sub-problem covariance can be represented  through a mixed-variable kernel, as discussed in Section \ref{MVBO}. The between sub-problems covariance, instead, can be represented as a discrete kernel defined on the search space of the (combinatorial) dimensional variable $w$.  
The global variable-size design space kernel can then be computed as the sum between the within sub-problem and between sub-problems covariances. This approach is referred to as Sub-Problem-Wise (SPW) decomposition in the remainder of this paper.

Let $q$ be one of the $N_p$ levels of the dimensional variable $w$ and let $\mathbf{x}_q$ and $\mathbf{z}_q$ be the continuous and discrete variables on which the objective and constraint functions associated to the sub-problem $q$ depend. The first proposed variable-size search space kernel is defined as follows:
\begin{equation}
\label{VarDimKernelSBW}
k\left((\mathbf{x},\mathbf{z},w),(\mathbf{x}',\mathbf{z}',w')\right) = \sum_{q = 1}^{N_p} k_{x_q} (\mathbf{x}_q,\mathbf{x}_q') \cdot k_{z_q}(\mathbf{z}_q,\mathbf{z}_q')  \cdot \delta_q(w,w') + k_{w} (w,w')
\end{equation}
where $\delta_q(w,w')$ is a simil-Kronecker function which yields $1$ in case both $w$ and $w'$ are equal to $q$, and $0$ otherwise:
\begin{equation}
\label{DeltaKernel}
\delta_q(w,w')  =  \begin{cases} 
1 \quad \mbox{if} \quad w = w' = q \\ 
0  \quad \mbox{else} 
\end{cases}
\end{equation}
This function can be seen as a kernel constructed in the following fashion:
\begin{equation}
\delta_q(w,w') = k_w(w,w') =  \langle \phi_{\delta_q}(w), \phi_{\delta_q}(w') \rangle =  \phi_{\delta_q}(w)\cdot \phi_{\delta_q}(w') 
\end{equation}
with $\phi_{\delta_q}$ defined as a mapping function returning $1$ for input values equal to $q$ and $0$ otherwise:
\begin{equation}
\phi_{\delta_q}(w)  =  \begin{cases} 
1 \quad \mbox{if} \quad w = q \\ 
0  \quad \mbox{else} 
\end{cases}
\end{equation}
The kernel defined in Eq. \ref{VarDimKernelSBW} is defined through sums and products of one-dimensional kernels. By consequence, the validity of the global kernel can be ensured as long as each one-dimensional kernel is properly defined, which is already discussed in Section \ref{MVBO}. The dimensional variable kernel $k_w$ characterizing the covariance between the sub-problems, instead, only depends on the dimensional variable which is shared by all data samples. It can be constructed with the same logic as for discrete variables and, for clarity purposes, the same parameterizations as for the discrete variables are used for $k_w$ in this work. However, in order to allow for a more refined modeling of the function trends, a heteroscedastic adaptation might be necessary.  In practice, Eq. \ref{VarDimKernelSBW} is defined so that in case the two considered individuals belong to the same sub-problem, one of the $N_p$ terms of the sum contributes to the covariance value while the remaining $N_p-1$ yield $0$. As a result, the total covariance value is computed as the sum between the within sub-problems covariance and the between sub-problem variance. If instead the two individuals do not belong to the same sub-problem, the covariance value is only computed as the between sub-problem covariance, as all the $N_p$ terms of the sum are null. 
\\[12pt]
Similarly to what is proposed by Roustant \emph{et al.} \cite{Roustant2018} for the modeling of discrete variables with large number of levels, the global variable-size design space kernel can also be represented under the form of a symmetric block matrix,
\begin{equation}
A = \left[
\begin{matrix}
W_1 + B_{1,1}  & B_{1,2} & \dots  & B_{1,N_p}\\
B_{2,1}  & W_2 + B_{2,2} & \ddots & \vdots \\
\vdots  & \ddots &  \ddots & B_{N_p-1,N_p} \\
B_{N_p,1}  & \dots & B_{N_p,N_p-1} & W_{N_p} + B_{N_p,N_p}
\end{matrix} \right]
\label{Mat}
\end{equation}
where $W_q$ is the within sub-problem covariance matrix associated to the sub-problem $q$ and computed as: 
\begin{equation}
W_q =k_{x_q} (\mathbf{x}_q,\mathbf{x}_q') \cdot k_{z_q}(\mathbf{z}_q,\mathbf{z}_q')
\end{equation}
while $B_{q,p}$ represents the covariance between the sub-problems $q$ and $p$ computed as:
\begin{equation}
B_{q,p} = k_{w} (w = q,w' = p) 
\end{equation}
The kernel described above allows to compute the covariance between data samples which belong to partially different search spaces by grouping the samples according to the sub-problem they belong to. However, this approach does not allow to exploit the information which may be provided by the fact that part of the design variables can be shared between data samples belonging to different sub-problems. In order to better exploit this additional information, an alternative variable-size design space kernel based on a decomposition by dimensional variable rather than by sub-problem is proposed in the following paragraphs.

\subsection{Dimensional variable-wise decomposition}
Without loss of generality, the generic VSDSP defined in Eq. \ref{PBDef} can formulated in such a way that each dimensional variable is associated to a specific number of continuous and discrete variables which may be active (or not) depending on its value. In the same way, it is also possible to formulate the problem in such a way that each continuous or discrete variable which does not always influence the problem function only depends from a single dimensional variable. In practice, this translates into having a search space and a feasibility domain that directly depend on the dimensional variables levels, rather than on the dimensional variables category (\ie dimensional variable level combinations). The global VSDSP can therefore be decomposed into groups of design variables which depend from a given dimensional variable, plus the continuous and discrete variables which are shared among all the sub-problems. Within the framework of complex system design, this decomposition of the problem can be seen as the variable-size modeling of the different components characterizing a given system.

If this formulation of VSDSP is considered, an alternative variable-size design space kernel can be defined by considering a separate and independent kernel for each dimensional variable (\ie sub-system) and the continuous and discrete variables which depend on it. In fact, it can be noticed that each dimensional variable related to a number of possibly active continuous and discrete design variables depending on its value is equivalent to a VSDSP characterized by a single dimensional variable, and can therefore be modeled through the SPW decomposition kernel described in the previous paragraph. The second variable-size design space kernel proposed in this paper can then be computed as a product of $n_w$ SPW kernel (\ie one for each dimensional variable).  This approach is referred to as Dimensional Variable-Wise (DVW) decomposition in the remainder of this work. Let $\mathbf{x}_{d_{l}}$ and $\mathbf{z}_{d_{l}}$ be the 'active' continuous and discrete variables associated to level $l$ of the $d_{th}$ dimensional variable (\emph{i.e.,} $w_d = l$). Let $\mathbf{x}_{s}$ and $\mathbf{z}_{s}$ be the continuous and discrete variables which are 'shared' between all the samples and do not depend on a dimensional variable. The DVW decomposition kernel can be computed as follows:

\begin{small}
\begin{eqnarray}
\label{VarDimKernelSSW}
k\left((\mathbf{x},\mathbf{z},\mathbf{w}),(\mathbf{x}',\mathbf{z}',\mathbf{w}')\right) = && \prod_{d = 1}^{n_w}  \left( \sum_{l = 1}^{l_{w_d}} k_{x_{d_{l}}} (\mathbf{x}_{d_{l}},\mathbf{x}_{d_{l}}') \cdot k_{z_{d_{l}}}(\mathbf{z}_{d_{l}},\mathbf{z}_{d_{l}}')  \cdot \delta_l(w_d,w_d') + k_{w_d} (w_d,w_d') \right) \nonumber \\ 
&&  k_{x_{s}}(\mathbf{x}_{s},\mathbf{x}_{s}') \cdot k_{z_{s}}(\mathbf{z}_{s},\mathbf{z}_{s}') 
\end{eqnarray}
\end{small}
\noindent where $l_{w_d}$ is the total number of levels characterizing the $d_{th}$ dimensional variable $w_d$. In this case, each term of the product of Eq. \ref{VarDimKernelSSW} can be schematically represented by an $n_{w_d}\times n_{w_d}$ matrix $A_d$ defined in the same way as the matrix presented in Eq. \ref{Mat}.
\begin{equation}
A_d = \left[
\begin{matrix}
W_{d_1} + B_{d_{1,1}}  & B_{d_{1,2}} & \dots  & B_{d_{1,n_{w_d}}}\\
B_{d_{2,1}}  & W_{d_2} + B_{d_{2,2}} & \ddots & \vdots \\
\vdots  & \ddots &  \ddots & B_{d_{n_{w_d}-1,n_{w_d}}} \\
B_{d_{n_{w_d},1}}  & \dots & B_{d_{n_{w_d},n_{w_d}-1}} & W_{d_{n_{w_d}}} + B_{d_{n_{w_d},n_{w_d}}}
\end{matrix} \right]
\end{equation}
The DVW decomposition kernel offers a more accurate modeling of VSDSP with respect to the SPW alternative as it allows to exploit the information which can be provided by design variables shared between data samples belonging to different sub-problems (\ie variables shared between different levels of a given dimensional variables). However, this kernel definition also results more complex and is usually characterized by a larger number of hyperparameters. As a consequence, the training of the GP model might result more difficult in case an insufficient amount of training data is provided.

\subsection{Infill criterion optimization}
In order to perform the BO of VSDSP, it is necessary to define and optimize an acquisition function within the variable-size design space allowing to determine the most promising locations of said search space. In this paper, the data sample $\{\mathbf{x}^n,\mathbf{z}^n,\mathbf{w}^n\}$ to be infilled at each iteration of the proposed VSDSP BO is computed by relying on the following infill criterion:
\begin{eqnarray}
\label{VSIC}
\mathbf{x}^n,\mathbf{z}^n,\mathbf{w}^n  & = & \mbox{argmax} (EI(\mathbf{x},\mathbf{z},\mathbf{w})) \\ \nonumber
& \mbox{s.t.} &  EV\left(g_c(\mathbf{x}, \mathbf{z},\mathbf{w}_q)\right) \leq t_c  \qquad \mbox{for} \qquad c = 1,...,n_g(\mathbf{w}) \\ \nonumber
& \text{w.r.t.} & \mathbf{x} \in F_x(\mathbf{w}) \subseteq \mathbb{R}^{n_{x}(\mathbf{w})}  \\ \nonumber
&& \mathbf{z}  \in \prod_{d = 1}^{n_{z}(\mathbf{w})} F_{z_d} \\ \nonumber
&& \mathbf{w}  \in F_w  \\ \nonumber
\end{eqnarray}
Similarly to the SOMVSP proposed in Section \ref{BudgetAllocationStrat}, the EV criterion is used in order to take into account the presence of constraints. It can be noticed that the optimization problem defined in Eq. \ref{VSIC} is also defined in the variable-size design space, and it is therefore necessary to solve an auxiliary VSDSP optimization. However, the objective and constraint functions of this optimization auxiliary problem present a negligible computational cost when compared to the actual problem. By consequence, the infill criterion optimization problem can be solved by relying on heuristic VSDSP optimization algorithms, such as the GA variant proposed in \cite{Nyew2015}. However, due to the fact that the implementation of said heuristic VSDSP algorithm requires problem specific inputs and parameterization, the results presented in this paper are obtained by optimizing a separate acquisition function for each sub-problem, and by selecting the data sample which yields the best acquisition function value among all the sub-problems:
\begin{align}
\label{FSIC}
\{ \mathbf{x}^n,\mathbf{z}^n,\mathbf{w}^n  \}  = & \mbox{argmax} \left\{\begin{array}{ll}
        \mbox{argmax} &  \left(EI(\mathbf{x},\mathbf{z},\mathbf{w}_q)\right) \\ \\
        \mbox{s.t.} &  EV\left(g_c(\mathbf{x}, \mathbf{z},\mathbf{w}_q)\right) < t_{c} \qquad \mbox{for} \qquad c = 1,...,n_g(\mathbf{w}_q) \\ \\
        \text{w.r.t.} & \mathbf{x} \in F_x(\mathbf{w}_q) \subseteq \mathbb{R}^{n_{x}(\mathbf{w}_q)} \\
       & \mathbf{z}  \in \prod_{d = 1}^{n_{z}(\mathbf{w}_q)} F_{z_d}
        \end{array}\right\} \\
        & \mbox{for} \ \ q = 1,\dots,N_p \nonumber
\end{align}
The $N_p$ optimization problems defined by Eq. \ref{FSIC} are standard mixed-variable optimization problems. By consequence, similarly to the SOMVSP discussed in Section \ref{BudgetAllocationStrat}, the results presented in this paper are obtained by optimizing the infill criterion with the help of a mixed continuous/discrete GA, similar to the one discussed in Section \ref{MVBO}.

\section{Applications and Results}
\label{VSDSPRes}
In order to assess and compare the performance of the two proposed BO-based approaches (\ie SOMVSP and variable-size design space kernel based BO) for the optimization of VSDSP, three test-cases with different characteristics are considered.  More specifically, two analytical test-cases and an aerospace engineering design problem are optimized. The proposed methods are applied to each benchmark with different discrete kernel parameterizations as well as different parameters (\ie SPW and DVW approaches for the variable-size design space kernel and different values of $a$ for the SOMVSP). Please note that for clarity purposes, the SOMVSP results are referred to with the acronym BA (Budget Allocation).
The objective is to assess the relative performance between the two proposed approaches as well as the impact of the considered algorithm parameters on the resulting convergence speed. The 2 discrete kernel parameterizations which are considered for these benchmarks are the Compound Symmetry (CS) and the Latent Variables (LV). Because of the implementation limitations and poor performance of the few existing heuristic algorithms allowing to solve VSDSP, the reference method which is considered for the presented benchmarks is defined as a separate and independent mixed-variable BO of each sub-problem characterizing the considered VSDSP (referred to as IO) by relying on the algorithm presented in Section \ref{MVBO}. The total computational budget allocated for each given test-case is distributed among the sub-problems proportionally to their total dimension. The initial data set which is provided to compare optimization algorithms is sampled on the global design space, \ie on every continuous and discrete design variable the considered problem depends on. The continuous variables are sampled through a single continuous LHS \cite{McKay1979}, while the discrete variables are drawn from a uniform discrete distribution. Finally, these data samples are randomly associated to each sub-problem (or to a dimensional variable category). The number of samples allocated to each sub-problem is proportional to its total dimension (\ie sum of continuous and discrete dimensions). In order to quantify and compensate the influence of the initial DoE random nature, each optimization problem is solved 10 times with different initial training data sets.

\subsection{Benchmark analysis}
In the following paragraphs, the 2 proposed alternative solutions for the optimization of VSDSP are tested on a number of test-cases with different kernel parameterizations and different parameters. These benchmarks present different characteristics in terms of number of sub-problems, combinatorial complexity of the dimensional variable design space, as well as in terms of sub-problem design space dimensions and sub-problem specific constraints. The main properties of the considered test-cases as well as the simulation details are provided below:

\textbf{Variable-size design space Goldstein function}
\begin{itemize}
\item 5 continuous variables, 4 discrete variables, 2  dimensional variables

\item 8 sub-problems 

\item 648 equivalent continuous problems

\item 1 constraint

\item Initial data set size: 104 samples ( \ie 2 samples per dimension of each sub-problem)

\item Number of infilled samples: 104

\item 10 repetitions

\item Compared methods: 

- Independent mixed-variable BO of each sub-problem  (IO) with both CS and LV kernels

- SOMVSP (BA) with CS and LV kernels and values of $a$ of $2$ and $3$

- Variable-size design space kernel BO with CS and LV kernels with both SPW and DVW approaches

\item Acquisition function: $EI$ under $EV$ constraints

\end{itemize}

\textbf{Variable-size design space Rosenbrock function}
\begin{itemize}
\item 8 continuous variables, 3 discrete variables, 2  dimensional variables

\item 4 sub-problems 

\item 32 equivalent continuous problems

\item 2 constraints

\item Initial data set size: 30 samples (\ie 1 sample per dimension of each sub-problem)

\item Number of infilled samples: 65

\item 10 repetitions

\item Compared methods: 

- Independent mixed-variable BO of each sub-problem  (IO) with both CS and LV kernels

- SOMVSP (BA) with CS and LV kernels and values of $a$ of $2$ and $3$

- Variable-size design space kernel BO with CS and LV kernels with both SPW and DVW approaches

\item Acquisition function: $EI$ under $EV$ constraints

\end{itemize}

\textbf{Multi-stage launch vehicle design}
\begin{itemize}
\item 18 continuous variables, 14 discrete variables, 3  dimensional variables

\item 6 sub-problems 

\item 29136 equivalent continuous problems

\item 19 constraints

\item Initial data set size: 122 samples (\ie 1.5 sample per dimension of each sub-problem)

\item Number of infilled samples: 58

\item 10 repetitions

\item Compared methods: 

- Independent mixed-variable BO of each sub-problem  (IO) with both CS and LV kernels

- SOMVSP (BA) with CS and LV kernels and values of $a$ of $3$

- Variable-size design space kernel BO with CS and LV kernels with both SPW and DVW approaches

\item Acquisition function: $EI$ under $EV$ constraints

\item Large number of continuous and discrete design variables. Large number of constraints.

\end{itemize}

\noindent \textbf{Implementation}

\noindent The results presented in the following paragraphs are obtained with the following implementation. The optimization routine overhead is written in Python 3.6. The GP models are created with the help of GPflow \cite{GPflow2017}, a Python based toolbox for GP-based modeling relying on the Tensorflow framework \cite{tensorflow2015-whitepaper} (version 1.13). The surrogate model training is performed with the help of a Bounded Limited memory Broyden - Fletcher - Goldfarb – Shanno (L-BFGS-B) algorithm \cite{byrd1995limited}, whereas the acquisition functions are optimized the help of a constraint domination based mixed continuous/discrete Genetic Algorithm (GA) \cite{Stelmack1998} implemented by relying on the Python based toolbox DEAP \cite{DEAP_JMLR2012}.

\subsection{Variable-size design space Goldstein function}
\begin{figure}[h!]
\centering
\includegraphics[width=.5\linewidth]{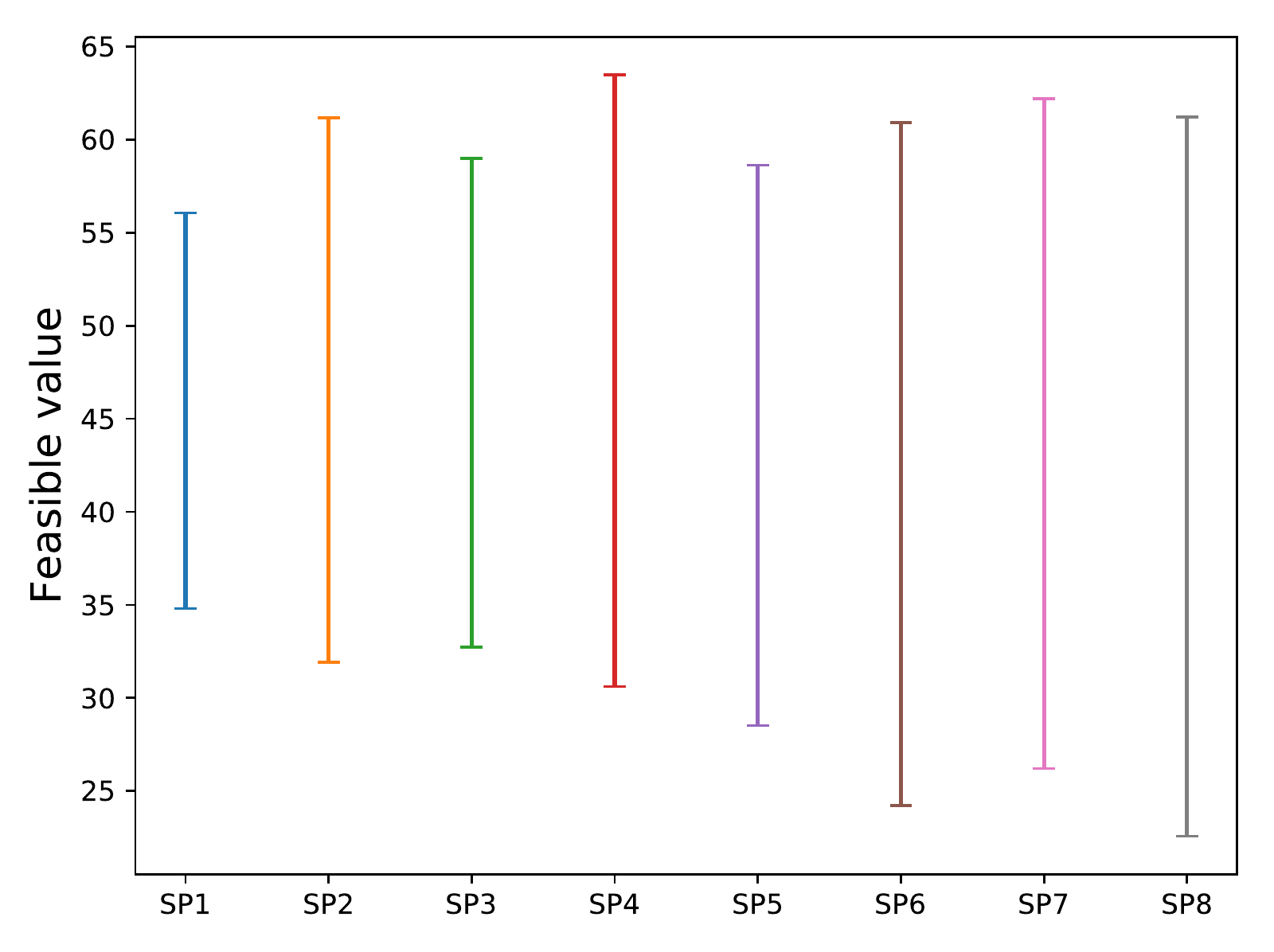}
\caption{Value range of the feasible objective function for each sub-problem of the variable-size design space Goldstein function.}
\label{FeasVal_Goldstein}
\end{figure}

The first analytical test-case that is considered in this paper is  a modified constrained variable-size design space version of the Goldstein \cite{picheny2013benchmark}. The global VSDSP is characterized by 5 continuous variables, 4 discrete variables and 2 dimensional variables. Depending on the dimensional variable values, 8 different sub-problems can be identified, with total dimensions of 6 or 7, ranging from 2 continuous variables and 4 discrete variables to 5 continuous variables and 2 discrete variables. All of the sub-problems are subject to a variable-size design space constraint. The resulting optimization problem can be defined as follows:

\begin{align}
 \min & \qquad  f(\textbf{x},\textbf{z},\textbf{w}) \\ 
\text{w.r.t.} &  \qquad \textbf{x} = \{x_1,\dots,x_5\} \ \mbox{ with } \  x_i \in [0,100] \ \mbox{ for }   i = 1,5  \nonumber \\
 &  \qquad \textbf{z} = \{z_1,\dots,z_4\} \ \mbox{ with } \  z_i \in \{ 0,1,2\}  \ \mbox{ for }   i = 1,4 \nonumber \\
 & \qquad \textbf{w}  = \{w_1,w_2\} \ \mbox{ with } \  w_1 \in \{ 0,1,2,3\}  \ \mbox{ and }  w_2 \in \{ 0,1\}  \nonumber \\
\text{s.t.:} &  \qquad g(\textbf{x},\textbf{z},\textbf{w}) \leq 0  \nonumber 
\end{align}

\noindent where:

\begin{equation}
  f(\textbf{x},\textbf{z},\textbf{w}) = \begin{cases}
               f_1(x_1,x_2,z_1,z_2,z_3,z_4) \hfill   \qquad  \mbox{ if }  w_1 = 0  \mbox{ and } w_2 = 0\\
               f_2(x_1,x_2,x_3,z_2,z_3,z_4) \hfill   \qquad  \mbox{ if }  w_1 = 1  \mbox{ and } w_2 = 0\\
               f_3(x_1,x_2,x_4,z_1,z_3,z_4) \hfill   \qquad  \mbox{ if }  w_1 = 2  \mbox{ and } w_2 = 0\\
               f_4(x_1,x_2,x_3,x_4,z_3,z_4) \hfill   \qquad  \mbox{ if }  w_1 = 3  \mbox{ and } w_2 = 0\\
                f_5(x_1,x_2,x_5,z_1,z_2,z_3,z_4) \hfill   \qquad  \mbox{ if }  w_1 = 0  \mbox{ and } w_2 = 1\\
               f_6(x_1,x_2,x_3,x_5,z_2,z_3,z_4) \hfill   \qquad  \mbox{ if }  w_1 = 1  \mbox{ and } w_2 = 1\\
               f_7(x_1,x_2,x_4,x_5,z_1,z_3,z_4) \hfill   \qquad  \mbox{ if }  w_1 = 2  \mbox{ and } w_2 = 1\\
               f_8(x_1,x_2,x_3,x_5,x_4,z_3,z_4) \hfill   \qquad  \mbox{ if }  w_1 = 3  \mbox{ and } w_2 = 1
            \end{cases}
\end{equation}

\noindent and:

\begin{equation}
  g(\textbf{x},\textbf{z},\textbf{w}) = \begin{cases}
               g_1(x_1,x_2,z_1,z_2) \hfill   \qquad  \mbox{ if }  w_1 = 0 \\
               g_2(x_1,x_2,z_2) \hfill   \qquad  \mbox{ if }  w_1 = 1 \\
               g_3(x_1,x_2,z_1) \hfill   \qquad  \mbox{ if }  w_1 = 2 \\
               g_4(x_1,x_2,z_3,z_4) \hfill   \qquad  \mbox{ if }  w_1 = 3 
            \end{cases}
    \end{equation}
For clarity purposes, the analytical definitions of $f_1(\cdot),\dots,f_8(\cdot)$ and $g_1(\cdot),\dots,g_4(\cdot)$ are not provided in this section but can be found in Appendix A. A synthesis of the characteristics of the various sub-problems comprising the considered variable-size design space function problem is presented in Table \ref{VarDimGoldSynth}. Furthermore, the value ranges of the feasible objective function for each sub-problem is shown in Figure \ref{FeasVal_Goldstein}. It can be seen that the feasible values of the different sub-problems overlap over a large part of their objective function range, thus making the identification of the optimal sub-problem challenging.
\begin{table}[h]
\centering
\begin{tabular}{|l|c|c|c|c|c|c|c|c|}
\hline
Sub-problem &  SP 1 & SP 2 & SP 3 & SP 4 & SP 5 & SP 6 & SP 7 & SP 8 \\ \hline
N$^\circ$ continuous variables & 2 & 3 & 3 & 4 & 3 & 4 & 4  & 5  \\ \hline
N$^\circ$ discrete variables & 4 & 3 & 3 & 2 & 4 & 3 & 3 & 2  \\ \hline
N$^\circ$ constraints & 1 & 1 & 1 & 1 & 1 & 1 & 1 & 1 \\  \hline \hline 
\multicolumn{9}{|c|}{Global VSDSP} \\ \hline 
 N$^\circ$ continuous variables & \multicolumn{8}{|c|}{5}  \\ \hline
N$^\circ$ discrete variables &\multicolumn{8}{|c|}{4} \\ \hline
N$^\circ$ dimensional variables &\multicolumn{8}{|c|}{2} \\ \hline
N$^\circ$ discrete categories & \multicolumn{8}{|c|}{648} \\ \hline 
N$^\circ$ constraint & \multicolumn{8}{|c|}{1} \\ \hline 
\end{tabular}
\caption{Defining characteristics of the sub-problems comprising the variable-size design space Goldstein function  optimization problem.}
\label{VarDimGoldSynth} 
\end{table}

\begin{figure}[h!]
\centering
\includegraphics[width=1.0\linewidth]{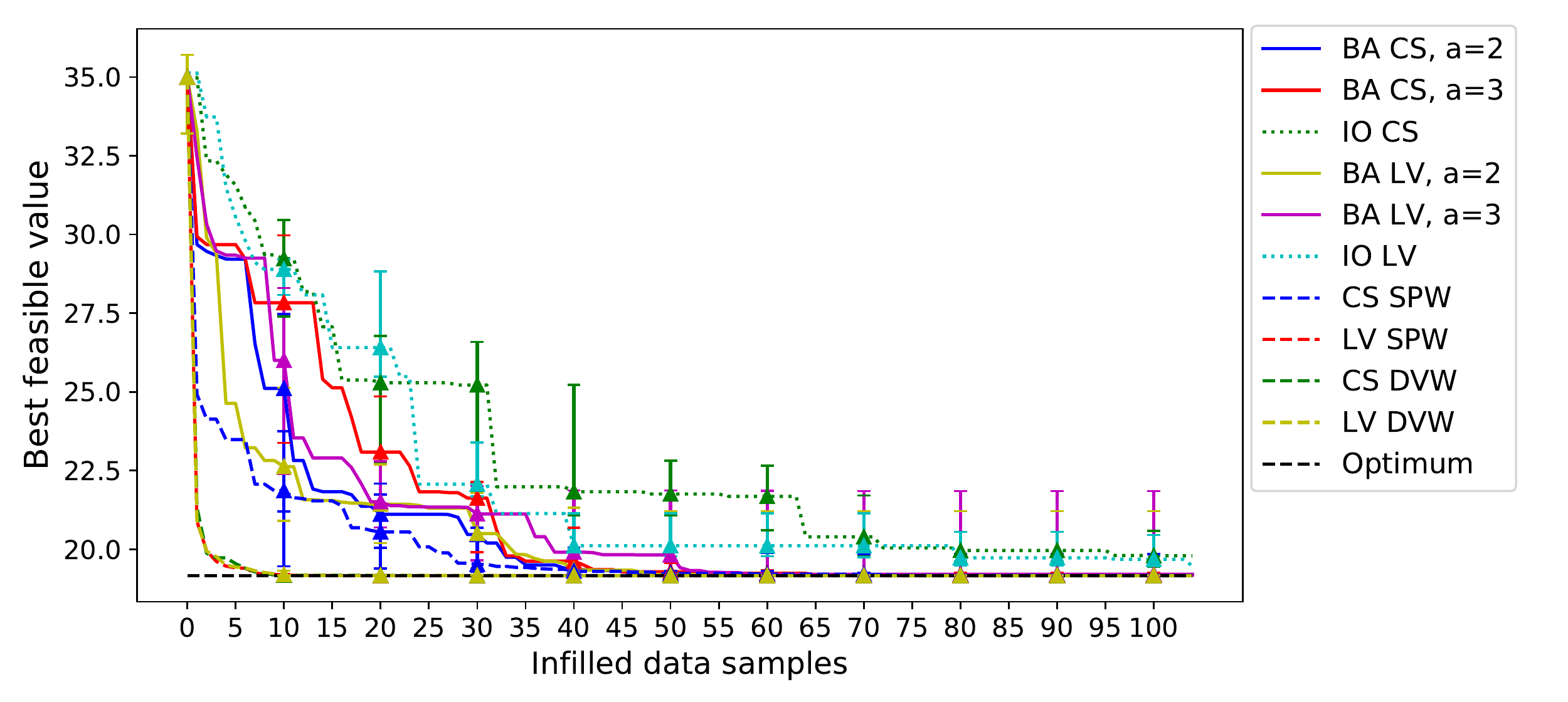}
\caption{Comparison of the convergence rate of various VSDSP optimization algorithms during the BO of the mixed-variable Goldstein function over 10 repetitions. The continuous lines represent the SOMVSP alternatives, the dashed lines represent the variable-size design space kernel BO alternatives and the dotted lines represent the independent optimization of each sub-problem.}
\label{Res_Opt_VarDimGoldstein1}
\end{figure}

\begin{figure}[h!]
\centering
\includegraphics[width=.7\linewidth]{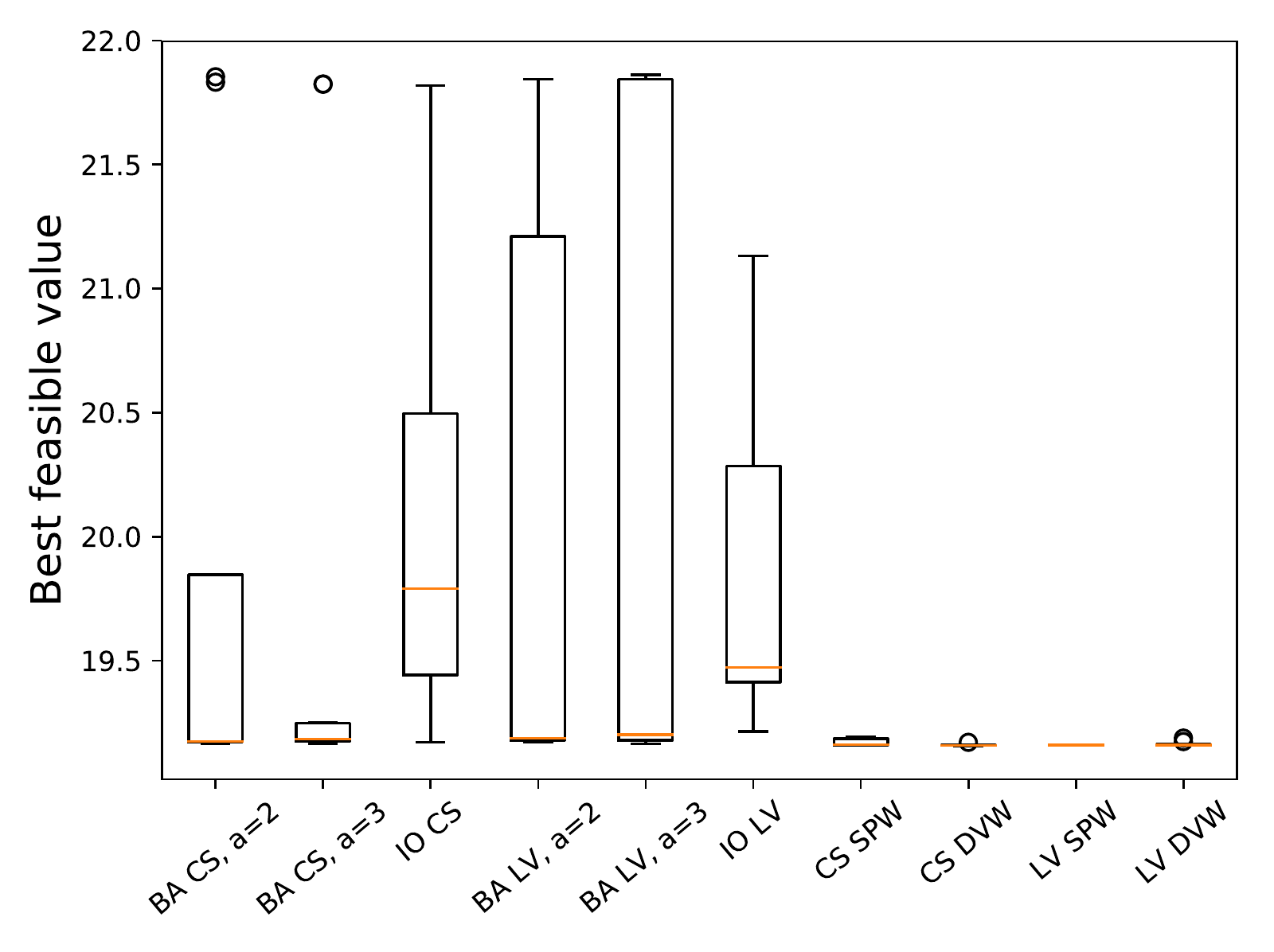}
\caption{Comparison of the convergence value of different VSDSP optimization algorithms on the variable-size design space Goldstein function over 10 repetitions.}
\label{Res_Opt_VarDimGoldstein2}
\end{figure}

\begin{figure}[h!]
\centering
\includegraphics[width=.65\linewidth]{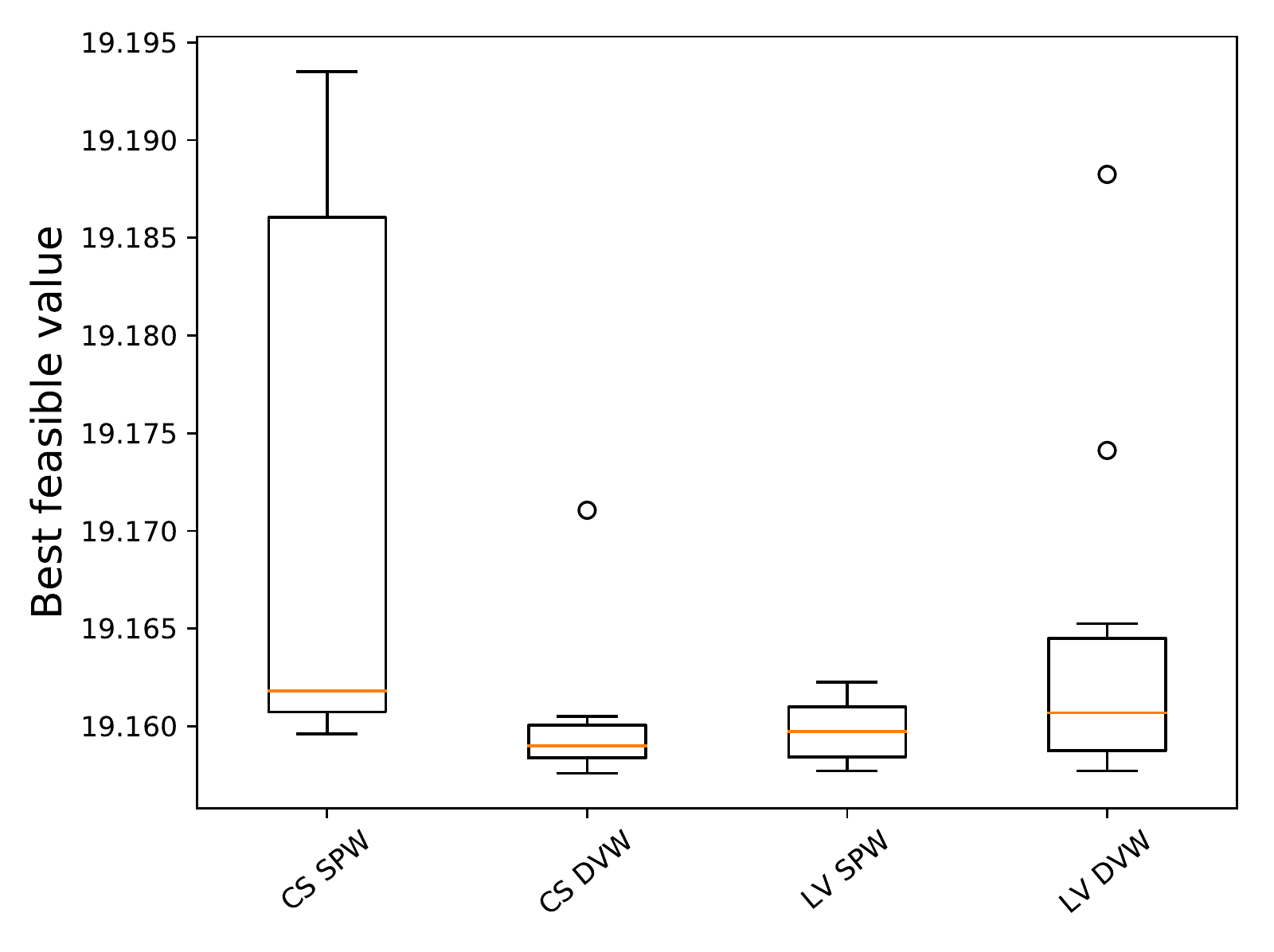}
\caption{Comparison of the convergence value of different VSDSP optimization algorithms on the variable-size design space Goldstein function over 10 repetitions. Focus on the variable-size design space kernel based BO algorithms.}
\label{Res_Opt_VarDimGoldstein3}
\end{figure}

\begin{figure}[h]
\centering
\includegraphics[width=1.0\linewidth]{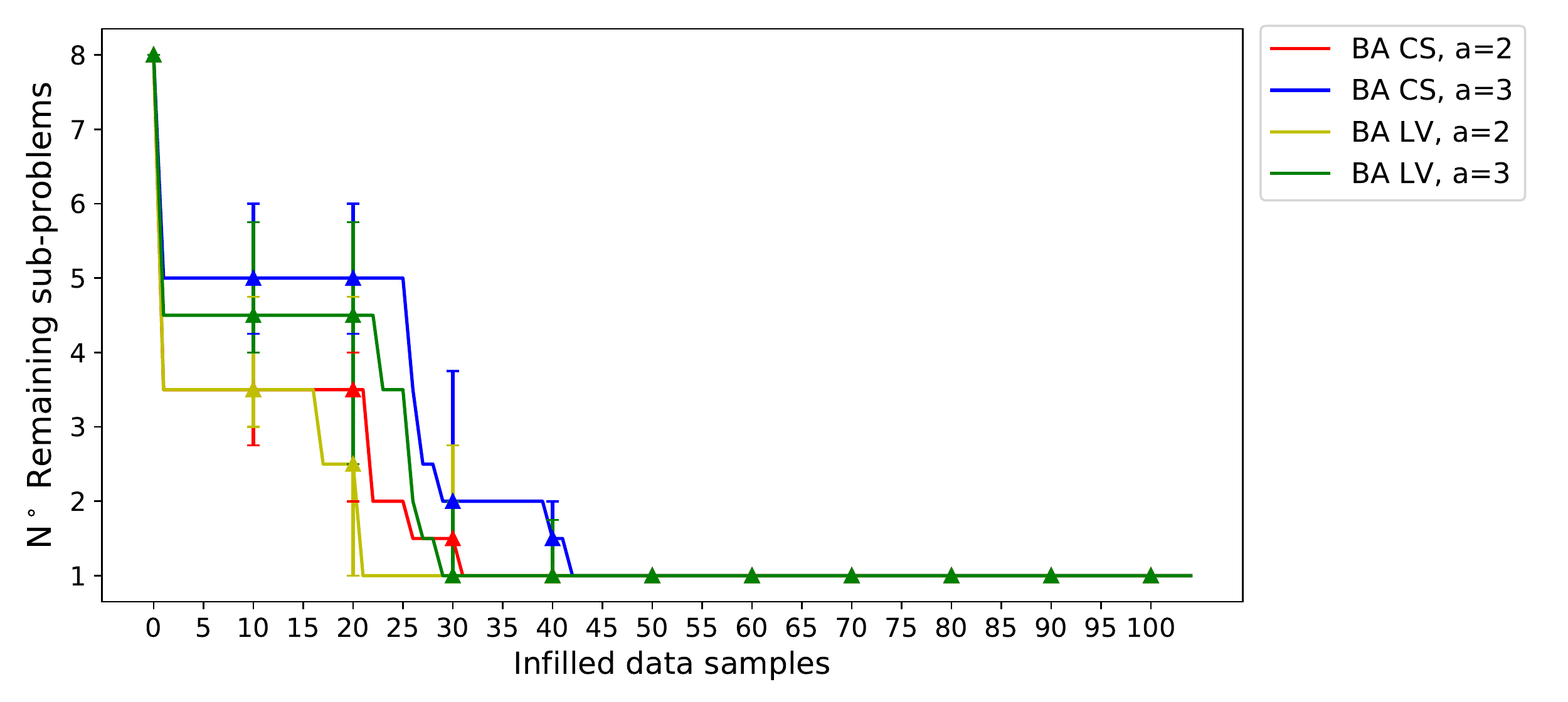}
\caption{Comparison of the number of remaining sub-problems along the optimization process for different parameterizations of the SOMVSP on the variable-size design space Goldstein function over 10 repetitions.}
\label{Res_Opt_VarDimGoldstein4}
\end{figure}

The compared  algorithms are initialized with a total data set of 104 data samples, which is equivalent to providing 2 samples for every dimension of  each of the 8 considered sub-problems. Subsequently, the optimizations are performed by relying on 104 additional function evaluations. The results obtained over 10 repetitions are provided in Figures \ref{Res_Opt_VarDimGoldstein1}, \ref{Res_Opt_VarDimGoldstein2}, \ref{Res_Opt_VarDimGoldstein3} and \ref{Res_Opt_VarDimGoldstein4}

\noindent Overall, the results show that all the variants of the proposed algorithms, namely the SOMVSP and the variable-size design space kernel based BO algorithm, provide a faster and more consistent convergence towards the considered VSDSP optimum when compared to the independent optimization of each sub-problem. For the SOMVSP, this difference can be explained by the fact that the proposed strategy allows to better focus the computational budget towards the most promising sub-problems through budget allocation and sub-problem discarding. By consequence, very few function evaluations are wasted onto the least promising sub-problems. For the variable-size design space kernel based BO algorithm, instead, the better performance can be explained by the fact that the GP models are created and trained by relying on the entirety of the available data (\ie the data set in the variable-size design space) and can therefore better exploit the available information. On the other hand, the IO of each sub-problem  relies on GP which are created by relying only on data sets specific to each sub-problem, and rely therefore on a lower amount of information. As a consequence, the resulting modeling performance of the problem functions is expected to be less accurate.

When considering the relative performance between the 2 proposed methods (and their variants), the results show a better convergence of the variable-size design space kernel based BO algorithm with respect to the SOMVSP, in terms of both convergence rate, as can be seen in Figure \ref{Res_Opt_VarDimGoldstein1}, as well as optimum value at convergence, as is shown in Figure \ref{Res_Opt_VarDimGoldstein2}. Even in the case in which the SOMVSP has identified the optimal sub-problem and discarded all the others, the variable-size design space BO may still converge faster thanks to the information  it can extrapolate from the data samples belonging to non-optimal sub-problems. 

When analyzing the results provided by the different SOMVSP variants, a difference in convergence rate can be noticed for different values of $a$ (but identical kernels). Indeed, the optimizations performed with a value of $a = 2$ provide a faster convergence when compared to the ones obtained with a value of $a = 3$. This is related to the fact that lower values of $a$ result in a more frequent discarding of sub-problems, which can then result in a better allocation of the computational resources. However, lower values of $a$ can also cause the proposed algorithm to discard the optimal sub-problems during the early stages of the optimization, when insufficient data is provided in order to make accurate choices. This phenomenon can, for instance, be noticed in the larger variance of the results obtained with the CS kernel and a value of $a=2$ when compared to the ones obtained with the CS kernel and a value of $a=3$. In order to better illustrate the effect of considering different values of $a$, the evolution of the remaining number of sub-problems along the optimization for different kernels is shown in Figure \ref{Res_Opt_VarDimGoldstein4}. As previously discussed, it can be seen that lower values of $a$ result in a faster discarding of sub-problems (if the same kernel is considered). Furthermore, it can also be noticed that relying on the LV kernel tends to result in a faster discarding of problems when compared to the CS, which is related to the different accuracy of the model of error. Finally, the figure also shows that from half of the optimization onward, all of the compared SOMVSP variants tend to have discarded all the problems but one.

When analyzing the results provided by the variable-size design space kernel based BO variants, the first noticeable result is a considerably faster convergence of the DVW variant when compared to the SPW one. Indeed, the DVW variant consistently converges around the 10th infilled point, whereas the SPW variant requires more or less half of the allocated budget in order to yield the same performance. This difference is related to the fact that the DVW kernel is defined in such a way that it can rely on the variables shared between the different sub-problems in order to better model the considered function, while the SPW kernel only computes the covariance between samples defined in different sub-problems through their dimensional variable values. This difference becomes more noticeable when the CS discrete kernel is considered, due to the fact that it relies on a single covariance value between all of the sub-problems. However, note that this is not true for the DVW kernel, as in this case the covariance between sub-problems is computed as the product between the covariances between the values of the two dimensional variables. Finally, no significant difference of performance between the variable-size design space kernel based BO variants can be noticed from the values at convergence shown in Figure \ref{Res_Opt_VarDimGoldstein3}, with the exception of a slightly larger variance for the CS based SPW kernel. This can be explained by the fact that all of these compared methods are provided with enough function evaluations in order to properly refine the incumbent solution.

\vspace{12pt}

\begin{figure}[H]
\centering
\includegraphics[width=1.05\linewidth]{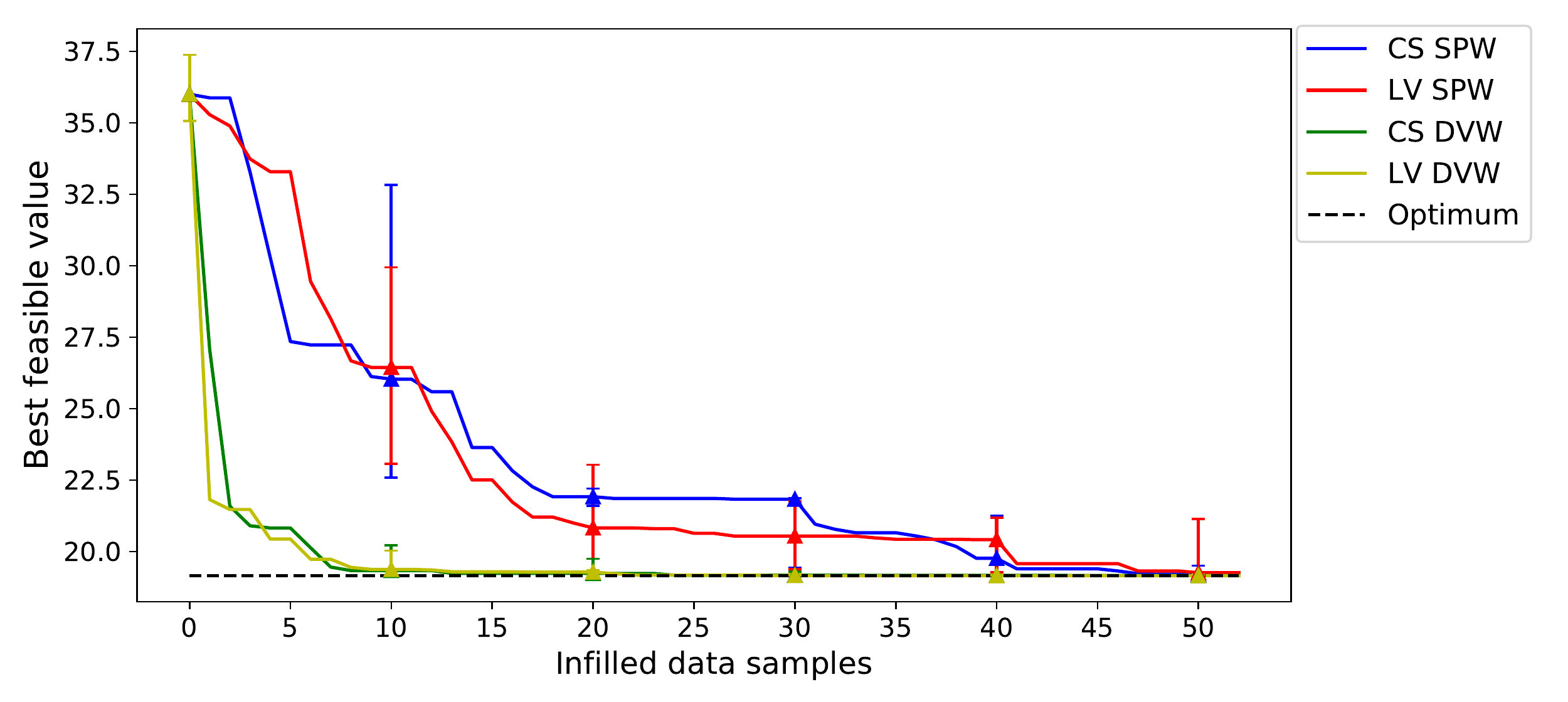}
\caption{Comparison of the convergence rate of different VSDSP optimization algorithms during the BO of the mixed-variable Goldstein function over 10 repetitions. Focus on the the variable-size design space kernel BO alternatives.}
\label{Res_Opt_VarDimGoldstein_HD_1}
\end{figure}

\begin{figure}[H]
\centering
\includegraphics[width=.7\linewidth]{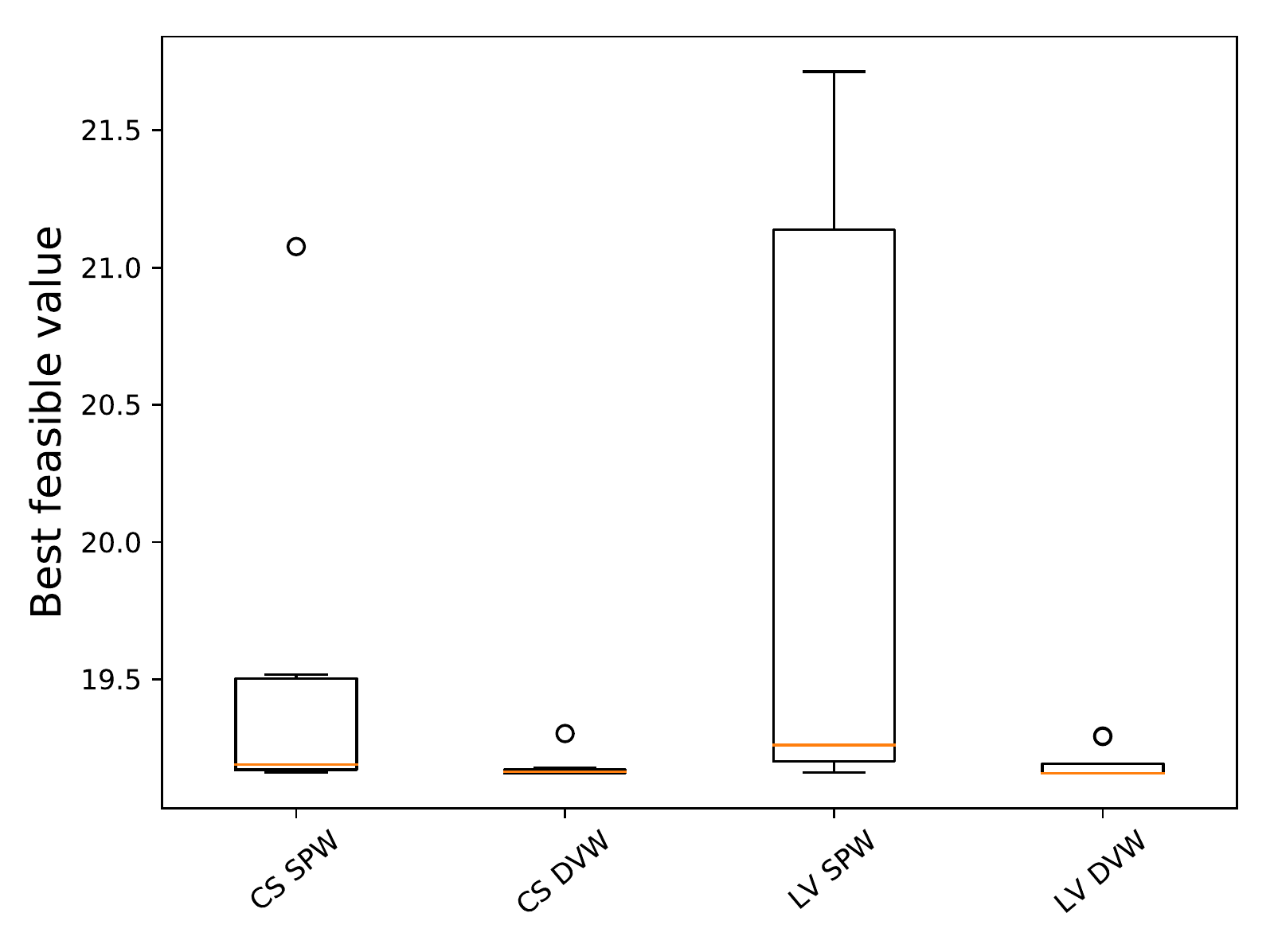}
\caption{Comparison of the convergence value of different VSDSP optimization algorithms on the variable-size design space Goldstein function over 10 repetitions. Focus on the the variable-size design space kernel BO alternatives.}
\label{Res_Opt_VarDimGoldstein_HD_2}
\end{figure}

In order to better compare the performance of the different  variable-size design space kernel based BO variants, the optimizations presented in the paragraph above are repeated by providing a smaller initial data set (\ie 52 data samples, 1 per dimension of each sub-problem) and a lower number of additional function evaluations (\ie  52 data samples). The results obtained over 10 repetitions are provided in Figures \ref{Res_Opt_VarDimGoldstein_HD_1} and \ref{Res_Opt_VarDimGoldstein_HD_2}.

As can be expected, the convergence rate obtained in this case is slower if compared to the previously considered optimizations. This can be explained by the fact that smaller initial data sets usually result in GP models characterized by a lower modeling accuracy. As a consequence, a larger number of function evaluations is necessary in order to sufficiently refine the models before being able to identify the global optimum neighborhood. Furthermore, the obtained results also show that the relative performance between the DVW and the SPW kernels remains the same when considering smaller initial data sets, as the DVW kernel provides a considerably faster convergence rate. As for the previous optimizations, this is due to the fact that the DVW kernel can rely on a larger amount of information in order to compute the covariance  between data samples characterized by different dimensional variable values.

\subsection{Variable-size design space Rosenbrock function}
The second analytical test-case that is considered in this paper is  a modified constrained variable-size design space version of the Rosenbrock function \cite{rosenbrock1960automatic}. The global VSDSP is characterized by 8 continuous variables, 3 discrete variables and 2 dimensional variables. Depending on the dimensional variable values, 4 different sub-problems can be identified, with total dimensions of 6 to 9, ranging from 4 continuous variables and 2 discrete variables to 6 continuous variables and 3 discrete variables. 

The resulting optimization problem can be defined as follows:

\begin{align}
 \min & \qquad  f(\textbf{x},\textbf{z},\textbf{w}) \\ 
\text{w.r.t.} &  \qquad \textbf{x} = \{x_1,\dots,x_8\} \ \mbox{ with } \  x_i \in [-1,0.5] \ \mbox{ for }   i = 1,3,5,7  \nonumber \\
& \qquad  \qquad \qquad \  \qquad \qquad \qquad x_i \in [0,1.5] \ \mbox{ for }   i = 2,4,6,8  \nonumber \\
 &  \qquad \textbf{z} = \{z_1,\dots,z_3\} \ \mbox{ with } \  z_i \in \{ 0,1\}  \ \mbox{ for }   i = 1,2 \nonumber \\
 &  \qquad \qquad \qquad \qquad \qquad \qquad \   z_i \in \{ 0,1,2\}  \ \mbox{ for }   i = 3 \nonumber \\
 & \qquad \textbf{w}  = \{w_1,w_2\} \ \mbox{ with } \  w_1 \in \{ 0,1\}  \ \mbox{ and }  w_2 \in \{ 0,1\}  \nonumber \\
\text{s.t.:} &  \qquad \mbox{g}(\textbf{x},\textbf{z},\textbf{w}) \leq 0  \nonumber 
\end{align}

\noindent where:

\begin{equation}
  f(\textbf{x},\textbf{z},\textbf{w}) = \begin{cases}
               f_1(x_1,x_2,x_3,x_4,z_1,z_2) \hfill   \qquad  \mbox{ if }  w_1 = 0  \mbox{ and } w_2 = 0\\
               f_2(x_1,x_2,x_5,x_6,z_1,z_2,z_3) \hfill   \qquad  \mbox{ if }  w_1 = 0  \mbox{ and } w_2 = 1\\
                f_3(x_1,x_2,x_3, x_4,x_7, x_8,z_1,z_2) \hfill   \qquad  \mbox{ if }  w_1 = 1  \mbox{ and } w_2 = 0\\
               f_4(x_1,x_2,x_5,x_6,x_7, x_8,z_1,z_2,z_3) \hfill   \qquad  \mbox{ if }  w_1 = 1  \mbox{ and } w_2 = 1
            \end{cases}
\end{equation}

\noindent and:

\begin{equation}
  g_1(\textbf{x},\textbf{z},\textbf{w}) = \begin{cases}
               g_{1_1}(x_1,x_2,x_3,x_4) \hfill   \qquad  \mbox{ if }  w_1 = 0  \mbox{ and } w_2 = 0\\
                g_{1_2}(x_1,x_2,x_5,x_6) \hfill   \qquad  \mbox{ if }  w_1 = 0  \mbox{ and } w_2 = 1\\
                 g_{1_3}(x_1,x_2,x_3, x_4,x_7, x_8) \hfill   \qquad  \mbox{ if }  w_1 = 1  \mbox{ and } w_2 = 0\\
                g_{1_4}(x_1,x_2,x_5,x_6,x_7, x_8) \hfill   \qquad  \mbox{ if }  w_1 = 1  \mbox{ and } w_2 = 1
            \end{cases}
\end{equation}

\noindent  and:

\begin{equation}
  g_2(\textbf{x},\textbf{z},\textbf{w}) = \begin{cases}
               g_{2_1}(x_1,x_2,x_3,x_4) \hfill   \qquad  \mbox{ if }  w_1 = 0  \mbox{ and } w_2 = 0\\
                 g_{2_2}(x_1,x_2,x_3, x_4,x_7, x_8) \hfill   \qquad  \mbox{ if }  w_1 = 0  \mbox{ and } w_2 = 1\\
                \mbox{Not active} \hfill   \qquad  \mbox{ if }  w_1 = 1  \mbox{ and } w_2 = 0\\
                \mbox{Not active}  \hfill   \qquad  \mbox{ if }  w_1 = 1  \mbox{ and } w_2 = 1
            \end{cases}
\end{equation}

\begin{figure}[h!]
\centering
\includegraphics[width=.5\linewidth]{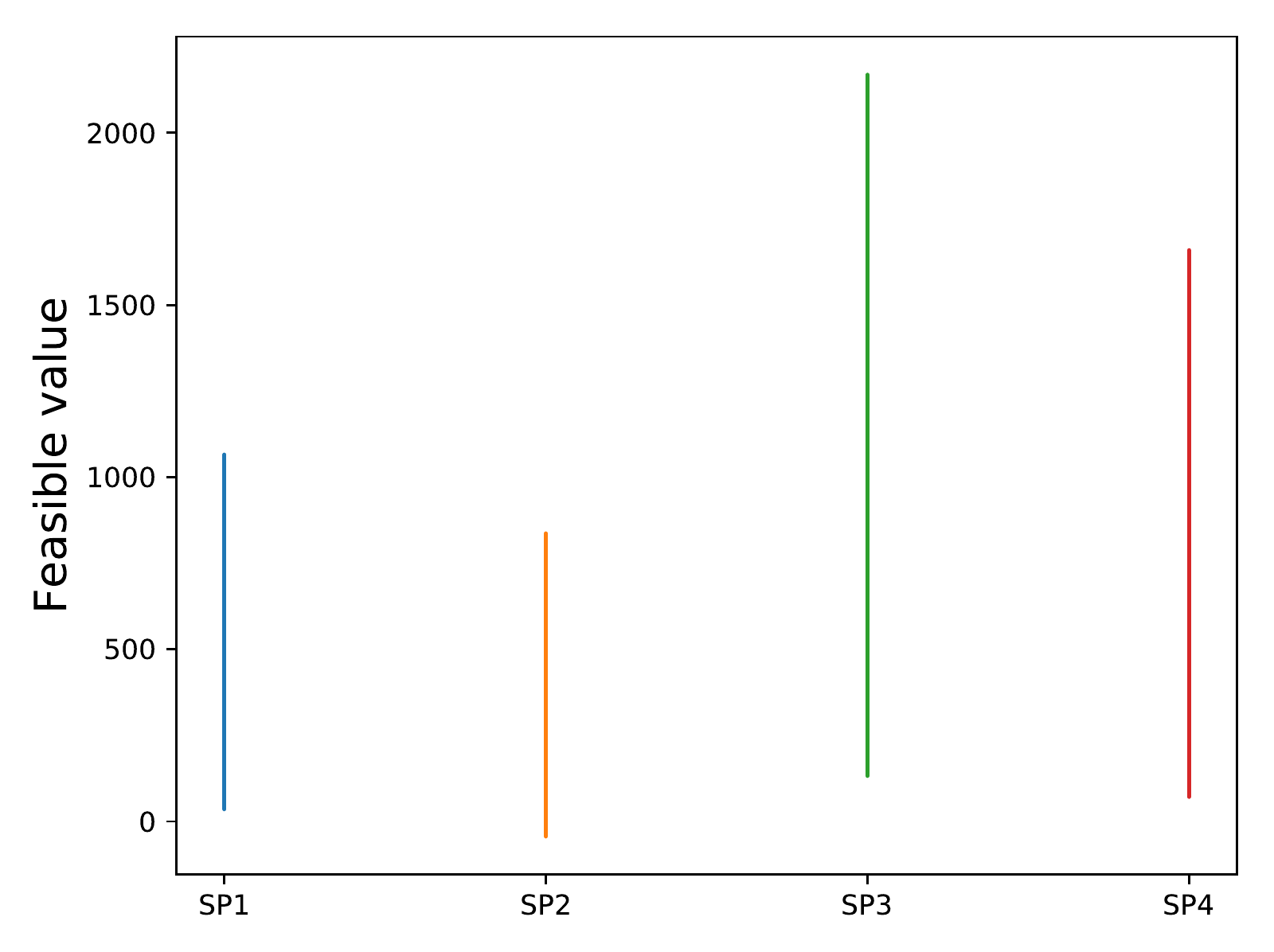}
\caption{Value range of the feasible objective function for each sub-problem of the variable-size design space Rosenbrock function.}
\label{FeasVal_Ros}
\end{figure}

As for the Goldstein test-case, the definition of the Rosenbrock objective and constraint functions can be found in Appendix B. A synthesis of the characteristics of the various sub-problems comprising the considered variable-size design space function problem is presented in Table \ref{VarDimRosSynth}. It can be noted that, compared to the previous test-case, the variable-size design space Rosenbrock function is characterized by a similar number of continuous and discrete design variables. However, due to the smaller number of levels characterizing the dimensional variables, it only presents half of the sub-problems. Furthermore, it is important to point out that this function is subject to 2 constraints, which are however not active in all sub-problems. Indeed, the second constraint is only present within sub-problems N${}^\circ$1 and 2. Finally, the value ranges of the feasible objective function for each sub-problem are shown in Figure \ref{FeasVal_Ros}. For this test-case as well, the feasible values of the different sub-problems overlap over a large part of their objective function range, thus making the identification of the optimal sub-problem challenging.
\begin{table}[h]
\centering
\begin{tabular}{|l|c|c|c|c|}
\hline
Sub-problem &  SP 1 & SP 2 & SP 3 & SP 4 \\ \hline
N$^\circ$ continuous variables & 4 & 4 & 6 & 6 \\ \hline
N$^\circ$ discrete variables & 2 & 3 & 2 & 3 \\ \hline
N$^\circ$ constraints & 2 & 2 & 1 & 1  \\  \hline \hline 
\multicolumn{5}{|c|}{Global VSDSP} \\ \hline 
 N$^\circ$ continuous variables & \multicolumn{4}{|c|}{8}  \\ \hline
N$^\circ$ discrete variables &\multicolumn{4}{|c|}{3} \\ \hline
N$^\circ$ dimensional variables &\multicolumn{4}{|c|}{2} \\ \hline
N$^\circ$ discrete categories & \multicolumn{4}{|c|}{32} \\ \hline 
N$^\circ$ constraint & \multicolumn{4}{|c|}{2} \\ \hline 
\end{tabular}
\caption{Defining characteristics of the sub-problems comprising the variable-size design space Rosenbrock function  optimization problem.}
\label{VarDimRosSynth} 
\end{table}

\begin{figure}[h!]
\centering
\includegraphics[width=1.0\linewidth]{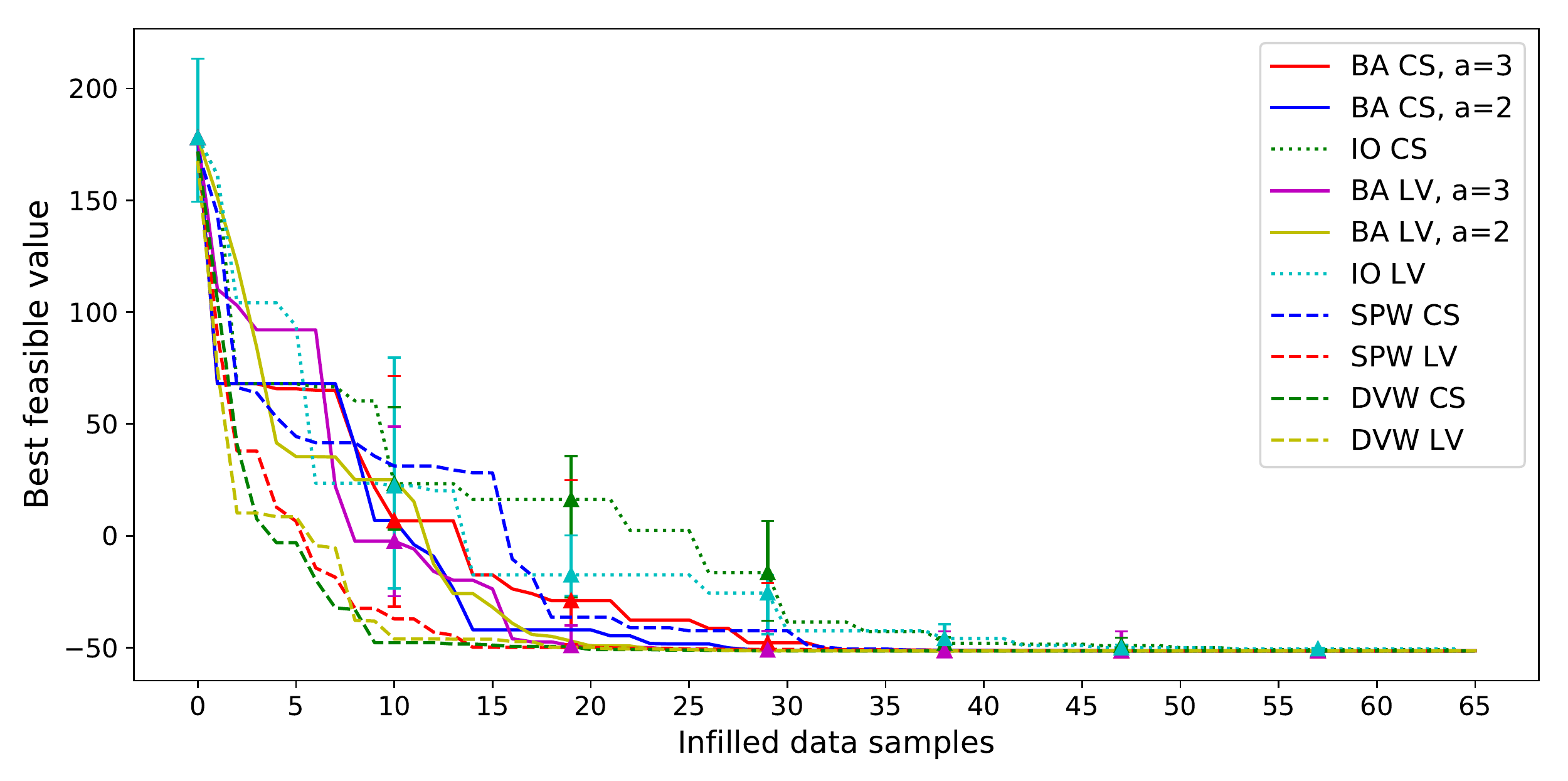}
\caption{Comparison of the convergence rate of various VSDSP optimization algorithms during the BO of the mixed-variable Rosenbrock function over 10 repetitions. The continuous lines represent the SOMVSP alternatives, the dashed lines represent the variable-size design space kernel BO alternatives and the dotted lines represent the independent optimization of each sub-problem.}
\label{Res_Opt_VarDimRos1}
\end{figure}

\begin{figure}[h!]
\centering
\includegraphics[width=.7\linewidth]{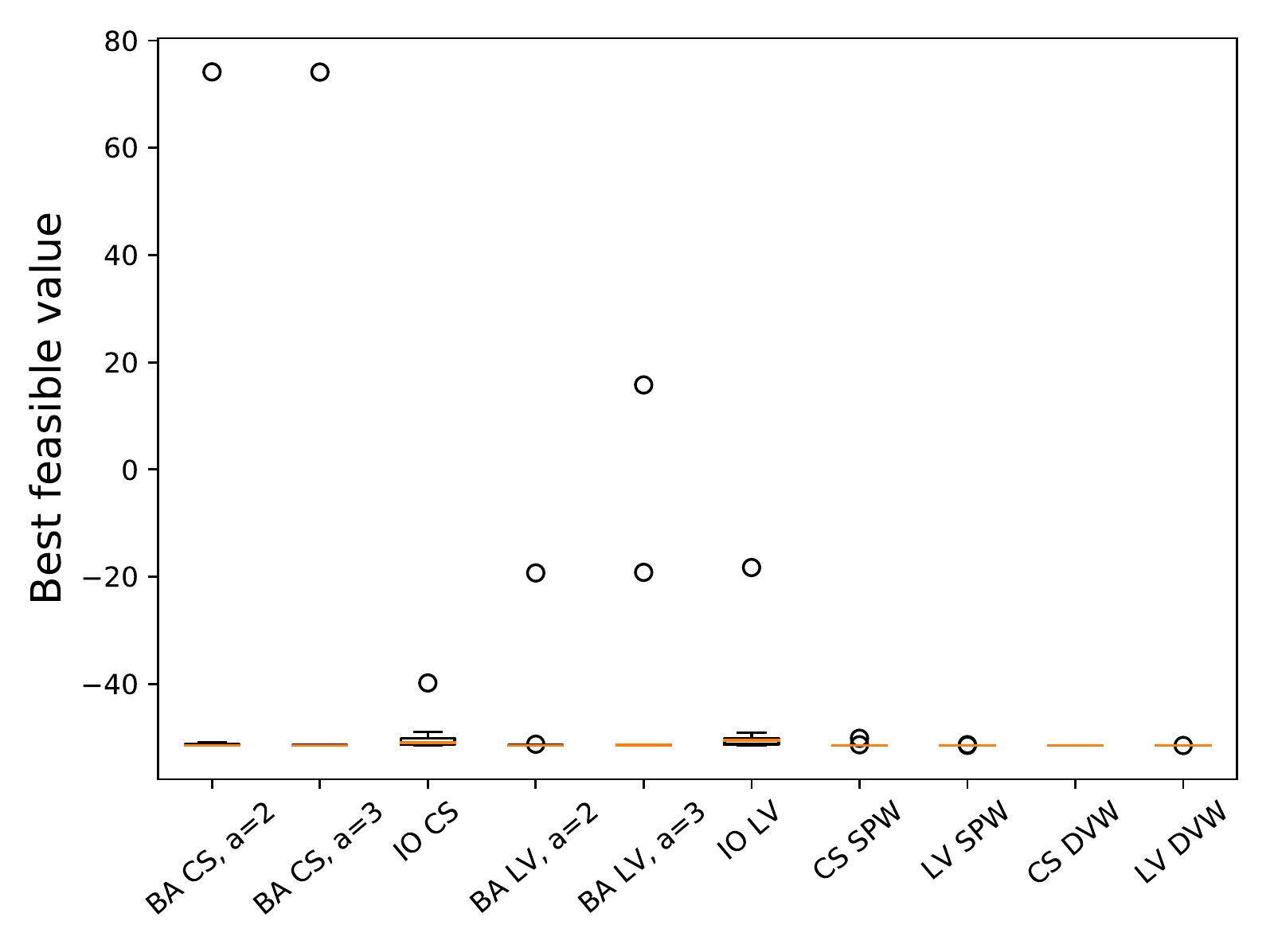}
\caption{Comparison of the convergence value of different VSDSP optimization algorithms on the variable-size design space Rosenbrock function over 10 repetitions.}
\label{Res_Opt_VarDimRos2}
\end{figure}

The compared  algorithms are initialized with a total data set of 30 data samples, which is equivalent to providing 1 samples for every dimension of  each of the 4 considered sub-problems. Subsequently, the optimizations are performed by relying on 65 additional function evaluations. The results obtained over 10 repetitions are provided in Figures \ref{Res_Opt_VarDimRos1} and  \ref{Res_Opt_VarDimRos2}.

Similarly to the variable-size design space Goldstein function, the optimization results show a faster convergence towards the problem optimum for both proposed approaches regardless of their parameterization (\ie type of kernel, SPW or DVW decomposition, value of $a$) if compared to the independent optimization of each sub-problem, thanks to the better exploitation of the limited information present within the training data set.  When analyzing the relative performance between the two proposed methods, it can be seen that also in this case the direct BO of the VSDSP through a variable-size design space kernel yields an overall faster convergence with respect to the SOMVSP. This difference is related to the fact that this approach can rely on the entirety of the data when computing the covariance between data samples, whereas the budget allocation approach performs a separate optimization of each sub-problem. Additionally, the results also show that for this test-case, the relative difference in performance between the various considered parameterizations of the compared algorithms (\ie type of kernel, SPW or DVW decomposition, value of $a$) has a relatively low influence if compared to what can be observed for the Goldstein function. This phenomenon is likely related to the smoother nature of the modeled objective and constraint functions which results in a relatively accurate modeling of the problem functions, independently from the kernel parameterization. The only exception to the previous statement is the variable-size design space kernel approach relying on the  CS SPW decomposition, which convergences considerably slower when compared to the other direct BO variants. This difference can be explained by the fact that this decomposition only relies on a single covariance value between the 4 sub-problems, which can result in an overly simplistic model. This effect is further accentuated by the fact that the SPW decomposition does not exploit the information provided by the variables which are shared between the sub-problems.

\begin{figure}[h]
\centering
\includegraphics[width=1.0\linewidth]{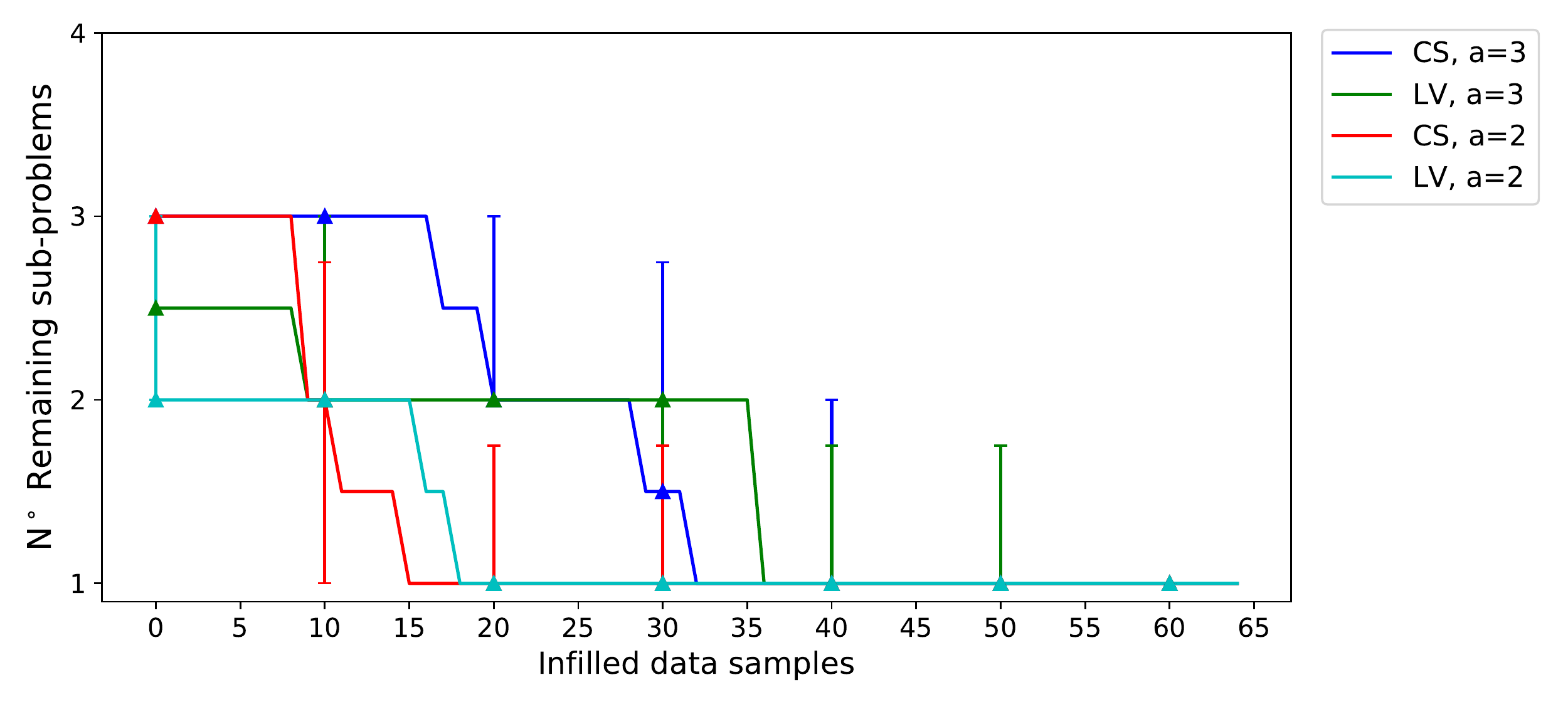}
\caption{Comparison of the number of remaining sub-problems along the optimization process for different parameterizations of the SOMVSP on the variable-size design space Rosenbrock function over 10 repetitions.}
\label{Res_Opt_VarDimRos3}
\end{figure}

Figure \ref{Res_Opt_VarDimRos2} also shows that, although both proposed methods manage to provide convergence towards the global optimum of the considered problem, the variable-size design space kernel approach manages doing so with a better consistency, whereas all of the SOMVSP variants present a few outliers, which represent optimization repetitions for which the global optimum was not found. This effect can be explained by the premature discarding of the optimal sub-problem by the SOMVSP due to the low amount of data the optimization algorithms are initialized with. This effect is better highlighted in Figure \ref{Res_Opt_VarDimRos3}, where the number of remaining sub-problems along the SOMVSP is shown. Indeed, it can be seen that on average over the 10 repetitions, one of the sub-problems is discarded at the very first iteration of the budget allocation strategy, when the amount of available data is the lowest, thus resulting in a possibly incorrect choice. Additionally, Figure \ref{Res_Opt_VarDimRos3} also shows a faster discarding of sub-problems when considering values of $a= 2$ instead of  $a= 3$, regardless of the kernel parameterization, which is consistent with the functioning principle of the algorithm. Indeed, higher values of $a$ tend to result in a safer, but slower, discarding process. Finally, the optimization results also show that similarly to the Goldstein test-case, the CS kernel-based SOMVSP tends to yield a faster discarding of sub-problems if compared to the LV kernel. This phenomenon is related to a difference in the error modeling accuracy.

\subsection{Multi-stage launch vehicle architecture optimization}
The second variable-size design space problem which is considered in this paper is the preliminary optimization of a multi-stage launch vehicle architecture. This design problem requires to simultaneously determine the optimal number of stages characterizing the system (\ie 2 or 3) and determine the most suitable type of propulsion (\ie solid or liquid) for each stage. Given that each propulsive alternative is characterized by different continuous and discrete design variables as well as different constraints, the resulting optimization problem presents a variable-size design space. The objective function allowing to assess the performance of the system is defined as the Gross Lift-Off Weight (GLOW), which is computed as the sum of the payload mass $M_{PL}$, as well as the dry mass ($M_d$) and propellant mass ($M_{prop}$) of each stage:
\begin{equation}
GLOW = M_{PL} + \sum_{i = 1}^{n_{stages}} \left(M_{d_i} + M_{prop_i}\right)
\end{equation}
The target mission for which the launch vehicle is designed is the injection of a 500 kg payload into an 800 km Low Earth Orbit (LEO) \cite{wertz2001mission}. In order to ensure that this orbit is reached, a number of constraints on the total velocity increment provided by the engines as well as  on the thrust-to-weight ratio of each stage must be taken into account.

\subsubsection{Liquid propulsion}
Each liquid propulsion stage is characterized by 2 continuous design variables, namely the stage specific propellant mass $M_{prop}$ and the thrust value $T$, as well as one discrete design variable characterizing the type of engine $Type_{eng}$ to be included in the design. Additionally, if the considered liquid propulsion stage is also the first stage of the launch vehicle, an additional discrete design variable representing the number of engines $N_{e}$ must be considered. The continuous and discrete design variables characterizing each liquid propulsion stage are detailed in Table \ref{VarLiquidStage}. Note that the continuous variable bounds vary as a function of the position of the considered architecture as well as on the overall system architecture.

\begin{table}[h]
\centering
\footnotesize
\begin{tabular}{|l|l|c|c|c|}
\hline
Variable & Nature & Levels \\ \hline
$T$ - Engine thrust [kN] & continuous & [-] \\ \hline
$M_{prop}$ – Propellant mass  [kg] & continuous & [-]\\ \hline
$Type_{eng}$ - Type of engine & discrete  &  type 1, type 2, type 3\\ \hline
$N_{e}$ - Number of engines & discrete  & 1,2  \\ \hline
\end{tabular}
\caption{Variables characterizing each liquid propulsion stage}
\label{VarLiquidStage}
\end{table}

\begin{figure}[h!]
\centering
\includegraphics[width=.65\linewidth]{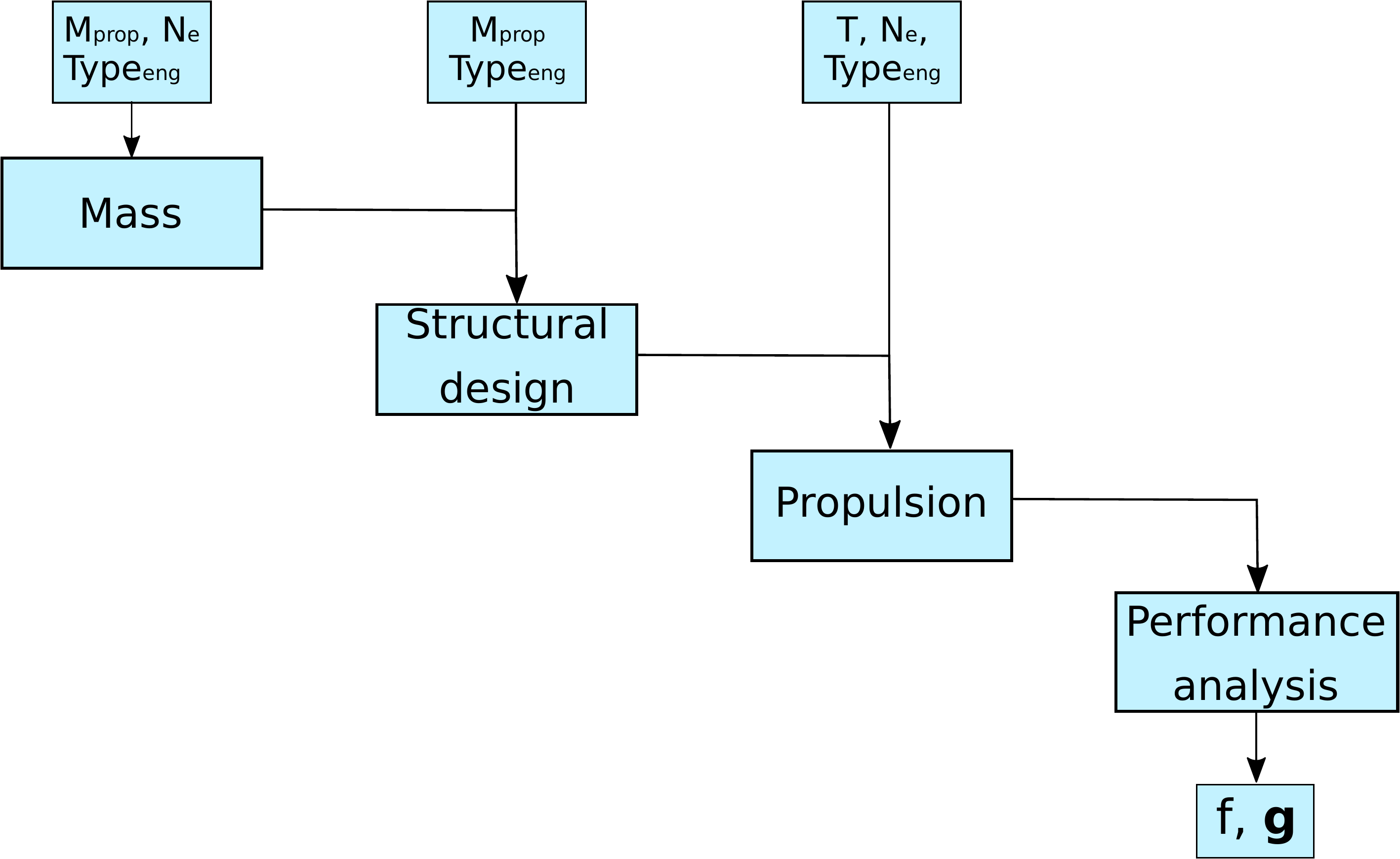}
\caption{Schematic MDO representation of the liquid propulsion stage.}
\label{Liquidprop}
\end{figure}

A schematic representation of the dependencies between the different disciplines characterizing the liquid propulsion stage as well as the continuous and discrete design variables they depend on is provided in Figure \ref{Liquidprop}.

\subsubsection{Solid propulsion}
Each solid propulsion stage is characterized as a function of 4 continuous design variables, namely the nozzle throat diameter $D_t$, the nozzle exit diameter $D_e$, the combustion chamber pressure $P_{comb}$ and the propellant mass $M_{prop}$. Additionally, the solid propulsion stage also depends on three discrete design variables: the type of propellant $Type_{prop}$, the type of material $Type_{mat}$ and the type of engine $ Type_{eng}$. Furthermore, if the first stage is considered, a discrete variable $N_b$ characterizing the number of boosters attached to the stage is also taken into account. The continuous and discrete design variables characterizing each solid propulsion stage are detailed in Table \ref{VarSolidStage}. Note that the continuous variable bounds vary as a function of the position of the considered architecture as well as on the overall system architecture.

\begin{table}[h]
\centering
\footnotesize
\begin{tabular}{|l|l|c|c|c|}
\hline
Variable & Nature & Levels \\ \hline
$D_t$ - Nozzle throat diameter [m] & continuous & [-] \\ \hline
$D_e$ - Nozzle exit diameter [m] & continuous   & [-]\\ \hline
$P_{comb}$ – Chamber pressure  [bar]  & continuous  & [-]\\ \hline
$M_{prop}$ – Propellant mass  [kg] & continuous & [-]\\ \hline
\multirow{ 2}{*}{$Type_{prop}$ - Type of propellant} & \multirow{ 2}{*}{discrete} & Butalite, Butalane, \\ 
& &  Nitramite, pAIM-120 \\ \hline
$Type_{mat}$ - Type of material & discrete  & Aluminium, Steel \\ \hline
$Type_{eng}$ - Type of engine & discrete  &  type 1, type 2, type 3\\ \hline
$N_{b}$ - Number of booster & discrete  & 1,\dots,8  \\ \hline
\end{tabular}
\caption{Variables characterizing each solid propulsion stage}
\label{VarSolidStage}
\end{table}

\begin{figure}[h!]
\centering
\includegraphics[width=.8\linewidth]{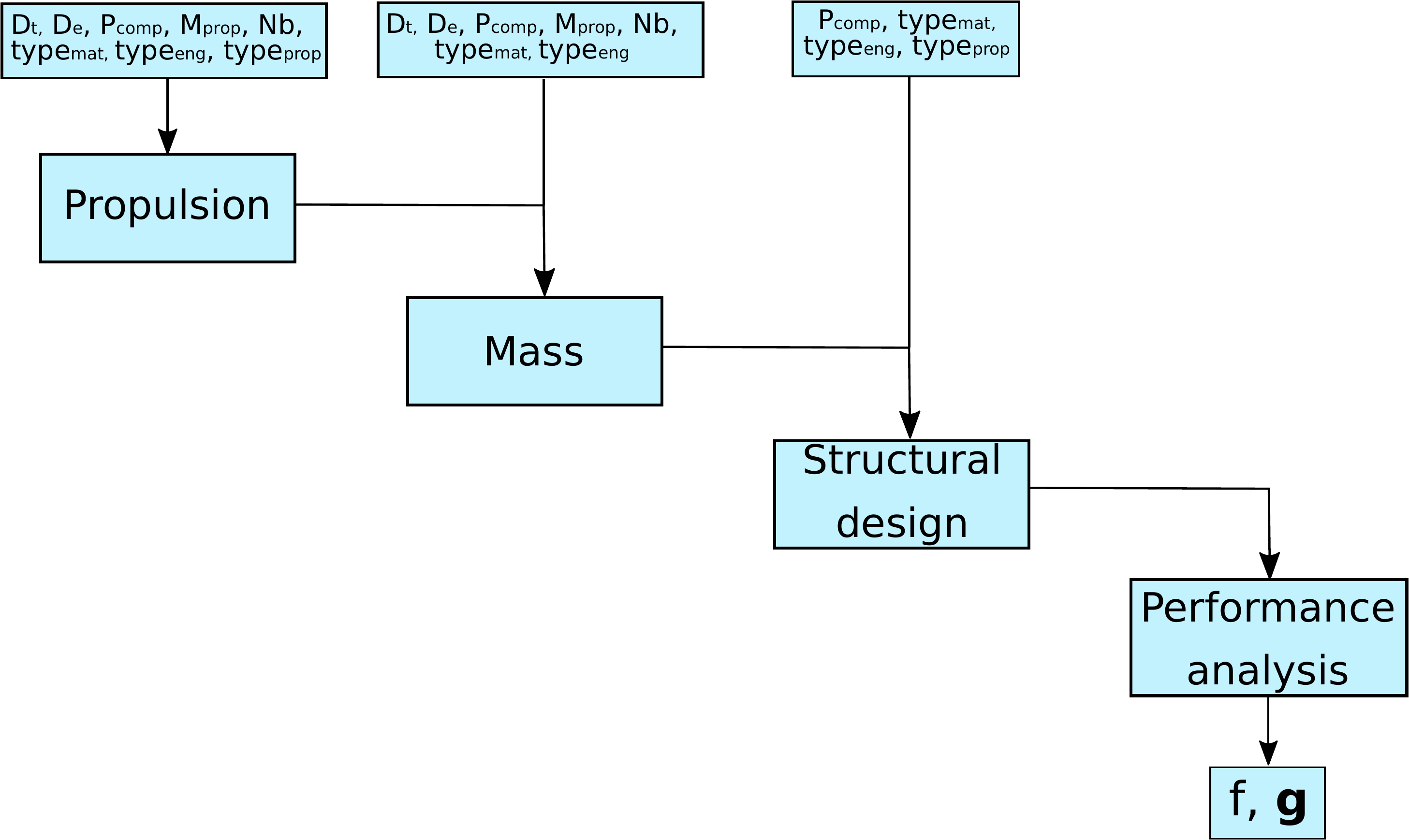}
\caption{Schematic MDO representation of the solid propulsion stage.}
\label{solidprop}
\end{figure}

A schematic representation of the dependencies between the different disciplines characterizing the solid propulsion stage as well as the continuous and discrete design variables they depend on is provided in Figure \ref{solidprop}

\subsubsection{Variable-size design space problem formulation}
\begin{figure}[h!]
\centering
\includegraphics[width=.40\linewidth]{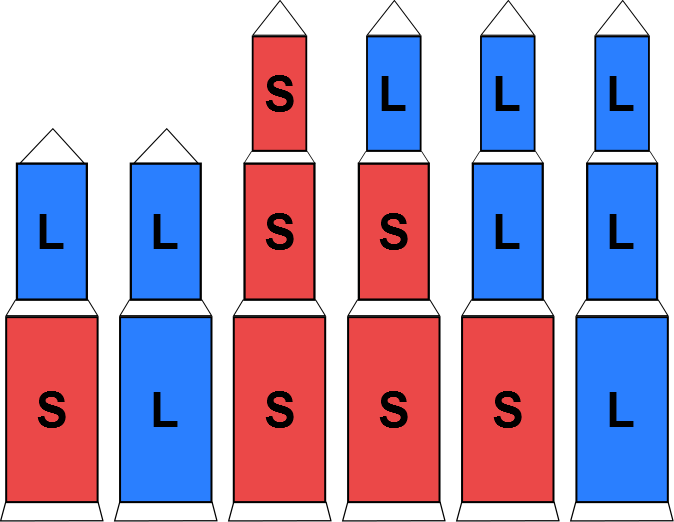}
\caption{Considered launch vehicle architectures (S: Solid propulsion, L: Liquid propulsion).}
\label{Architecture}
\end{figure}

The multi-stage launch vehicle architecture optimization described above can be formulated under the form of a variable-size design space problem by considering 3 dimensional variables, each one representing the type of propulsion to be included in one of the 3 considered stages. In order to simplify the problem, the unfeasible configurations are not taken into account, thus reducing the total number of sub-problems from 12 to 6. The remaining configurations are schematically represented in Figure \ref{Architecture}.

The resulting dimensional variables $w_1, w_2, w_3$ are respectively characterized by 2, 2 and 3 levels. $w_1$ and $w_2$ have the purpose of determining whether the first and second stage are characterized by liquid or solid propulsion, whereas $w_3$ determines whether the first stage is characterized by solid propulsion, liquid propulsion (\ie 3 stage architecture) or whether it is not included in the architecture  (\ie 2 stage architecture).

\begin{figure}[h!]
\centering
\includegraphics[width=.55\linewidth]{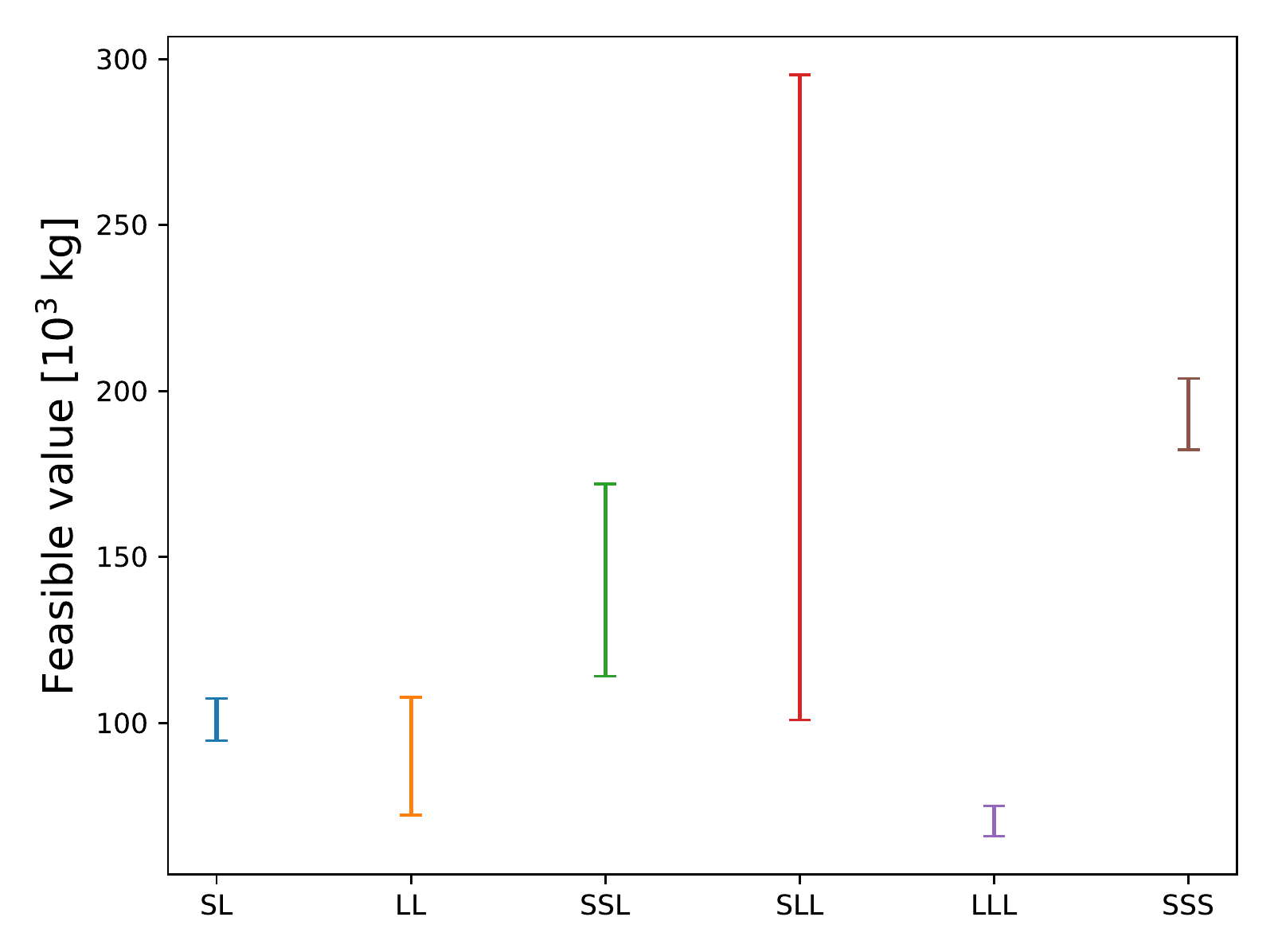}
\caption{Value range of the feasible objective function for each sub-problem of the multi-stage launch vehicle architecture design problem.}
\label{FeasVal_Launcher}
\end{figure}

The resulting VSDSP is characterized by a total of 18 continuous variables, 14 discrete variables and 3 dimensional variables, thus resulting in 29136 discrete categories. Furthermore, the problem problem is subject to 19 constraints. It can be formulated as follows:
\begin{align}
 \min & \qquad  GLOW(S_1,L_1,S_2,L_2,S_3,L_3,w_1,w_2,w_3) \\ 
\text{w.r.t.} &  \qquad S_1,L_1,S_2,L_2,S_3,L_3,w_1,w_2,w_3 \nonumber \\
\text{s.t.:} & \qquad g_{\Delta V}(S_1,L_1,S_2,L_2,S_3,L_3,w_1,w_2,w_3) \leq 0   \nonumber \\
& \qquad g_{TW_{i}}(S_1,L_1,S_2,L_2,S_3,L_3,w_1,w_2,w_3) \leq 0  \quad \mbox{for i } =1,2,3  \nonumber \\
& \qquad g_{g_{1_i}}(S_1,L_1,S_2,L_2,S_3,L_3,w_1,w_2,w_3) \leq 0  \quad \mbox{for i } =1,2,3  \nonumber \\
& \qquad g_{g_{2_i}}(S_1,L_1,S_2,L_2,S_3,L_3,w_1,w_2,w_3) \leq 0  \quad \mbox{for i } =1,2,3  \nonumber \\
& \qquad g_{g_{3_i}}(S_1,L_1,S_2,L_2,S_3,L_3,w_1,w_2,w_3) \leq 0  \quad \mbox{for i } =1,2,3  \nonumber \\
& \qquad g_{g_{4i}}(S_1,L_1,S_2,L_2,S_3,L_3,w_1,w_2,w_3) \leq 0  \quad \mbox{for i } =1,2,3  \nonumber \\
& \qquad g_{e_{i}}(S_1,L_1,S_2,L_2,S_3,L_3,w_1,w_2,w_3) \leq 0  \quad \mbox{for i } =1,2,3  \nonumber
\end{align}
where $S_n$ and $L_n$ respectively represent the design variables associated to the $n$-th solid and liquid propulsion stages, respectively. It is important to note that some constraints, such as the thrust-to-weight ratio ($g_{TW}$) and the geometrical ($g_{g}$) constraints, must be separately considered for each stage that is included in the launch vehicle architecture. A synthesis of the characteristics of the six sub-problems comprising the considered multi-stage launch vehicle architecture design problem is presented in Table \ref{VarDimLaunchTab}. Furthermore, the value range of the feasible objective function for each sub-problem is shown in Figure \ref{FeasVal_Launcher}. Differently than for the variable-size design space Goldstein function, not all the feasible objective function ranges of the different sub-problems overlap. For instance, the best feasible performance value of the triple solid propulsion stage (SSS) results considerably larger than any feasible performance of the triple liquid propulsion stage (LLL). It is expected that the optimization algorithms should mostly explore the areas of the design space associated to the sub-problems SL, LL and LLL.

\begin{table}[h]
\centering
\begin{tabular}{|l|c|c|c|c|c|c|}
\hline
Sub-problem &  SL & LL & SSS & SSL & SLL & LLL \\ \hline
N$^\circ$ continuous variables &6 & 4 & 12 & 10 & 8 & 6  \\ \hline
N$^\circ$ discrete variables & 5 & 3 & 10 & 8 & 6 & 4  \\ \hline
N$^\circ$ discrete categories & 48 & 32 & 27648 & 1152 & 192 & 64  \\ \hline
N$^\circ$ constraints & 8 & 3 & 19 & 14 & 9 & 4 \\  \hline \hline 
\multicolumn{7}{|c|}{Global VSDSP} \\ \hline 
 N$^\circ$ continuous variables & \multicolumn{6}{|c|}{18}  \\ \hline
N$^\circ$ discrete variables &\multicolumn{6}{|c|}{14} \\ \hline
N$^\circ$ dimensional variables &\multicolumn{6}{|c|}{3} \\ \hline
N$^\circ$ discrete categories & \multicolumn{6}{|c|}{29136} \\ \hline 
N$^\circ$ constraints & \multicolumn{6}{|c|}{19} \\ \hline 
\end{tabular}
\caption{Defining characteristics of the sub-problems comprising the multi-stage launch vehicle architecture optimization problem. (S: Solid propulsion, L: Liquid propulsion)}
\label{VarDimLaunchTab} 
\end{table}

\subsubsection{Optimization results}
\begin{figure}[h]
\centering
\includegraphics[width=1\linewidth]{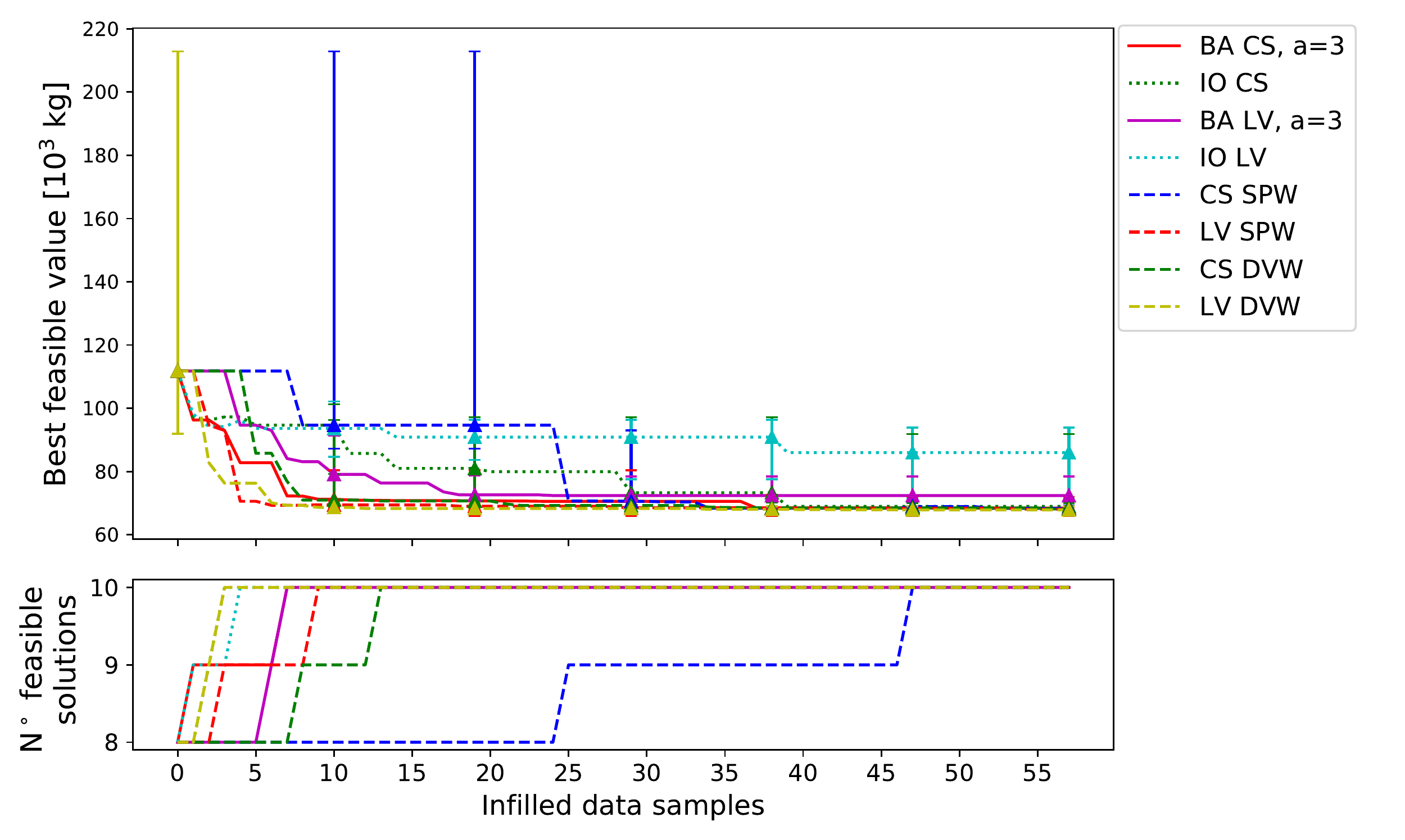}
\caption{Comparison of the convergence rate of various discrete kernels during the BO of the multi-stage launch vehicle architecture design over 10 repetitions. The bottom part of the plots represents the number of repetitions which have found a feasible solution at a given iteration.}
\label{Res_Opt_VarDimLauncher}
\end{figure}

\begin{figure}[h]
\centering
\includegraphics[width=.7\linewidth]{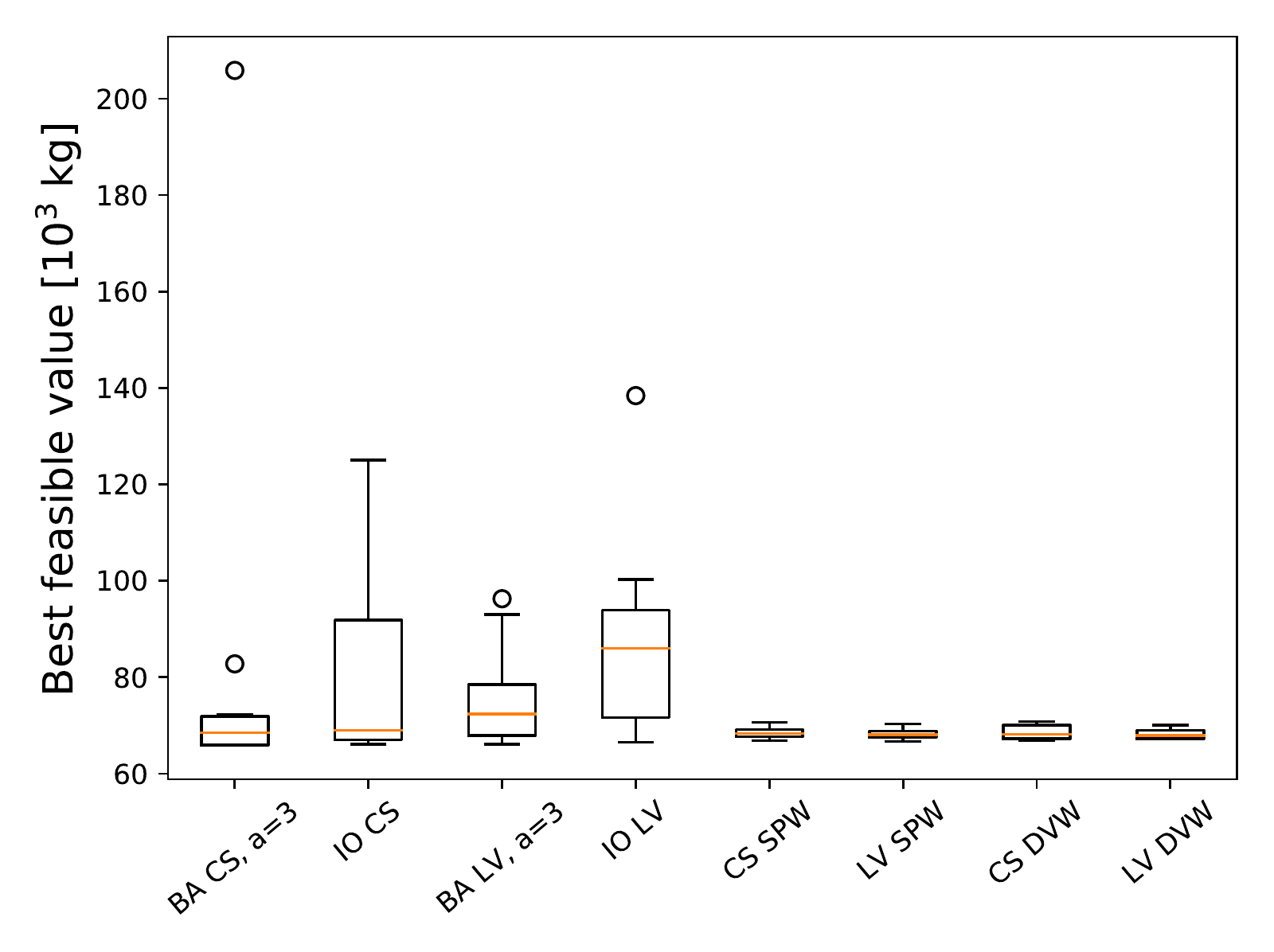}
\caption{Comparison of the convergence value of different VSDSP optimization algorithms for the BO of the multi-stage launch vehicle architecture design over 10 repetitions.}
\label{Res_Opt_VarDimLauncher2}
\end{figure}

The compared  algorithms are initialized with a total data set of 122 data samples. Subsequently, the optimizations are performed by relying on 58 additional function evaluations. Please note that because of the complexity of the considered problem, the SOMVSP methods are only applied with  values of $a$ equal to 3, in order to reduce the chances of having optimal sub-problems discarded.  The results obtained over 10 repetitions are provided in Figures \ref{Res_Opt_VarDimLauncher} and \ref{Res_Opt_VarDimLauncher2}. 

\begin{figure}[h!]
\centering
\includegraphics[width=1\linewidth]{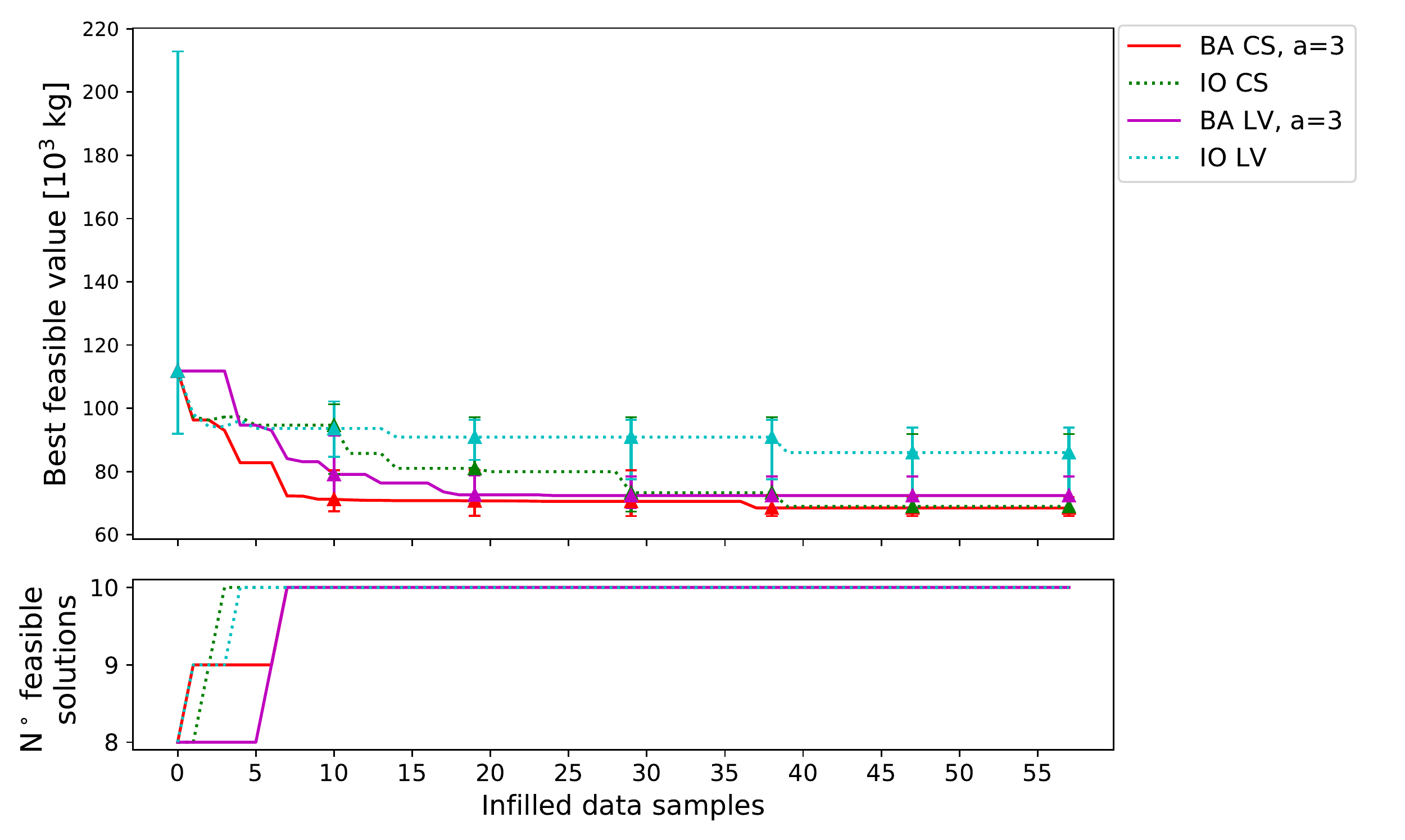}
\caption{Comparison of the convergence rate of different VSDSP optimization algorithms during the BO of the multi-stage launch vehicle architecture design over 10 repetitions. Focus on the SOMVSP methods.  The bottom part of the plots represents the number of repetitions which have found a feasible solution at a given iteration.}
\label{Res_Opt_VarDimLauncher5}
\end{figure}

\begin{figure}[h!]
\centering
\includegraphics[width=1\linewidth]{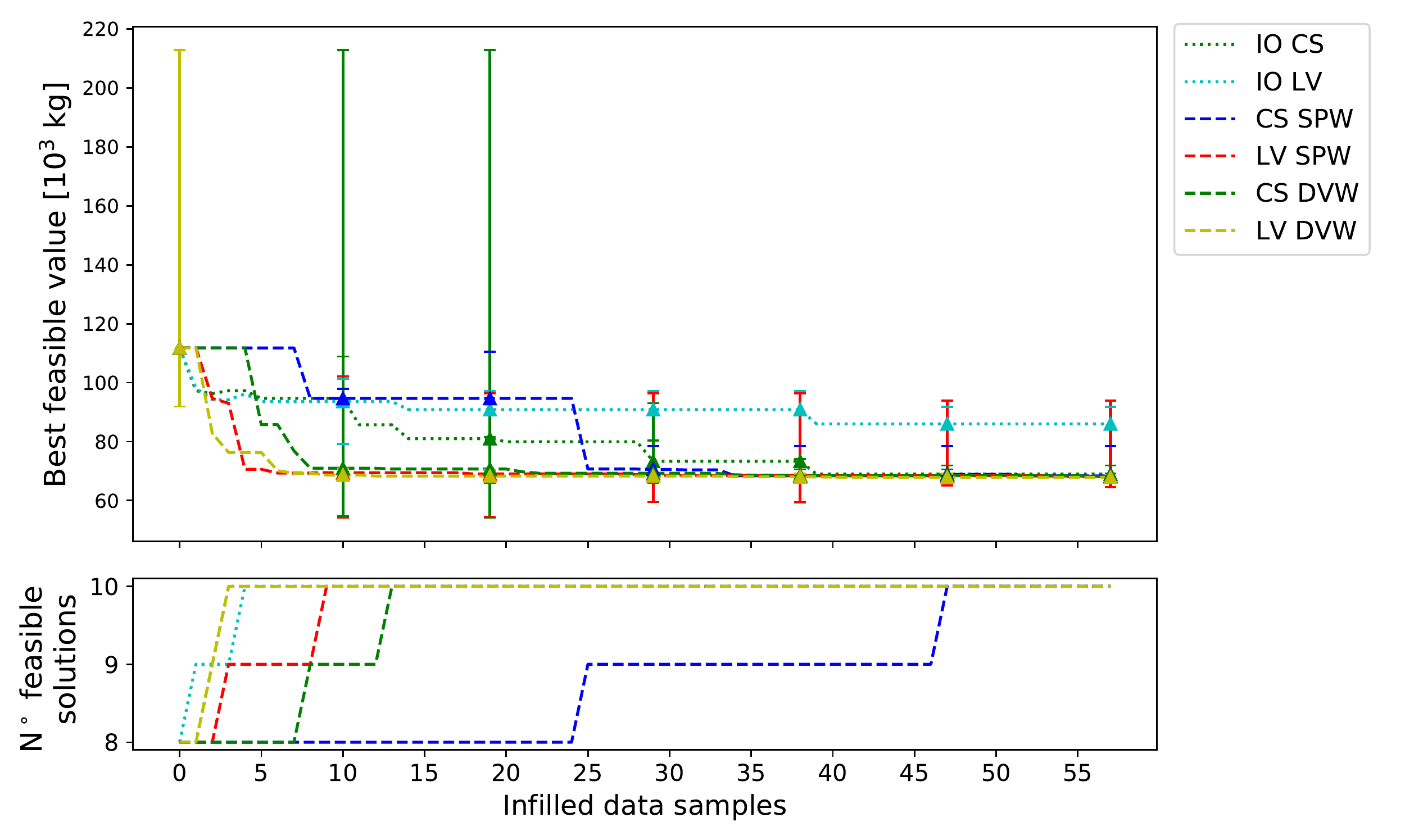}
\caption{Comparison of the convergence rate of different VSDSP optimization algorithms during the BO of the multi-stage launch vehicle architecture design over 10 repetitions. Focus on the variable-size design space kernel methods.  The bottom part of the plots represents the number of repetitions which have found a feasible solution at a given iteration.}
\label{Res_Opt_VarDimLauncher6}
\end{figure}

A global analysis of the presented results shows that both the SOMVSP methods and the variable-size design space kernel based methods provide better optimization results when compared to the independent optimization of each sub-problem. This difference can be seen in terms of convergence speed as well as final optimal value obtained at the end of the optimization process. Indeed, when considering the results obtained at the end of the optimization process, it can be noticed that the independent optimization of each sub-problem is often not sufficient in order to identify the neighborhood of the global optimum, whereas for most of the repetitions, the proposed methods are able to do so. This can be explained by the fact that the proposed methods allow to better exploit the information provided by the initial data set in order to identify the most promising areas of the design space, as well as providing a more efficient and focused use of the available computational budget.
In order to better compare the performance of the proposed methods, the convergence rate obtained with the SOMVSP and with the variable-size design space kernel BO are separately presented in Figures \ref{Res_Opt_VarDimLauncher5} and \ref{Res_Opt_VarDimLauncher6}.

Similarly to what is obtained for the variable-size design space Goldstein test-case, Figure \ref{Res_Opt_VarDimLauncher6} shows a better performance of the Dimensional Variable-Wise approach when compared to the Sub-Problem-Wise approach in terms of convergence speed. This difference is due to the fact that the DVW kernel enables to exploit a larger amount of information when computing the covariance between samples which belong to different sub-problems. Furthermore, a slightly better performance of the LV kernel with respect to the CS one can be identified when considering the variable-size design space kernel approach, as it allows to model more complex trends by relying on a larger number of hyperparameters. An opposite trend can instead be identified when analyzing the results obtained with the SOMVSP method  variants as well as for the independent optimization of each sub-problem. Indeed, for these approaches the CS kernel yields a faster and more consistent convergence if compared to the LV kernel. This can be explained by the fact that both approaches rely on the independent optimization of each sub-problem, which is performed with a considerably smaller amount of data. This makes the training of the LV hyperparameters more challenging, thus resulting in less accurate surrogate models.

\begin{figure}[h!]
\centering
\includegraphics[width=.85\linewidth]{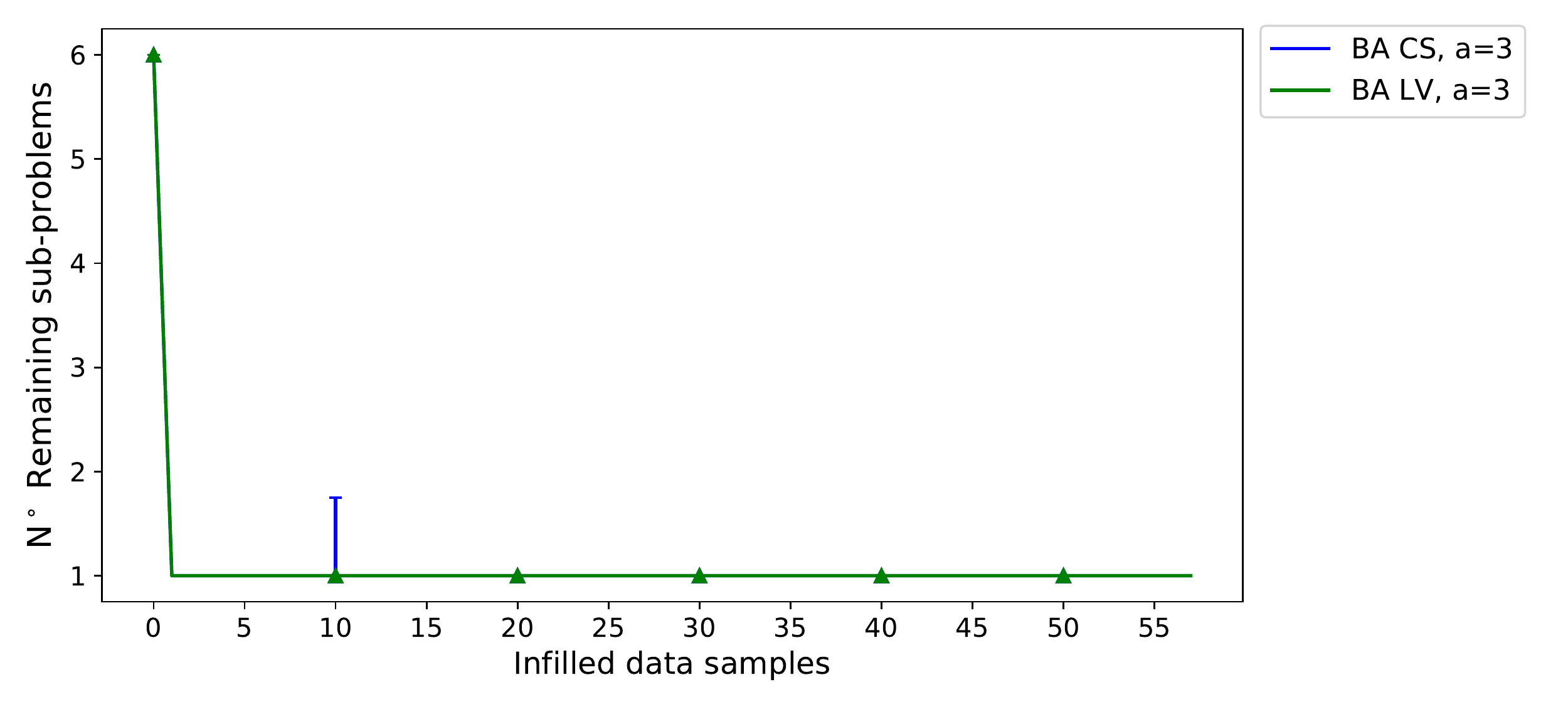}
\caption{Comparison of the remaining number of sub-problems along the optimization process for different SOMVSP algorithms during the BO of the multi-stage launch vehicle architecture design over 10 repetitions.}
\label{Res_Opt_VarDimLauncher3}
\end{figure}

A second noticeable difference between the two families of VSDSP optimization methods is related to the speed at which the the compared algorithms are able to identify the feasible areas of the search space. Indeed, it can be seen that for the repetitions which are not initialized with a feasible value within the data set, the SOMVSP methods are able to identify the feasible domain more rapidly if compared to the variable-size design space kernel based algorithms. This can be better explained by analyzing the number of remaining sub-problems along the SOMVSP process, as shown in Figure \ref{Res_Opt_VarDimLauncher3}.
It can be seen that for both the considered SOMVSP variants, all of the sub-problems but one are discarded at the first iteration of nearly every repetition. Subsequently, having to deal with a single mixed-variable problem rather than with the global VSDSP allows the SOMVSP to more easily identify the feasible areas of the design space, thus explaining the faster identification of the feasible areas of the design space.

\begin{figure}[h!]
\centering
\includegraphics[width=.75\linewidth]{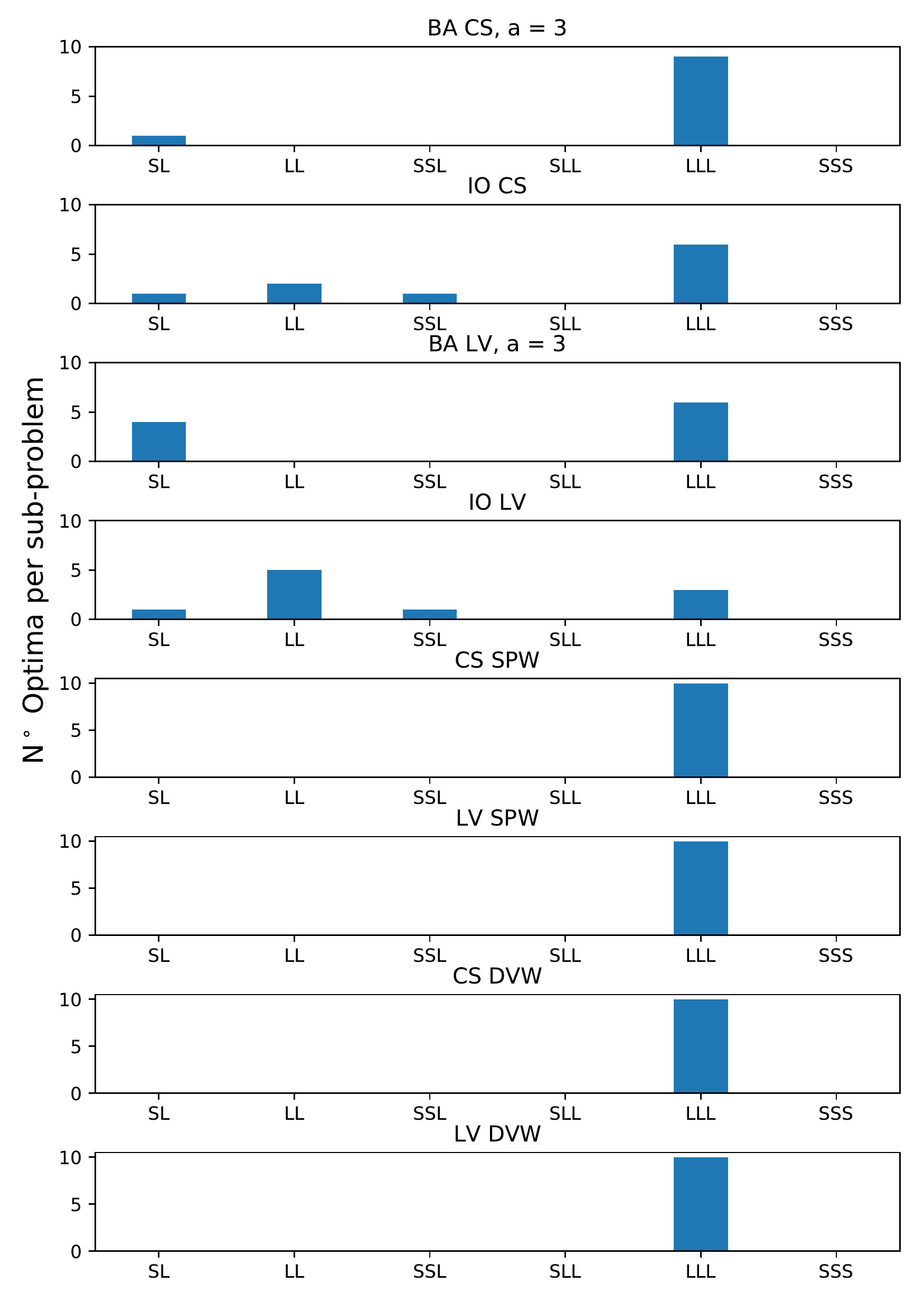}
\caption{Sub-problems towards which the compared VSDSP algorithms convergence over the 10 repetitions of the multi-stage launch vehicle architecture design optimization.}
\label{Res_Opt_VarDimLauncher4}
\end{figure}

The obtained results show that two families of VSDSP optimization algorithms yield similar convergence speeds, with a slight better performance for variable-size design space kernel based BO. However, Figure \ref{Res_Opt_VarDimLauncher2} shows that the SOMVSP methods provide a less robust convergence towards the actual optimum of the considered problem with respect to the initial data set. This can be explained with the help of Figure \ref{Res_Opt_VarDimLauncher4}, in which the sub-problems towards which the different algorithms converge over the 10 repetitions are shown. It can be seen that the variable-size design space kernel based BO methods consistently convergence towards the triple liquid propulsion sub-problem, which contains the global optimum (see Figure \ref{FeasVal_Launcher}). The SOMVSP methods, instead, also happen to converge towards non optimal sub-problems, such as the 2-stage solid/liquid propulsion architecture (SL) and the 2-stage liquid propulsion architecture (LL).

By combining the information provided by Figure  \ref{Res_Opt_VarDimLauncher3} and Figure \ref{Res_Opt_VarDimLauncher4}, it can be deduced that the SOMVSP method variants are not provided with sufficient data samples at the beginning of the optimization process. As a consequence, the sub-problems surrogate models are not accurate enough and the optimal sub-problems are discarded during the initial phases of the process, thus resulting in optimization runs which do not converge to the global optimum. A larger size initial  data set could theoretically improve the robustness of the proposed method with  and ensure a convergence towards the correct sub-problem. However, this would be incoherent with the scope of the paper, which is related to the optimization of computationally intensive problems. Alternatively, the maximum accepted EV value ($t$) for the challenging constraints could manually be tuned in order to allow for a larger violation tolerance at the beginning of the optimization process, and subsequently be reduced once sufficient data is provided. However, this would require prior information on the considered problem and time-consuming tuning in order to ensure a consistent and fast convergence of the optimization process.

\subsection{Result synthesis}
The results obtained on the considered test-cases show that the two alternative methods proposed in order to solve VSDSP provide overall a considerably faster and more consistent convergence towards the problem optimum when compared to the independent optimization of each sub-problem. This difference is due to the fact that the proposed methods allow to  exploit more efficiently the information provided by the initial data set, and are therefore able to focus the available computational budget towards the areas of interest of the design space. Additionally, the results show a slightly better performance of the variable-size design space kernel approach when compared to the SOMVSP one, both in terms of convergence speed and robustness with respect to the initial data set. Again, this can be explained by the fact that the variable-size design space kernel based  BO can exploit the entirety of the data set in order to identify the most promising areas of the design space, whereas the GP models used by the SOMVSP are independently defined over the various sub-problems and rely therefore on smaller data sets. 

The proposed variable-size design space kernel based methods (\ie SPW and VDW) present a fast convergence speed as well as a good robustness with respect to the initial data set for both considered test-cases. Furthermore, the obtained results also show that between the two variants, the dimensional variable-wise decomposition yields a faster convergence speed due to the fact that a larger amount of information can be exploited when computing the covariance between samples which belong to different sub-problems. The sub-problem-wise decomposition approach, instead, only relies on a kernel defined with respect to the (scalar) dimensional variable levels in order to compute the same type of covariance value. 

The proposed SOMVSP method shows a promising convergence speed as well, if compared to the independent optimization of each sub-problem. However, the obtained results also indicate that this approach is less robust with respect to the initial provided data set, which results in a larger variance of the determined optimum  values.  This issue could be partially avoided by providing larger data sets at the beginning of the process and/or by having the EV threshold vary along the optimization as a function of the constraint models accuracy. However, these solutions would result computationally expensive and/or would require problem specific and time consuming tuning. 

\section{Conclusions \& Perspectives}
In this paper, two alternative extensions of the mixed-variable Bayesian optimization method are proposed in order to perform the optimization of variable-size design space problems. The first approach is a budget allocation strategy based on the independent BO of each sub-problem characterizing the considered problem. The second approach, instead, is based on the definition of a variable-size design space kernel allowing to compute the covariance between samples characterized by partially different sets of variables. Overall, the obtained results show that the proposed algorithms provide a better optimization performance in terms of convergence speed as well as robustness with respect to the initial data set when compared to the independent optimization of each sub-problem. As a consequence, they allow dealing with complex system design problems which depend on technological choices and which are characterized by computationally intensive simulation codes in their most general formulation (\ie without simplifications and/or assumptions).

However, the proposed algorithms are only defined within the context of single-objective optimization problems, whereas real-life design optimization problems are usually characterized by the presence of multiple antagonistic objectives, such as the simultaneous requirement of low production costs and high performance. For this reason, it could be valuable to assess the possibility of adapting and extending existing multi-objective acquisition functions, such as the expected hypervolume improvement, in order to enable the optimization of multi-objective mixed-variable problems and variable-size design space problems. 

Finally, it is also worth mentioning that if the proposed methods were to be applied to actual complex system design problems, the quality of the optimal solution could be improved by coupling the presented global optimization algorithms with local gradient-based optimization algorithms in order to help refining the incumbent solution and thus provide a better final design.


%
%

\bibliographystyle{spphys}       
\bibliography{References}

\section*{Appendix A: Variable-size design space Goldstein function} 

The variable-size design space variant of the Goldstein function which is considered for the testing discussed in Section  \ref{VSDSPRes} is characterized by a global design space with 5 continuous design variables, 4 discrete design variables and 2 dimensional design variables. Depending on the dimensional variable values, 8 different sub-problems can be identified, with total dimensions of 6 or 7, ranging from 2 continuous variables and 4 discrete variables to 5 continuous variables and 2 discrete variables. All of the sub-problem are subject to a variable-size design space constraint.

The resulting optimization problem can be defined as follows:

\begin{align}
 \min & \qquad  f(\textbf{x},\textbf{z},\textbf{w}) \\ 
\text{w.r.t.} &  \qquad \textbf{x} = \{x_1,\dots,x_5\} \ \mbox{ with } \  x_i \in [0,100] \ \mbox{ for }   i = 1,5  \nonumber \\
 &  \qquad \textbf{z} = \{z_1,\dots,z_4\} \ \mbox{ with } \  z_i \in \{ 0,1,2\}  \ \mbox{ for }   i = 1,4 \nonumber \\
 & \qquad \textbf{w}  = \{w_1,w_2\} \ \mbox{ with } \  w_1 \in \{ 0,1,2,3\}  \ \mbox{ and }  w_2 \in \{ 0,1\}  \nonumber \\
\text{s.t.:} &  \qquad g(\textbf{x},\textbf{z},\textbf{w}) \leq 0  \nonumber 
\end{align}

\noindent where:

\begin{equation}
  f(\textbf{x},\textbf{z},\textbf{w}) = \begin{cases}
               f_1(x_1,x_2,z_1,z_2,z_3,z_4) \hfill   \qquad  \mbox{ if }  w_1 = 0  \mbox{ and } w_2 = 0\\
               f_2(x_1,x_2,x_3,z_2,z_3,z_4) \hfill   \qquad  \mbox{ if }  w_1 = 1  \mbox{ and } w_2 = 0\\
               f_3(x_1,x_2,x_4,z_1,z_3,z_4) \hfill   \qquad  \mbox{ if }  w_1 = 2  \mbox{ and } w_2 = 0\\
               f_4(x_1,x_2,x_3,x_4,z_3,z_4) \hfill   \qquad  \mbox{ if }  w_1 = 3  \mbox{ and } w_2 = 0\\
                f_5(x_1,x_2,x_5,z_1,z_2,z_3,z_4) \hfill   \qquad  \mbox{ if }  w_1 = 0  \mbox{ and } w_2 = 1\\
               f_6(x_1,x_2,x_3,x_5,z_2,z_3,z_4) \hfill   \qquad  \mbox{ if }  w_1 = 1  \mbox{ and } w_2 = 1\\
               f_7(x_1,x_2,x_4,x_5,z_1,z_3,z_4) \hfill   \qquad  \mbox{ if }  w_1 = 2  \mbox{ and } w_2 = 1\\
               f_8(x_1,x_2,x_3,x_5,x_4,z_3,z_4) \hfill   \qquad  \mbox{ if }  w_1 = 3  \mbox{ and } w_2 = 1
            \end{cases}
\end{equation}

\noindent and:

\begin{equation}
  g(\textbf{x},\textbf{z},\textbf{w}) = \begin{cases}
               g_1(x_1,x_2,z_1,z_2) \hfill   \qquad  \mbox{ if }  w_1 = 0 \\
               g_2(x_1,x_2,z_2) \hfill   \qquad  \mbox{ if }  w_1 = 1 \\
               g_3(x_1,x_2,z_1) \hfill   \qquad  \mbox{ if }  w_1 = 2 \\
               g_4(x_1,x_2,z_3,z_4) \hfill   \qquad  \mbox{ if }  w_1 = 3 
            \end{cases}
    \end{equation}
        
The objective functions $f_1(\cdot),\dots,f_8(\cdot)$ are defined as follows:

\begin{equation}
\begin{aligned}
    f_1(x_1,x_2,z_1,z_2,z_3,z_4)  = \quad & 53.3108+ 0.184901x_1-5.02914 x_1^3 \cdot 10^{-6} + 7.72522x_1^{z_3} \cdot 10^{-8}-  \\
        & 0.0870775x_2 - 0.106959 x_3 +  7.98772x_3^{z_4} \cdot 10^{-6} +  \\
     & 0.00242482 x_4 +  1.32851 x_4^3 \cdot 10^{-6} - 0.00146393 x_1 x_2 -  \\
       & 0.00301588 x_1 x_3 - 0.00272291 x_1 x_4+ 0.0017004 x_2 x_3 +  \\ 
       & 0.0038428 x_2 x_4 -  0.000198969 x_3x_4 +  1.86025 x_1x_2x_3 \cdot 10^{-5} -  \\
       & 1.88719 x_1 x_2 x_4 \cdot 10^{-6}+  2.50923 x_1 x_3 x_4 \cdot 10^{-5} - \\ & 5.62199 x_2 x_3 x_4 \cdot 10^{-5}
        \end{aligned}
\end{equation}

\noindent where $x_3$ and $x_4$ are defined as a function of $z_1$ and $z_2$ according to the relations defined in Table \ref{Gold1}.
\begin{table}[h!]
\centering
\begin{tabular}{|c|c|c|c|}
\hline
& $z_1 = 0 $ & $ z_1 = 1$  & $  z_1 = 2 $ \\ \hline
$ z_2 = 0$  & $x_3 =  20, x_4 = 20$ & $ x_3 =  50, x_4 = 20$  & $x_3 =  80, x_4 = 20 $ \\ \hline
$ z_2 = 1$  & $x_3 =  20, x_4 = 50$ & $ x_3 =  50, x_4 = 50$  & $x_3 =  80, x_4 = 50 $ \\ \hline
$ z_2 = 2$  & $x_3 =  20, x_4 = 80$ & $ x_3 =  50, x_4 = 80$  & $x_3 =  80, x_4 = 80 $ \\ \hline
\end{tabular}
\caption{Characterization of the variable-dimension search space Goldstein function sub-problem N${}^\circ$1 discrete categories}
\label{Gold1}
\end{table}

\begin{equation}
\begin{aligned}
    f_2(x_1,x_2,x_3,z_2,z_3,z_4)   = \quad & 53.3108+ 0.184901x_1-5.02914 x_1^3 \cdot 10^{-6} + 7.72522x_1^{z_3} \cdot 10^{-8}-  \\
        & 0.0870775x_2 - 0.106959 x_3 +  7.98772x_3^{z_4} \cdot 10^{-6} +  \\
     & 0.00242482 x_4 +  1.32851 x_4^3 \cdot 10^{-6} - 0.00146393 x_1 x_2 -  \\
       & 0.00301588 x_1 x_3 - 0.00272291 x_1 x_4+ 0.0017004 x_2 x_3 +  \\ 
       & 0.0038428 x_2 x_4 -  0.000198969 x_3x_4 +  1.86025 x_1x_2x_3 \cdot 10^{-5} -  \\
       & 1.88719 x_1 x_2 x_4 \cdot 10^{-6}+  2.50923 x_1 x_3 x_4 \cdot 10^{-5} - \\ & 5.62199 x_2 x_3 x_4 \cdot 10^{-5}
        \end{aligned}
\end{equation}

\noindent where $x_4$ is defined as a function of $z_2$ according to the relations defined in Table \ref{Gold2}.
\begin{table}[h!]
\centering
\begin{tabular}{|c|c|c|}
\hline
$z_2 = 0 $ & $ z_2 = 1$  & $  z_2 = 2 $ \\ \hline
$x_4 =  20$ & $ x_4 =  50$  & $x_4 =  80$ \\ \hline
\end{tabular}
\caption{Characterization of the variable-dimension search space Goldstein function sub-problem N${}^\circ$2 discrete categories}
\label{Gold2}
\end{table}

\begin{equation}
\begin{aligned}
    f_3(x_1,x_2,x_4,z_1,z_3,z_4)   = \quad & 53.3108+ 0.184901x_1-5.02914 x_1^3 \cdot 10^{-6} + 7.72522x_1^{z_3} \cdot 10^{-8}-  \\
        & 0.0870775x_2 - 0.106959 x_3 +  7.98772x_3^{z_4} \cdot 10^{-6} +  \\
     & 0.00242482 x_4 +  1.32851 x_4^3 \cdot 10^{-6} - 0.00146393 x_1 x_2 -  \\
       & 0.00301588 x_1 x_3 - 0.00272291 x_1 x_4+ 0.0017004 x_2 x_3 +  \\ 
       & 0.0038428 x_2 x_4 -  0.000198969 x_3x_4 +  1.86025 x_1x_2x_3 \cdot 10^{-5} -  \\
       & 1.88719 x_1 x_2 x_4 \cdot 10^{-6}+  2.50923 x_1 x_3 x_4 \cdot 10^{-5} - \\ & 5.62199 x_2 x_3 x_4 \cdot 10^{-5}
        \end{aligned}
\end{equation}

\noindent where $x_3$ is defined as a function of $z_1$ according to the relations defined in Table \ref{Gold3}.
\begin{table}[h!]
\centering
\begin{tabular}{|c|c|c|}
\hline
$z_1 = 0 $ & $ z_1 = 1$  & $  z_1 = 2 $ \\ \hline
$x_3 =  20$ & $ x_3 =  50$  & $x_3 =  80$ \\ \hline
\end{tabular}
\caption{Characterization of the variable-dimension search space Goldstein function sub-problem N${}^\circ$2 discrete categories}
\label{Gold3}
\end{table}

\begin{equation}
\begin{aligned}
    f_4(x_1,x_2,x_3,x_4,z_3,z_4   = \quad & 53.3108+ 0.184901x_1-5.02914 x_1^3 \cdot 10^{-6} + 7.72522x_1^{z_3} \cdot 10^{-8}-  \\
        & 0.0870775x_2 - 0.106959 x_3 +  7.98772x_3^{z_4} \cdot 10^{-6} +  \\
     & 0.00242482 x_4 +  1.32851 x_4^3 \cdot 10^{-6} - 0.00146393 x_1 x_2 -  \\
       & 0.00301588 x_1 x_3 - 0.00272291 x_1 x_4+ 0.0017004 x_2 x_3 +  \\ 
       & 0.0038428 x_2 x_4 -  0.000198969 x_3x_4 +  1.86025 x_1x_2x_3 \cdot 10^{-5} -  \\
       & 1.88719 x_1 x_2 x_4 \cdot 10^{-6}+  2.50923 x_1 x_3 x_4 \cdot 10^{-5} - \\ & 5.62199 x_2 x_3 x_4 \cdot 10^{-5}
        \end{aligned}
\end{equation}

\begin{equation}
\begin{aligned}
    f_5(x_1,x_2,z_1,z_2,z_3,z_4)  = \quad & 53.3108+ 0.184901x_1-5.02914 x_1^3 \cdot 10^{-6} + 7.72522x_1^{z_3} \cdot 10^{-8}-  \\
        & 0.0870775x_2 - 0.106959 x_3 +  7.98772x_3^{z_4} \cdot 10^{-6} +  \\
     & 0.00242482 x_4 +  1.32851 x_4^3 \cdot 10^{-6} - 0.00146393 x_1 x_2 -  \\
       & 0.00301588 x_1 x_3 - 0.00272291 x_1 x_4+ 0.0017004 x_2 x_3 +  \\ 
       & 0.0038428 x_2 x_4 -  0.000198969 x_3x_4 +  1.86025 x_1x_2x_3 \cdot 10^{-5} -  \\
       & 1.88719 x_1 x_2 x_4 \cdot 10^{-6}+  2.50923 x_1 x_3 x_4 \cdot 10^{-5} - \\ & 5.62199 x_2 x_3 x_4 \cdot 10^{-5} + 5 \cos(2 \pi \frac{x_5}{100})-2
        \end{aligned}
\end{equation}

\noindent where $x_3$ and $x_4$ are defined as a function of $z_1$ and $z_2$ according to the relations defined in Table \ref{Gold5}.
\begin{table}[h!]
\centering
\begin{tabular}{|c|c|c|c|}
\hline
& $z_1 = 0 $ & $ z_1 = 1$  & $  z_1 = 2 $ \\ \hline
$ z_2 = 0$  & $x_3 =  20, x_4 = 20$ & $ x_3 =  50, x_4 = 20$  & $x_3 =  80, x_4 = 20 $ \\ \hline
$ z_2 = 1$  & $x_3 =  20, x_4 = 50$ & $ x_3 =  50, x_4 = 50$  & $x_3 =  80, x_4 = 50 $ \\ \hline
$ z_2 = 2$  & $x_3 =  20, x_4 = 80$ & $ x_3 =  50, x_4 = 80$  & $x_3 =  80, x_4 = 80 $ \\ \hline
\end{tabular}
\caption{Characterization of the variable-dimension search space Goldstein function sub-problem N${}^\circ$5 discrete categories}
\label{Gold5}
\end{table}

\begin{equation}
\begin{aligned}
    f_6(x_1,x_2,x_3,z_2,z_3,z_4)   = \quad & 53.3108+ 0.184901x_1-5.02914 x_1^3 \cdot 10^{-6} + 7.72522x_1^{z_3} \cdot 10^{-8}-  \\
        & 0.0870775x_2 - 0.106959 x_3 +  7.98772x_3^{z_4} \cdot 10^{-6} +  \\
     & 0.00242482 x_4 +  1.32851 x_4^3 \cdot 10^{-6} - 0.00146393 x_1 x_2 -  \\
       & 0.00301588 x_1 x_3 - 0.00272291 x_1 x_4+ 0.0017004 x_2 x_3 +  \\ 
       & 0.0038428 x_2 x_4 -  0.000198969 x_3x_4 +  1.86025 x_1x_2x_3 \cdot 10^{-5} -  \\
       & 1.88719 x_1 x_2 x_4 \cdot 10^{-6}+  2.50923 x_1 x_3 x_4 \cdot 10^{-5} - \\ & 5.62199 x_2 x_3 x_4 \cdot 10^{-5} + 5 \cos(2 \pi \frac{x_5}{100})-2
        \end{aligned}
\end{equation}

\noindent where $x_4$ is defined as a function of $z_2$ according to the relations defined in Table \ref{Gold6}.
\begin{table}[h!]
\centering
\begin{tabular}{|c|c|c|}
\hline
$z_2 = 0 $ & $ z_2 = 1$  & $  z_2 = 2 $ \\ \hline
$x_4 =  20$ & $ x_4 =  50$  & $x_4 =  80$ \\ \hline
\end{tabular}
\caption{Characterization of the variable-dimension search space Goldstein function sub-problem N${}^\circ$6 discrete categories}
\label{Gold6}
\end{table}

\begin{equation}
\begin{aligned}
    f_7(x_1,x_2,x_4,z_1,z_3,z_4)   = \quad & 53.3108+ 0.184901x_1-5.02914 x_1^3 \cdot 10^{-6} + 7.72522x_1^{z_3} \cdot 10^{-8}-  \\
        & 0.0870775x_2 - 0.106959 x_3 +  7.98772x_3^{z_4} \cdot 10^{-6} +  \\
     & 0.00242482 x_4 +  1.32851 x_4^3 \cdot 10^{-6} - 0.00146393 x_1 x_2 -  \\
       & 0.00301588 x_1 x_3 - 0.00272291 x_1 x_4+ 0.0017004 x_2 x_3 +  \\ 
       & 0.0038428 x_2 x_4 -  0.000198969 x_3x_4 +  1.86025 x_1x_2x_3 \cdot 10^{-5} -  \\
       & 1.88719 x_1 x_2 x_4 \cdot 10^{-6}+  2.50923 x_1 x_3 x_4 \cdot 10^{-5} - \\ & 5.62199 x_2 x_3 x_4 \cdot 10^{-5} + 5 \cos(2 \pi \frac{x_5}{100})-2
        \end{aligned}
\end{equation}

\noindent where $x_3$ is defined as a function of $z_1$ according to the relations defined in Table \ref{Gold7}.
\begin{table}[h!]
\centering
\begin{tabular}{|c|c|c|}
\hline
$z_1 = 0 $ & $ z_1 = 1$  & $  z_1 = 2 $ \\ \hline
$x_3 =  20$ & $ x_3 =  50$  & $x_3 =  80$ \\ \hline
\end{tabular}
\caption{Characterization of the variable-dimension search space Goldstein function sub-problem N${}^\circ$7 discrete categories}
\label{Gold7}
\end{table}

\begin{equation}
\begin{aligned}
    f_8(x_1,x_2,x_3,x_4,z_3,z_4   = \quad & 53.3108+ 0.184901x_1-5.02914 x_1^3 \cdot 10^{-6} + 7.72522x_1^{z_3} \cdot 10^{-8}-  \\
        & 0.0870775x_2 - 0.106959 x_3 +  7.98772x_3^{z_4} \cdot 10^{-6} +  \\
     & 0.00242482 x_4 +  1.32851 x_4^3 \cdot 10^{-6} - 0.00146393 x_1 x_2 -  \\
       & 0.00301588 x_1 x_3 - 0.00272291 x_1 x_4+ 0.0017004 x_2 x_3 +  \\ 
       & 0.0038428 x_2 x_4 -  0.000198969 x_3x_4 +  1.86025 x_1x_2x_3 \cdot 10^{-5} -  \\
       & 1.88719 x_1 x_2 x_4 \cdot 10^{-6}+  2.50923 x_1 x_3 x_4 \cdot 10^{-5} - \\ & 5.62199 x_2 x_3 x_4 \cdot 10^{-5} + 5 \cos(2 \pi \frac{x_5}{100})-2
        \end{aligned}
\end{equation}

\vspace{20pt}

The constraints  $g_1(\cdot),\dots,g_4(\cdot)$ are defined are defined as follows:

\begin{equation}
  g_1(x_1,x_2,z_1,z_2)=  - (x_1-50)^2 - (x_2-50)^2 + (20+c_1*c_2)^2
\end{equation}

\noindent where $c_1$ and $c_2$ are defined as a function of $z_1$ and $z_2$ according to the relations defined in Table \ref{GoldConst1}.
\begin{table}[h!]
\centering
\begin{tabular}{|c|c|c|c|}
\hline
& $z_1 = 0 $ & $ z_1 = 1$  & $  z_1 = 2 $ \\ \hline
$ z_2 = 0$  & $c_1 =  3, c_2 = 0.5$ &  $c_1 =  2, c_2 = 0.5$  &  $c_1 =  1, c_2 = 0.5$ \\ \hline
$ z_2 = 1$  & $c_1 =  3, c_2 = -1$ & $c_1 =  2, c_2 = -1$   & $c_1 =  1, c_2 = -1$ \\ \hline
$ z_2 = 1$  & $c_1 =  3, c_2 = -2$ & $c_1 =  2, c_2 = -2$   & $c_1 =  1, c_2 = -2$ \\ \hline
\end{tabular}
\caption{Characterization of the variable-dimension search space Goldstein function constraint}
\label{GoldConst1}
\end{table}

\begin{equation}
  g_2(x_1,x_2,z_2) =  - (x_1-50)^2 - (x_2-50)^2 + (20+c_1*c_2)^2
\end{equation}

\noindent where $c_1 = 0.5$ and $c_2$ is defined as a function of $z_2$ according to the relations defined in Table \ref{GoldConst2}.
\begin{table}[h!]
\centering
\begin{tabular}{|c|c|c|c|}
\hline
$z_2 = 0 $ & $ z_2 = 1$  & $  z_2 = 2 $ \\ \hline
$c_2 = 0.5$ &  $c_2 = -1$  &  $c_2 = -2$ \\ \hline
\end{tabular}
\caption{Characterization of the variable-dimension search space Goldstein function constraint}
\label{GoldConst2}
\end{table}

\begin{equation}
  g_3(x_1,x_2,z_1) =  - (x_1-50)^2 - (x_2-50)^2 + (20+c_1*c_2)^2
\end{equation}

\noindent where $c_2 = 0.7$ and $c_1$ is defined as a function of $z_1$  according to the relations defined in Table \ref{GoldConst3}.
\begin{table}[h!]
\centering
\begin{tabular}{|c|c|c|c|}
\hline
$z_1 = 0 $ & $ z_1 = 1$  & $  z_1 = 2 $ \\ \hline
$c_1 =  3$ &  $c_1 =  2$  &  $c_1 =  1$ \\ \hline
\end{tabular}
\caption{Characterization of the variable-dimension search space Goldstein function constraint}
\label{GoldConst3}
\end{table}

\begin{equation}
  g_4(x_1,x_2,z_3,z_4) =  - (x_1-50)^2 - (x_2-50)^2 + (20+c_1*c_2)^2
\end{equation}

\noindent where $c_1$ and $c_2$ are defined as a function of $z_3$ and $z_4$ according to the relations defined in Table \ref{GoldConst4}.
\begin{table}[h!]
\centering
\begin{tabular}{|c|c|c|c|}
\hline
& $z_3 = 0 $ & $ z_3 = 1$  & $  z_3 = 2 $ \\ \hline
$ z_4 = 0$  & $c_1 =  3, c_2 = 0.5$ &  $c_1 =  2, c_2 = 0.5$  &  $c_1 =  1, c_2 = 0.5$ \\ \hline
$ z_4 = 1$  & $c_1 =  3, c_2 = -1$ & $c_1 =  2, c_2 = -1$   & $c_1 =  1, c_4 = -1$ \\ \hline
$ z_4 = 2$  & $c_1 =  3, c_2 = -2$ & $c_1 =  2, c_2 = -2$   & $c_1 =  1, c_2 = -2$ \\ \hline
\end{tabular}
\caption{Characterization of the variable-dimension search space Goldstein function constraint}
\label{GoldConst4}
\end{table}

\section*{Appendix B: Variable-size design space Rosenbrock function} 

The variable-size design space variant of the Rosenbrock function which is considered for the testing discussed in Section  \ref{VSDSPRes} is characterized by a global design space with 8 continuous design variables, 3 discrete design variables and 2 dimensional design variables. Depending on the dimensional variable values, 4 different sub-problems can be identified, with total dimensions of 6 to 9, ranging from 4 continuous variables and 2 discrete variables to 6 continuous variables and 3 discrete variables. One of the 2 constraints is only active for half of the sub-problems

The resulting optimization problem can be defined as follows:

\begin{equation}
  f(\textbf{x},\textbf{z},\textbf{w}) = \begin{cases}
               f_1(x_1,x_2,x_3,x_4,z_1,z_2) \hfill   \qquad  \mbox{ if }  w_1 = 0  \mbox{ and } w_2 = 0\\
               f_2(x_1,x_2,x_5,x_6,z_1,z_2,z_3) \hfill   \qquad  \mbox{ if }  w_1 = 0  \mbox{ and } w_2 = 1\\
                f_3(x_1,x_2,x_3, x_4,x_7, x_8,z_1,z_2) \hfill   \qquad  \mbox{ if }  w_1 = 1  \mbox{ and } w_2 = 0\\
               f_4(x_1,x_2,x_5,x_6,x_7, x_8,z_1,z_2,z_3) \hfill   \qquad  \mbox{ if }  w_1 = 1  \mbox{ and } w_2 = 1
            \end{cases}
\end{equation}

\noindent and:

\begin{equation}
  g_1(\textbf{x},\textbf{z},\textbf{w}) = \begin{cases}
               g_{1_1}(x_1,x_2,x_3,x_4) \hfill   \qquad  \mbox{ if }  w_1 = 0  \mbox{ and } w_2 = 0\\
                g_{1_2}(x_1,x_2,x_5,x_6) \hfill   \qquad  \mbox{ if }  w_1 = 0  \mbox{ and } w_2 = 1\\
                 g_{1_3}(x_1,x_2,x_3, x_4,x_7, x_8) \hfill   \qquad  \mbox{ if }  w_1 = 1  \mbox{ and } w_2 = 0\\
                g_{1_4}(x_1,x_2,x_5,x_6,x_7, x_8) \hfill   \qquad  \mbox{ if }  w_1 = 1  \mbox{ and } w_2 = 1
            \end{cases}
\end{equation}

\noindent  and:

\begin{equation}
  g_2(\textbf{x},\textbf{z},\textbf{w}) = \begin{cases}
               g_{2_1}(x_1,x_2,x_3,x_4) \hfill   \qquad  \mbox{ if }  w_1 = 0  \mbox{ and } w_2 = 0\\
                 g_{2_2}(x_1,x_2,x_3, x_4,x_7, x_8) \hfill   \qquad  \mbox{ if }  w_1 = 0  \mbox{ and } w_2 = 1\\
                \mbox{Not active} \hfill   \qquad  \mbox{ if }  w_1 = 1  \mbox{ and } w_2 = 0\\
                \mbox{Not active}  \hfill   \qquad  \mbox{ if }  w_1 = 1  \mbox{ and } w_2 = 1
            \end{cases}
\end{equation}

The objective functions $f_1(\cdot),\dots,f_4(\cdot)$ as well as the constraints $g_{1_1}(\cdot),\dots,g_{1_4}(\cdot)$ and $g_{2_1}(\cdot),\dots,g_{2_4}(\cdot)$ are defined as follows:

\begin{equation}
 f_1(x_1,x_2,x_3,x_4,z_1,z_2) =  \begin{cases}
 100z_0 + \sum_i a_1*a_2(x_{i+1} - x_i)^2+(a_1+a_2)/10(1-x_i)^2 \quad \hfill \mbox{ if } z_2 = 0 \\
  100z_0 +  \sum_i 0.7a_1*a_2(x_{i+1} - x_i)^2+(a_1-a_2)/10(1-x_i)^2 \quad \hfill  \mbox{ if } z_2 = 1 
 \end{cases}
\end{equation}
with $a_1 = 7$ and $a_2 = 9$

\begin{equation}
 f_2(x_1,x_2,x_5,x_6,z_1,z_2,z_3) =  \begin{cases}
 100z_0 - 35z_3 + \sum_i a_1*a_2(x_{i+1} - x_i)^2+(a_1+a_2)/10(1-x_i)^2 \quad \hfill \mbox{ if } z_2 = 0 \\
  100z_0 - 35z_3  +  \sum_i 0.7a_1*a_2(x_{i+1} - x_i)^2+(a_1-a_2)/10(1-x_i)^2 \quad \hfill  \mbox{ if } z_2 = 1 
 \end{cases}
\end{equation}
with $a_1 = 7$ and $a_2 = 6$

\begin{equation}
 f_3(x_1,x_2,x_3,x_4,x_7,x_8,z_1,z_2) =  \begin{cases}
 100z_0 + \sum_i a_1*a_2(x_{i+1} - x_i)^2+(a_1+a_2)/10(1-x_i)^2 \quad \hfill \mbox{ if } z_2 = 0 \\
  100z_0 +  \sum_i 0.7a_1*a_2(x_{i+1} - x_i)^2+(a_1-a_2)/10(1-x_i)^2 \quad \hfill  \mbox{ if } z_2 = 1 
 \end{cases}
\end{equation}
with $a_1 = 10$ and $a_2 = 9$

\begin{equation}
 f_2(x_1,x_2,x_5,x_6,x_7,x_8,z_1,z_2,z_3) =  \begin{cases}
 100z_0 - 35z_3 + \sum_i a_1*a_2(x_{i+1} - x_i)^2+(a_1+a_2)/10(1-x_i)^2 \quad \hfill \mbox{ if } z_2 = 0 \\
  100z_0 - 35z_3  +  \sum_i 0.7a_1*a_2(x_{i+1} - x_i)^2+(a_1-a_2)/10(1-x_i)^2 \quad \hfill  \mbox{ if } z_2 = 1 
 \end{cases}
\end{equation}
with $a_1 = 10$ and $a_2 = 6$

\begin{equation}
 g_{1_1}(x_1,x_2,x_3,x_4) =  \sum_i -(x_i-1)^3 +  x_{i+1} - 2.6
\end{equation}

\begin{equation}
 g_{1_2}(x_1,x_2,x_5,x_6) =  \sum_i -(x_i-1)^3 +  x_{i+1} - 2.6
\end{equation}

\begin{equation}
 g_{1_3}(x_1,x_2,x_3,x_4,x_7,x_8) =  \sum_i -(x_i-1)^3 +  x_{i+1} - 2.6
\end{equation}

\begin{equation}
 g_{1_4}(x_1,x_2,x_5,x_6,x_7,x_8) =  \sum_i -(x_i-1)^3 +  x_{i+1} - 2.6
\end{equation}

\begin{equation}
 g_{2_1}(x_1,x_2,x_3,x_4) =  \sum_i -x_i - x_{i+1} + 0.4
\end{equation}

\begin{equation}
 g_{2_2}(x_1,x_2,x_5,x_6) =  \sum_i -x_i - x_{i+1} + 0.4
\end{equation}

\end{document}